\documentclass[12pt,reqno]{amsart}

\usepackage[colorlinks=true,urlcolor=blue,
bookmarks=true,bookmarksopen=true,citecolor=blue
]{hyperref}

\usepackage{color,soul}
\setstcolor{red}

\usepackage{pdfsync}
\usepackage[all, cmtip]{xy}
\usepackage{graphicx}

\usepackage{dsfont}
\usepackage{float}
\usepackage{bm}
\usepackage{mathtools}
\usepackage{amssymb}
\usepackage[T1]{fontenc}

\usepackage{epsfig}
\usepackage{amscd}
\usepackage[mathscr]{eucal}
\usepackage{amssymb}
\usepackage{bm}
\usepackage{amsxtra}
\usepackage{amsmath}
\usepackage{latexsym} 
\usepackage[all]{xy}
\usepackage{enumerate}

\usepackage{enumitem}

\usepackage{color}
\usepackage{xcolor, soul}
\sethlcolor{green}

\usepackage{bbm}

\theoremstyle{plain}

\newtheorem{thm}{Theorem}[subsection]

\newtheorem{cor}[thm]{Corollary}

\newtheorem{lem}[thm]{Lemma}
\newtheorem{prop}[thm]{Proposition}

\theoremstyle{definition}
\newtheorem{dfn}[thm]{Definition}

\newtheorem*{claim-nonum}{Claim}

\theoremstyle{remark}
\newtheorem{remark}[thm]{Remark}
\newtheorem{rem}[thm]{Remark}

\newtheorem{ex}[thm]{Example}

\theoremstyle{plain}

\newtheorem{lemma}[thm]{Lemma}
\def\C{\mathcal{C}}
\def\A{\mathcal{A}}

\def\D{\mathcal{D}}

\DeclareMathOperator{\coker}{coker}

\def\k{{\bf k}}
\def\F{\mathcal F}

\newcommand{\Qed}{\hfill \qedsymbol \medskip}

\newcommand{\id}{\textnormal{id}}

\newcommand{\R}{\mathbb{R}}
\newcommand{\Z}{\mathbb{Z}}

\newcommand{\N}{\mathbb{N}}
\newcommand{\la}{\lambda}
\newcommand{\La}{\Lambda}

\newcommand{\fuk}{\mathcal{F}uk}

\newcommand{\tor}{{\textnormal{Tor}}}

\renewcommand{\k}{\mathbf{k}}

\newcommand{\ac}{\mathscr{A}}
\newcommand{\fil}{\mathscr{F}}
\newcommand{\ch}{\mathbf{Ch}}
\newcommand{\KK}{K}
\newcommand{\fch}{\fil \ch}
\newcommand{\fchfg}{\fil\ch^{\text{fg}}}
\newcommand{\fcht}{\fil\ch^{\text{t}}}
\newcommand{\Ob}{\text{Obj}}
\newcommand{\ho}{\hom}
\newcommand{\ta}{\text{t}}
\newcommand{\fg}{\text{fg}}

\newcommand{\laideal}{\Lambda_P^{(0)}}
\newcommand{\md}{\text{mod}}

\newcommand{\CCF}{{CF'}}
\newcommand{\HHF}{{HF'}}

\newcommand{\pbred}[1]{#1}

\newcommand{\pbnote}[1]{#1}

\newcommand{\jznote}[1]{#1}

\newcommand{\zjnote}[1]{#1}

\newcommand{\zjred}[1]{#1}

\newcommand{\pbaddress}{biran@math.ethz.ch}
\newcommand{\ocaddress}{cornea@dms.umontreal.ca}

\begin{document}

\title{Persistence K-theory}

\date{\today}

\thanks{The second author was supported by an individual NSERC
  Discovery grant. The third author is \zjred{partially supported by National Key R\&D Program of China No.~2023YFA1010500, NSFC No.~12301081, NSFC No.~12361141812, and USTC Research Funds of the Double First-Class Initiative.}}

\author{Paul Biran, Octav Cornea and Jun Zhang}

\address{Paul Biran, Department of Mathematics, ETH-Z\"{u}rich,
  R\"{a}mistrasse 101, 8092 Z\"{u}rich, Switzerland}
\email{\pbaddress}
 
\address{Octav Cornea, Department of Mathematics and Statistics,
  University of Montreal, C.P. 6128 Succ.  Centre-Ville Montreal, QC
  H3C 3J7, Canada} \email{\ocaddress}

\address{Jun Zhang, The Institute of Geometry and Physics, University of Science and Technology of China, 96 Jinzhai Road, Hefei, Anhui, 230026, China} \email{jzhang4518@ustc.edu.cn}

\keywords{Triangulated persistence categories, Grothendieck group,
  Fukaya category, Symplectic manifold, Lagrangian submanifold, Floer
  homology} \subjclass[2020]{55N31 53D12 (Primary); 18F30 19E99 18G80
  53D40 53D37 (Secondary)}

\bibliographystyle{plain}

\bibliographystyle{alphanum}

%

%

\begin{abstract}
 This paper studies the basic $K$-theoretic properties of a triangulated persistence category (TPC). This notion was introduced in \cite{BCZ23} and it 
is a type of category that can be viewed as a refinement of a triangulated category in the  sense that the morphisms sets of a TPC  are persistence modules. We calculate the  $K$-groups in some basic examples and discuss an application to Fukaya categories and to the topology of exact Lagrangian submanifolds.

 \end{abstract}

\maketitle

%
%

\tableofcontents

\section{Introduction}
The main aim of the paper is to study $K$-theoretic properties of
triangulated persistence categories (TPCs) and provide some
applications of these properties in the setting of symplectic
topology.  We recall from \zjnote{\cite{BCZ23}} (or its earlier version \cite{BCZ21}) that a TPC is,
roughly, a persistence category $\C$ - that is a category with
\zjnote{morphisms that} are persistence modules
\zjnote{$\hom_{\C}(A,B)=\{\hom_{\C}^{\alpha}(A,B)\}_{\alpha\in\R}$,
  with} composition compatible with the persistence structure - such
that the $0$-level category $\C_{0}$, with the same objects as $\C$
but with only the $0$-level morphisms
$\hom_{\C_{0}}(-,-)=\hom_{C}^{0}(-,-)$, is triangulated. The full
definition and all relevant notions are recalled in
\S\ref{subsec:tpc}. The fundamental property of such a TPC $\C$ is
that the $\infty$-level category, $\C_{\infty}$, again with the same
objects as $\C$ but with morphisms the $\infty$-limit of
$\hom_{C}(-,-)$ is also triangulated and that its exact triangles each
carry a weight $r\geq 0$. Using this weight one can define a family of
fragmentation pseudo-distances on the objects of $\C$, by taking into
account the weights required to decompose objects by means of exact
triangles.

\

The basic $K$-theoretic constructions associated to a TPC $\C$ appear
in Section \ref{sec-k-TPC}.  Recall that the \zjnote{Grothendieck group
$K(\mathcal{D})$} of a triangulated category $\mathcal{D}$ is defined
as the abelian group freely generated by the objects of the category
modulo the relations $B=A+C$ for each exact triangle
$A\to B\to C\to TA$ in $\mathcal{D}$, where $T$ is the translation
functor of $\mathcal{D}$.

Returning to the triangulated persistence category $\C$, we define its
$K$-group by $$K(\C):=K(\C_{0})~.~$$ This is a linearization of $\C$
that transforms iterated cones involving {\em only \zjnote{$0$-weight}
  triangles} in linear relations. Thus the $K$-class $[X]\in K(\C)$ of
an object $X\in\mathrm{Obj}(\C)$ can be viewed as \zjnote{an algebraic}
invariant of $X$ in the rough sense that it is \zjnote{invariant under the $0$-isomorphisms (see Lemma \ref{prop-tpc})}.  At the same time, because higher weight
triangles are not used in the quotient giving $K(\C)$, this $K$-group
continues to carry relevant metric information.

\

The more precise aim of the paper is to understand more of the
structure of this group, in particular the interplay between the
\zjnote{algebraic} and metric aspects.

\

A first structural result relates $K(\C)$ to the $K$-group of the
$\infty$-level of $\C$, $K_{\infty}(\C):=K(\C_{\infty})$.  We notice
in \S\ref{sec-k-TPC} that there exists a long exact sequence of
$\La_{P}$-modules: \zjnote{ \begin{equation}\label{eq:ex-seq1} 0\to
    {\rm Tor} K(\C)\to K(\ac\C) \to K(\C)\to K_{\infty}(\C)\to 0.
  \end{equation}}
Here $\La_{P}$ is the universal Novikov polynomial
ring \zjnote{with integral coefficients} (we recall its
definition in \S\ref{subsubsec:nov}), and the $K$-groups of TPC's are
naturally $\La_{P}$-modules.  The category $\mathcal{AC}$ is the full
subcategory of $\C$ whose objects are the objects of $\C$ that are
isomorphic to $0$ in $\C_{\infty}$, this is called the acyclic
\zjnote{subcategory} of $\C$. This $\ac\C$ is also a TPC thus the
$K$-group, \pbred{$K(\ac\C) := K(\ac\C_0)$}, is defined \zjnote{as above}.  The term \zjnote{${\rm Tor}K(-)$} is a torsion type functor.
To put this in some perspective, recall from~\cite{BCZ23} that
$\C_{\infty}$ is the Verdier localization  \zjnote{$\C_{0}/\mathcal{AC}_{0}$}. 

\

In Section \ref{sec-fil-cc} we compute the relevant $K$-groups for
categories of filtered chain complexes. We will focus in this
introduction on just one example, $\C^{\fg}=H_{0}(\fchfg)$, the
homotopy category of finitely generated, filtered chain
complexes. Filtered chain complexes (of different types) provide the
most common examples of TPCs.  Additionally, these examples are also
universal in certain ways. For instance, we show that for
$\C=\C^{\fg}$ the exact sequence~\eqref{eq:ex-seq1} can be identified
with the following one:
$$0\to 0\to \La_{P}^{(0)}\to \La_{P}\to \Z\to 0,$$
where $\Lambda_P^{(0)} \subset \Lambda_P$ is the ideal consisting of
the Novikov polynomials $P(t)$ that vanish at $t=1$, i.e.~$P(1)=0$.
The most interesting part of the identification of~\eqref{eq:ex-seq1}
with the latter sequence is the isomorphism
\begin{equation}\label{eq:la1}
  K(\C^{\fg}) \xrightarrow{\; \; \cong \; \;} \La_{P}
\end{equation}
which shows that, in this case, the $\La_{P}$-module $K(\C)$ coincides
with $\La_{P}$ itself. We will denote this isomorphism by $\lambda$
and outline its definition soon. The isomorphism
$K_{\infty}(\C^{\fg}) \xrightarrow{\; \; \cong\; \;} \Z$ is induced by
the Euler characteristic.

The isomorphism
$\la : K(\C^{\fg}) \xrightarrow{\; \; \cong \; \;} \La_{P}$ is induced
by the following way to associate a polynomial to a filtered, finitely
generated chain complex $C$. First associate to $C$ its \zjnote{(graded) persistence
homology $H_*(C)$} (see \S\ref{subsec:pres-mod-pre} for a rapid review of
the definitions). Then associate to \zjnote{$H_*(C)$} its graded barcode
\zjnote{$\mathcal{B}_{H_*(C)}$}.  This consists of a finite collection of
intervals (a.k.a.~bars) of the form $[a,b)$ and $[c,\infty)$, where
each interval has a certain degree \zjnote{in $\Z$}. To an interval $[a,b)$
of degree $k$ we associate the polynomial
$(-1)^{k}(t^{a}-t^{b})\in \La_{P}$, and to an interval $[c,\infty)$ of
degree $k$ we associate $(-1)^{k}t^{c}$. Finally, we sum up the
previously defined polynomials over each interval in the graded
barcode \zjnote{$\mathcal{B}_{H_*(C)}$} and obtain a Novikov polynomial
\zjnote{$\la_{H_*(C)}$}.  The map \zjnote{$C \longmapsto \la_{H_*(C)}\in \La_{P}$} induces
the isomorphism $\la$. An additional point to note here is that
$K(\C)$ is a $\La_{P}$-algebra in case $\C$ has a monoidal structure
(as defined in~\S\ref{subsubsec:ring-K}). This is the case for
$\C^{\fg}$, the monoidal operation being the tensor product, and the
isomorphism $\la$ is a ring isomorphism.

\

In Section \ref{sb:diff-Nov} we investigate a natural way to associate
to a graded persistence module, or to a barcode, a real function. As a
by-product we reinterpret the map $C\to [C]\in K(\C^{\fg})$ that
associates to a filtered chain complex its $K$-class as a filtered
Euler characteristic.  Given a persistence module
$M=\{M^{\alpha}\}_{\alpha\in\R}$ consider the \zjnote{integral-valued} function
$\sigma_{M}(\alpha)=\dim (M^{\alpha})$.  With the definition of
persistence modules considered in this paper (see
\ref{subsec:pres-mod-pre}), and assuming some finite type constraints,
$\sigma_{M}:\R\to \R$ is a bounded step function belonging to a group
of such functions, \zjnote{denoted by $LC_{B}$ and} defined in
\S\ref{subsubsec:iso-rings}. If $M_{\bullet}$ is a graded persistence
module consisting of persistence modules $M_{i}$ for each degree $i$
it is natural to consider the Euler characteristic function
$\bar{\chi}_{M_{\bullet}}\in LC_{B}$,
$\bar{\chi}_{M_{\bullet}}=\sum_{k} (-1)^{k}\sigma_{M_{k}}$. We show
that $LC_{B}$ admits a remarkable ring structure based on convolution
of functions and, with this ring structure, it is isomorphic to the
ring $\La_{P}$. Moreover, under this ring isomorphism
$LC_{B}\cong \La_{P}$, the map
$M_{\bullet}\to \bar{\chi}_{M_{\bullet}}\in LC_{B}$ is another
expression of $\la$ from (\ref{eq:la1}). This chapter also contains - in Proposition \ref{perp-Morse-inequality} - a 
persistence version of the Morse inequalities as well as an application to counting 
bars, in Corollary \ref{cor:count-bars}.

\

In Section \ref{sec-pairing} we return to the general setting and
consider a triangulated persistence category $\C$. Under some finite
type constraints we show that there exists a bilinear pairing
\begin{equation}\label{eq:pair11}
\kappa : \overline{K(\C)}\otimes _{\La_{P}} K(\C)\to \La_{P}
\end{equation}
Here $\overline{K(\C)}$ is \zjnote{the same abelian group as $K(\C)$ but} with the
$\La_{P}$-module structure $\bar{\cdot}$ such that
$t\, \bar{\cdot}\, x = t^{-1}\cdot x$ where $\cdot $ is the module
structure on $ K(\C)$. The pairing $\kappa$ is induced by the map that
associates to a pair of objects $X,Y$ the Novikov polynomial
$\la_{\hom_{C}(X,T^{\ast}(Y))}$ where $\hom_{\C}(X,T^{\ast}(Y))$ is
the graded persistence module having $\hom_{\C}(X,T^{i}Y)$ in degree
$i$ (\zjnote{recall that} $T$ is the translation functor in $\C_{0}$). In view of the
discussion above, notice that $\kappa$ can be viewed as an Euler
pairing.  The pairing is \zjnote{particularly} simple in the case of the
category $\C^{\fg}$.  Using the identification
$K(\C^{\fg})\cong \La_{P}$, we have \zjred{$\kappa(P,Q)(t)=P(t^{-1})Q(t)$} for any
two polynomials $P,Q\in \La_{P}$.  There are conditions - of
Calabi-Yau duality type - implying that the pairing $\kappa$ is skew
symmetric in the sense that
$\kappa(x,y)(t)=\kappa(y,x)(t^{-1})$. These conditions are satisfied
for the category $\C^{\fg}$.

\

Section \ref{S:bar-m-p} is dedicated to analyzing ``measurements'' on
$K$-groups. We see that fragmentation pseudo-distances on the objects
of a TPC $\C$ can be used to define non-trivial group semi-norms on
$K(\C)$. We then discuss other ways to extract quantitative
information from barcodes by using their representation in $LC_{B}$.

\

Section \ref{sec-app} deals with \zjnote{applications to symplectic
  topology.} We discuss some results indicating that the pairing
$\kappa$ defined on the $K$-group of the triangulated persistence
Fukaya category introduced in Chapter 3 of \cite{BCZ23} distinguishes,
in a quantitative way, embedded from immersed Lagrangians. The main
results are in Corollary \ref{thm-emb-char} and Theorem
\ref{thm:dist-immersed}. The corollary shows that the pairing
$\kappa(A,A)$, where $A$ is the $K$-class of an embedded Lagrangian,
is a constant Novikov polynomial equal to the Euler characteristic of
the Lagrangian (up to sign). The theorem provides \zjnote{a} way to
estimate the minimal ``energy'' required to transform certain immersed
Lagrangians into embedded ones. The estimate is somewhat subtle as it
involves also estimates on the number of bars in certain associated
\zjnote{barcodes} that are deduced from the persistence Morse
inequalities in \S\ref{sb:diff-Nov}.

\subsubsection*{Acknowledgments} \pbnote{We would like to thanks the
  anonymous referee for a careful reading of the paper and for
  numerous comments that helped to improve the exposition.}

\section{Prerequisites}
We list here some of the conventions used further in the paper and we
recall some of the notions required to make the reading of this paper
more accessible.
\subsection{Novikov rings}\label{subsubsec:nov}
We will use the following versions of the Novikov ring. We denote by
\begin{equation} \label{eq:Nov} \Lambda = \left\{ \sum_{k=0}^{\infty} n_k
  t^{a_k}\, \bigg| \, n_k \in \mathbb{Z}, a_k \in \mathbb{R}, \lim_{k \to
    \infty} a_k = \infty\right\}
\end{equation}
the Novikov ring with coefficients in $\mathbb{Z}$ (and formal
variable $t$), endowed with the standard multiplication. Denote also
by $\Lambda_{P} \subset \Lambda$ the subring of Novikov polynomials,
whose elements consist of {\em finite} sums
$\sum n_k t^{a_k} \in \Lambda$.

\subsection{Persistence modules} \label{subsec:pres-mod-pre} The
bibliography on persistence modules contains a large number of
variants of the definition. We fix here the setting that we will use
in this paper.

A persistence module $M$ is a family of vector spaces over the field
$\k$, indexed by the reals, $M= \{M^{s}\}_{s\in \R}$, together with
maps $i_{s,t}:M^{s}\to M^{t}$ for all $s \leq t$ satisfying
$i_{t,r}\circ i_{s,t}=i_{s,r}$ \zjnote{for all
  $s \leq t \leq r \in \R$} and \zjnote{$i_{s, s} = \mathds{1}$}
\zjnote{for any $s \in \R$}. We will always assume that the following
three conditions apply to the persistence modules $M$ considered here.
\begin{itemize}
\item[(i)](lower semi-continuity) For any $s\in \R$ and any $t\geq s$
  sufficiently close to $s$, the map $i_{s,t}:M^{s}\to M^{t}$ is an
  isomorphism.
\item[(ii)](lower bounded) \zjnote{For} $s$ sufficiently small we have
  $M^{s}=0$.
\item[(iii)] \zjnote{(tame) For every $s \in \mathbb{R}$,\, we have} 
\begin{equation} \label{eq:finite-dim}
\dim_{\k} (M^s) < \infty.
\end{equation}
\end{itemize}

We will often use graded persistence modules. Namely, $M_{\bullet}$ is
a graded persistence module if
$ M_{\bullet}= \{M_i\}_{i \in \mathbb{Z}}$ is a sequence of
persistence modules indexed by $\mathbb{Z}$ with the index indicating
the degree. The tameness condition in this case becomes:
\begin{equation}\label{eq:finite-dim-graded}
  \mathrm{for\ every}\  s \in \mathbb{R},
  \  \dim_{\k}(\oplus_{i \in \mathbb{Z}} M_i^s) < \infty~.~
\end{equation}

The homology of filtered chain complexes provides some of the most
typical examples of persistence modules.  For our purposes here, a
filtered chain complex is a triple $(C_*, \partial, \ell)$ where
$(C_*, \partial)$ is a chain complex and
$\ell: \oplus_i C_i \to \R \cup \{-\infty\}$ is a filtration function
(see~\cite[Definition~2.2]{UZ16}). \zjnote{In particular,} $\partial$ and $\ell$ are
compatible in the sense that
\[ \ell(\partial x) \leq \ell(x) \,\,\,\,\,\mbox{for any $x \in
    C_*$}. \] \zjnote{It follows that} for any $r \in \R$ the truncation
\zjnote{$(C_*^{\leq r}, \partial)$} where
$C_*^{\leq r} =\{x\in C \ | \ \ell(x)\leq r\}$ is a well-defined chain
complex. The same holds also for the ``strict'' truncations
$C_*^{< r}$.

We will consider mainly two types of filtered chain complexes:
\label{p:types-ch}
\begin{itemize}
\item[(i)] {\em finitely generated} - chain complexes for which the
  total space $\oplus_{k\in \Z}C_{k}$ is finite dimensional, and
  moreover $C_{\ast}^{\leq r}=0$ for $r$ sufficiently small.
\item[(ii)] {\em tame} - chain complexes such that
  $\oplus_{k\in \Z}C^{\leq r}_{k}$ is finite dimensional for all
  $r\in \R$, and moreover $C_{\ast}^{\leq r}=0$ for $r$ sufficiently
  small.
\end{itemize}
If a filtered chain complex $C$ is tame, then
$\{C_{\ast}^{\leq r}\}_{r \in \mathbb{R}}$ forms a \zjnote{graded} persistence module
in the class discussed above. Moreover, the homology of these
truncations, $\{H_{\ast}(C^{\leq r})\}_{r\in\R}$, also forms a
\zjnote{graded} persistence module in the same class. We will denote this persistence
module by \zjnote{$H_*(C)$}.

\subsection{Barcodes}\label{subsec:barfuc}
In this paper we will use the following version of barcodes. A barcode
$\mathcal{B} = \{(I_j, m_j)\}_{j \in \mathcal{J}}$ is a collection of
pairs consisting of intervals $I_j \subset \mathbb{R}$ and positive
integers $m_j \in \mathbb{Z}_{>0}$, indexed by a set $\mathcal{J}$,
and satisfying the following {\em admissibility} conditions:
\begin{enumerate}
\item[(i)] $\mathcal{J}$ is assumed to be either finite or
  $\mathcal{J} = \mathbb{Z}_{\geq 0}$.
\item[(ii)] Each interval $I_j$ is of the type $I_j = [a_j, b_j)$, with
  $-\infty < a_j < b_j \leq \infty$.
\item[(iii)] In case $\mathcal{J} = \mathbb{Z}_{\geq 0}$ we assume that
  \zjnote{$a_j \rightarrow \infty$ as $j \rightarrow \infty$}.
\end{enumerate}
The intervals $I_j$ are called bars and for each $j$, $m_j$ is called
the multiplicity of the bar $I_j$.

To every barcode $\mathcal{B}$ we can associate a persistence module
in the following way. For an interval $I = [a,b)$ as above, denote by
$V(I)$ the persistence module defined by
\begin{equation} \label{eq:VI} V^s(I) = 
  \begin{cases}
    \k, & \text{if $s \in I$;} \\
    0, & \text{if $s \notin I$,}
  \end{cases}
\end{equation}
and with persistence maps $i_{s,t}:V^s(I) \longrightarrow V^t(I)$
being the identity map $\id_{\k}$ whenever $s,t \in I$. We call $V(I)$
the interval module corresponding to $I$.

By taking direct sums we associate to every barcode
$\mathcal{B} = \{(I_j, n_j)\}_{j \in \mathcal{J}}$ a persistence
module, as follows:
\begin{equation} \label{eq:mod-B} V(\mathcal{B}) := \bigoplus_{j \in
    \mathcal{J}} V(I_j)^{\oplus m_j}.
\end{equation}
It is easy to see that $V(\mathcal{B})$ satisfies the conditions in \S
\ref{subsec:pres-mod-pre}.  Conversely, it is well known that this
construction can be reversed.  More precisely, with the conventions
from \S \ref{subsec:pres-mod-pre}, we have that for every persistence
module $M$ there exists a unique barcode $\mathcal{B}$ as above and
such that $M \cong V(\mathcal{B})$. In other words, barcodes classify
isomorphism types of persistence modules. See~\cite{UZ16,PRSZ20} for more
details. We will denote the barcode corresponding to a persistence
module $M$ by $\mathcal{B}_M$.

\

We will need also the notion of graded barcodes. By this we mean a
sequence of barcodes
$\mathcal{B} = \{\mathcal{B}_i\}_{i \in \mathbb{Z}}$ indexed by the
integers and satisfying the following additional assumption. Denote by
$\mathcal{B}^{\text{tot}}$ the union of all the intervals in all the
$\mathcal{B}_i$'s, together with their multiplicities and adding up
multiplicities whenever bars from different $\mathcal{B}_i$'s
coincide). We assume that $\mathcal{B}^{\text{tot}}$ is an admissible
barcode, in the sense that all the multiplicities appearing in it are
finite, and moreover that it satisfies the admissibility conditions
(1), (2), (3) above.

We will often regard a single (ungraded) barcode $\mathcal{B}$ as a
graded one, simply by placing it in degree $0$.

Let $\mathcal{B} = \{(I_j, m_j)\}_{j \in \mathcal{J}}$ be an ungraded
barcode. We associate to $\mathcal{B}$ an element
$\lambda_{\mathcal{B}} \in \Lambda$ of the Novikov ring as follows:
\begin{equation} \label{eq:lambda-B} \lambda_{\mathcal{B}} := \sum_{j
    \in \mathcal{J}} m_j \lambda_{I_j},
\end{equation}
where
\begin{equation} \label{eq:lambda-I}
  \lambda_{I_j} := 
  \begin{cases}
    t^{a_j} - t^{b_j}, & \text{if $I_j = [a_j, b_j)$;} \\
    t^{a_j},          & \text{if $I_j = [a_j, \infty)$.}
  \end{cases}
\end{equation}
Let $M$ be a persistence module with barcode $\mathcal{B}_{M}$.
Applying the construction above, we obtain an element of the Novikov
ring
\begin{equation} \label{eq:lambda-M}
  \lambda_{M} := \lambda_{\mathcal{B}_M} \in \Lambda.
\end{equation}

The definition of $\la_{M}$ can be extended to the graded case.
Namely, let $M_{\bullet} = \{M_i\}_{i \in \mathbb{Z}}$ be a graded
persistence module. \zjnote{Assume further the tameness condition in
  (\ref{eq:finite-dim-graded}), then define:}
\begin{equation} \label{eq:lambda-M-blt}
  \lambda_{M_{\bullet}} := \sum_{i \in \mathbb{Z}} (-1)^i
  \lambda_{M_i}  \ \in \ \La ~.~
\end{equation}
A similar definition is available for graded barcodes.

\subsection{Triangulated persistence categories}\label{subsec:tpc}

Triangulated persistence categories were introduced \zjnote{in
  \cite{BCZ21,BCZ23}}.  We will only roughly review here some of the basics to
make the current paper more accessible.

A TPC is a category $\C$ with four properties \zjnote{summarized below
  (we refer to \cite{BCZ23} for the full definitions).}
\begin{itemize}
\item[(i)] $\C$ is a persistence category. This means that for each
  two objects $A,B$ of $\C$, the morphisms $\hom_{\C}(A,B)$ form a
  persistence module $\{\hom_{\C}^{s}(A,B)\}_{s \in \mathbb{R}}$ and
  that composition of morphisms is compatible with the persistence
  structure in the obvious way. In the current paper we will assume
  that these persistence modules are always lower-bounded and
  lower-semicontinuous - see \S \ref{subsec:pres-mod-pre}.
\item[(ii)] $\C$ carries a ``shift functor'' $\Sigma$. This is a
  structure composed of shift functors $\Sigma^{r}:\C \to \C$
  \zjnote{where $\Sigma^r \circ \Sigma^s = \Sigma^{r+s}$ for all
    $r,s \in \R$ and} \zjnote{$\Sigma^{0}=\mathds{1}$} that are
  compatible with the persistence structure, and of natural
  isomorphisms $\eta_{r,s}:\Sigma^{r}\to \Sigma^{s}$ satisfying some
  natural compatibility relations and with $\eta_{r,r}=id$,
  $\eta_{r,s}\in \hom^{s-r}(-,-)$.

  To keep an example in mind, the simplest TPC is the (homotopy)
  category of filtered chain complexes.  In that case $\Sigma^{r}C$ is
  obtained from the filtered chain complex $C$ by shifting the
  filtration up by $r$ and the natural transformation
  $\Sigma^{r}C \stackrel{\eta_{r,s}}{\longrightarrow} \Sigma^{s}C$ is
  the identity at the vector space level, when the filtration is
  forgotten.

\item[(iii)] The category $\C_{0}$ with the same objects as $\C$ and
  whose \zjred{morphisms} are the $0$-level morphisms of $\C$,
  $\hom_{\C}^{0}(- , -)$, is triangulated, $\Sigma^{r}$ is a
  triangulated functor (on $\C_{0}$), and $\eta_{r,s}$ are compatible
  with the additive structure on $\C_{0}$. Note that being
  triangulated, $\mathcal{C}_0$ has a translation functor $T$. This
  functor should not be confused with the shift functor $\Sigma$ which
  has to do with shifting the persistence level. (Unfortunately the
  names ``shift'' and ``translation'' functor are both used for the
  same thing in the context of triangulated categories, often related
  to a shift or translation in grading. However in our framework these
  two notions are absolutely different from each other.)
\item[(iv)]  To formulate the last property that characterizes a TPC we
  need a couple more notions that will also be useful later in the
  paper.  We say that $f,g\in \hom_{\C}^{s}(A,B)$ are $r$-{\em
    equivalent} if $i_{s,s+r}(f-g)=0$. Here
  $i_{s, s+r} :\hom^{s}(-,-)\to \hom^{s+r}(-,-)$ are the \zjnote{persistence}
  structural maps. We write $f\simeq_{r} g$.  An object $K$ is called
  $r$-{\em acyclic} for some $r\geq 0$ if $\mathds{1}_{K}\simeq_{r} 0$.

 \zjnote{For all $r\geq 0$, t}he natural transformation
  $\eta_{r,0}:\Sigma^{r} \to \Sigma^0 = \mathds{1}$, specialized to
  $X \in \Ob(\mathcal{C})$ belongs to $\hom^{-r}(\Sigma^r X, X)$ so we
  can apply to it
  $i_{-r,0}:\hom^{-r}(\Sigma^{r}X,X)\to \hom^{0}(\Sigma^{r}X, X)$.
  The result is denoted by $\eta_{r}=i_{-r,0}(\eta_{r,0})$ and thus
  $\eta_{r}\in \hom_{\C_{0}}$.

 \zjnote{Then t}he fourth property that defines a TPC is the following: For each
  $r\geq 0$ and object $X$ of $\C$ there is an exact triangle in
  $\C_{0}$
  $$\Sigma^{r}X\stackrel{\eta_{r}}{\longrightarrow}
  X\to K\to T\Sigma^{r} X$$ with $K$ being $r$-acyclic.
\end{itemize}

\begin{dfn} \label{def:tpc-tame} (a) A triangulated persistence
  category $\C$ is called {\em tame} if for any two objects
  $A,B \in \mathrm{Obj}(\C)$ the persistence module $\hom_{\C}(A,B)$
  satisfies condition~\eqref{eq:finite-dim} and is called of {\em
    finite type} if the barcode of the persistence module
  $\ho_{\C}(A,B)$ has only finitely many bars.
 
 (b) The triangulated persistence category $\C$ is called {\em
    bounded} if for every two objects
  $X, Y \in \mathrm{Obj}(\mathcal{C})$ there exists $l \in \mathbb{N}$
  with the property that $\ho_{\mathcal{C}}(X, T^iY) = 0$ for every
  $|i| \geq l$. Here, $T$ is the translation functor of the
  triangulated category $\mathcal{C}_0$.
\end{dfn}

Obviously, if a TPC is tame and $X,Y\in \mathrm{Obj}(\C)$ then
$\la_{\ho_{\C}(X,Y)}\in \La$ is well-defined. Further, given
$X,Y\in \mathrm{Obj}(\C)$, consider the graded persistence module
$\ho(X,T^{\bullet}Y)$ defined by:
\begin{equation}\label{eq:pers-grad-mor}
  \ho(X,T^{\bullet}Y)_{i} = \ho_{\mathcal{C}}(X,T^iY), \; i \in
  \mathbb{Z}~.~
\end{equation}
If $\C$ is tame and bounded, then the element
\begin{equation}\label{eq:mor-pers-la}
  \la_{\ho_{\C}(X,T^{\bullet}Y)}\in \La
\end{equation}
is well defined. This element belongs to $\La_{P}$ in case $\C$ is
of finite type and bounded.

\
  
A useful notion in a TPC is that of an $r$-{\em isomorphism}. This is
a map $f\in \hom^{0}_{\C}(A,B)$ that fits into an exact triangle
$A\stackrel{f}{\longrightarrow} B\longrightarrow K\longrightarrow TA$
in $\C_{0}$ with $K$ $r$-acyclic.  An $r$-isomorphism $f$ admits a
right $r$-inverse, which is a map
$g\in \hom^{0}_{\C}(\Sigma ^{r} B, A)$ with $f\circ g=\eta_{r}$. There
are also left $r$-inverses with a similar definition. But left and
right inverses do not generally agree.
  
\
  
The key structure that is available in a TPC is that of {\em strict
  exact triangle} of weight $r\geq 0$.  This is a pair
$\widetilde{\Delta}=(\Delta, r)$, where $r \in [0,\infty)$ and
$\Delta$ is a triangle in $\C_{0}$,
$$\Delta \ \ : \ \ A\stackrel{\bar{u}}{\longrightarrow}
B\stackrel{\bar{v}}{\longrightarrow}
C\stackrel{\bar{w}}{\longrightarrow}\ \Sigma^{-r}TA,$$ which can be
completed to a commutative diagram in $\C_{0}$
\begin{equation}\label{dfn-set-2}
  \xymatrix{
    & & \Sigma^r C \ar[d]_-{\psi} \ar[rd]^-{\Sigma^{r} \bar{w}} & \\
    A \ar[r]^-{u} & B \ar[r]^-{v} \ar[rd]_-{\bar{v}}
    & C' \ar[r]^-{w} \ar[d]^-{\phi} & TA\\
    & & C &} 
\end{equation}
with $u = \bar{u}$ and such that the following holds: the triangle
$A \stackrel{u}{\longrightarrow} B\stackrel{v}{\longrightarrow}
C'\stackrel{w}{\longrightarrow} TA$ is exact in $\C_{0}$ and $\phi$ is
an $r$-isomorphism with an $r$-right inverse $\psi$.  The {\em weight}
of the triangle $\Delta$ is denoted by $w(\Delta)= r$.

\

The category $\C_{\infty}$ associated to $\C$ has the same objects as
$\C$ but its morphisms are: \label{p:C-infty}
$$\ho_{\C_{\infty}}(A,B)=\lim_{\longrightarrow} \ho_{\C}^{\alpha}(A,B),$$
where the direct limit (or colimit) is taken along
$\alpha \in \mathbb{R}$ with respect to the standard order on
$\mathbb{R}$ (see \zjnote{\cite[Definition~2.9 (ii)]{BCZ23}}).  

The set of all $r$-acylic objects in $\C$, $r\geq 0$, are the objects
of a \zjnote{full subcategory} $\mathcal{A}\C$ of $\C$. It is easy to
see that $\mathcal{A}\C$ is itself a TPC and that its $0$-level,
$\mathcal{A}\C_{0}$, is simply the \zjnote{full subcategory} of
$\C_{0}$ with the same objects as $\mathcal{A}\C$.  It is shown in
\zjnote{\cite[Proposition 2.38]{BCZ23}} that $\C_{\infty}$ is the
Verdier localization of $\C_{0}$ with respect to
$\mathcal{A}\C_{0}$. This implies \zjred{that} $\C_{\infty}$ is
triangulated.
  
The exact triangles of $\C_{\infty}$ are related to \zjred{the} strict
exact triangles of $\C$ as follows. First, given any diagram in $\C$
$$\Delta :  A\stackrel{u}{\longrightarrow}
B\stackrel{v}{\longrightarrow} C \stackrel{w}{\longrightarrow} D$$ we
denote by $\Sigma^{s_{1},s_{2},s_{3},s_{4}}\Delta$ the diagram
$$\Sigma^{s_{1},s_{2},s_{3},s_{4}}\Delta :
\Sigma^{s_{1}} A\stackrel{u'}{\longrightarrow} \Sigma ^{s_{2}}
B\stackrel{v'}{\longrightarrow}\Sigma^{s_{3}} C
\stackrel{w'}{\longrightarrow} \Sigma^{s_{4}}D$$ where the maps
$u', v',w'$ are obtained from $u,v,w$ by composing with the
appropriate $\eta_{r,s}$'s.

\

With this convention, a triangle
$\Delta: A \xrightarrow{u} B \xrightarrow{v} C \xrightarrow{w} TA$ in
$\C_{\infty}$ is called {\em exact} if there exists a \pbnote{diagram}
$\bar{\Delta}: A \xrightarrow{\bar{u}} B \xrightarrow{\bar{v}} C
\xrightarrow{\bar{w}} TA$ in $\C$ that represents $\Delta$ (in the
sense that $u$ is the $\infty$-image of $\bar{u}$ and similarly for
the other morphisms), such that the shifted triangle
\begin{equation} \label{extri-shifted} \widetilde{\Delta} = \Sigma^{0,
    -s_{1}, -s_{2}, -s_{3}} \bar{\Delta}
\end{equation}
is strict exact (of weight $s_{3}$) in $\C$ and $0\leq s_{1}\leq s_{2}\leq s_{3}$. The
{\em unstable weight} of $\Delta$, $w_{\infty}(\Delta)$, is the
infimum of the strict weights $w(\widetilde{\Delta})$ of all the
strict exact triangles $\widetilde{\Delta}$ as above. Finally, the
{\em weight} of $\Delta$, $\bar{w}(\Delta)$, is given by:
$$\overline{w}(\Delta)=\inf\{ w_{\infty}(\Sigma^{s,0,0,s}\Delta)
\ | \ s\geq 0\}~.~$$
If there does not exist $\widetilde{\Delta}$ as in the definition,
then we put $w_{\infty}(\Delta) =\infty$.

The two notions of weight that are of most interest are \zjnote{$w$ for
strict exact triangles in $\C$} and $\overline{w}_{\C}=\overline{w}$
for exact triangles in $\C_{\infty}$.
  
\subsection{Triangular weights and fragmentation
  pseudo-metrics}\label{subsec:tr-weight}
We discuss here the definition of triangular weights and their
application to the definition of a family of pseudo-metrics on the
objects of a triangulated category. Our main example of interest is
the persistence weight $\overline{w}_{\C}$ that was recalled
in~\S\ref{subsec:tpc}.
 
Let $\mathcal{D}$ be a (small) triangulated category and denote by
$\mathcal{T}_{\mathcal{D}}$ its class of exact triangles. A {\em
  triangular weight} $w$ on $\mathcal{D}$ is a function
$$w:\mathcal{T}_{\mathcal{D}}\to [0,\infty)$$
with two properties, one is a weighted form of the octahedral axiom
(we refer to \zjnote{\cite[Definition 2.1 (i)]{BCZ23}} for the precise
formulation) and the other is a normalization axiom requiring that all
exact triangles have weight at least equal to some
$w_{0}\in [0,\infty)$ and this $w_{0}$ is attained for all exact
triangles of the form $0\to X\xrightarrow{\mathds{1}_{X}}X\to 0$,
$X\in\mathrm{Obj}(\mathcal{D})$, and their rotations.
  
An additional property one may require from a triangular weight is
subadditivity. This means that taking direct sums with triangles of
the type $0\to X\xrightarrow{\mathds{1}_{X}}X\to 0$ does not increase
the weight. See~\cite[Chapter~2, Section~2.1]{BCZ23} for the precise
definitions.




The weight $\bar{w}_{\C}$ recalled in~\S\ref{subsec:tpc} is a
triangular weight on $\C_{\infty}$ called the {\em persistence weight}
(induced from $\C$).  It is subadditive with $(\bar{w}_{\C})_{0}=0$.

\

Fix a triangulated category $\mathcal{D}$ together with a triangular
weight $w$ on $\mathcal{D}$.  Let $X$ be an object of
$\mathcal{D}$. An {\em iterated cone decomposition} $D$ of $X$ with
{\em linearization} $\ell(D) = (X_1,X_{2}, ..., X_n)$ consists of a
family of exact triangles in $\mathcal{D}$:
\begin{equation}\label{eq:iterated-dec}
  \left\{
    \begin{array}{ll}
      \Delta_{1}: \, \, & X_{1}\to 0\to Y_{1}\to  TX_{1}\\
      \Delta_2: \,\, & X_2 \to Y_1 \to Y_2 \to  TX_2\\
      \Delta_3: \,\, & X_3 \to Y_2 \to Y_3 \to  TX_3\\
                        &\,\,\,\,\,\,\,\,\,\,\,\,\,\,\,\,\,\,\,\,\vdots\\
      \Delta_n: \,\, & X_n \to Y_{n-1} \to Y_{n} \to  TX_n
    \end{array}
  \right.
\end{equation}
with $X=Y_{n}$ (to include the case $n=1$ we set $Y_0=0$).  The weight
of such a cone decomposition is defined by:
\begin{equation}\label{eq:weight-cone} 
w(D)=\sum_{i=1}^{n}w(\Delta_{i})-w_{0}~.~
\end{equation}
\

Let $\mathcal{F} \subset {\rm Obj}(\mathcal{D})$.  For two objects
$X, X'$ of $\mathcal{D}$, define
\begin{equation} \label{frag-met-0} \delta^{\mathcal F} (X, X') =
  \inf\left\{ w(D) \, \Bigg| \,
    \begin{array}{ll} \mbox{$D$ is an
      iterated cone decomposition} \\
      \mbox{of $X$ with linearization}
      \ \mbox{$(F_1, ..., T^{-1}X', ..., F_k)$}\\
      \mbox { where $F_i \in \mathcal
      F$, $k \geq 0$}
    \end{array}
  \right\}
\end{equation}

and symmetrize:
\begin{equation}\label{eq:frag-metr2}
  d^{\F}(X, X') = \max\{\delta^{\F}(X, X'), \delta^{\F}(X', X)\}. 
\end{equation}
 
    
The weighted octahedral axiom implies that $\delta^{F}$ satisfies the
triangle inequality, hence $d^{\mathcal{F}}$ is a pseudo-metric,
called {\rm the fragmentation pseudo-metric} associated to $w$ and
$\mathcal{F}$. See~\cite{BCZ23}.


For a TPC $\C$, as defined in~\S\ref{subsec:tpc}, we will be
interested in the fragmentation pseudo-metrics on
$\mathrm{Obj}(\C)=\mathrm{Obj}(\C_{\infty})$ that are associated to
the persistence weight $\bar{w}_{\C}$ defined on $\C_{\infty}$.

In addition, there are some simpler decompositions that are often
useful, which we call {\em reduced decompositions}. Let
$X \in \Ob(\mathcal{C})$. A reduced cone-decomposition of $X$ in
$\mathcal{C}$ of weight $r$ is a pair $D = (\bar{D},\phi)$ where
$\bar{D}$ is a cone-decomposition in $\mathcal{C}_{\infty}$ as in
(\ref{eq:iterated-dec}), with all triangles $\Delta_{i}$ of weight
$0$, and with
$$\phi : Y_{n}\to X$$  an $r$-isomorphism in $\C$.
The linearization of such a decomposition is by definition the
linearization of $\bar{D}$. We denote by $w(D) := r$ the weight of
$D$.

It is easy to see that $\Delta_{n}$ combined with the
$r$-isomorphism $\phi$ can be used to define an exact triangle in
$\mathcal{C}_{\infty}$ of the form:
$$X_{n}\to Y_{n-1} \to X\to TX_{n}$$ of weight $r$.   
Thus, if $X$ admits a reduced cone-decomposition of weight $r$ in
$\C$, then it also admits a cone-decomposition $D'$ in
$\mathcal{C}_{\infty}$ with $\bar{w}_{\C}(D')=r$.

Reduced decompositions are much easier to handle in applications
compared to general ones. Moreover, it is shown in the proof of
\zjnote{Corollary 2.86 in \cite{BCZ23}} that given any decomposition
$D$ of $X$ as in (\ref{eq:iterated-dec}) there exists a reduced
decomposition $D'$ of $X$, of weight no more than $4w(D)$, and with
linearization
$\ell(D')=(\Sigma^{r_1}X_{1},\Sigma^{r_{2}}X_{2},\ldots,
\Sigma^{r_{n}}X_{n})$ for some $r_{i}\in \R$.

Note that for reduced decomposition the analogous measurement
to~\eqref{frag-met-0} does not satisfy the triangle inequality, hence
such decompositions do not seem to lead to pseudo-metrics. However, as
we will see in~\S\ref{sb:k-lag} reduced decompositions turn out to be
useful for estimates involving $K$-theoretic invariants, weights and
barcodes.


\subsection{Fukaya triangulated persistence categories}\label{subsec:Fuk-TPC}
The notion of triangulated persistence category is inspired, in part,
by natural filtrations in Floer theory. Indeed, it is shown in
\zjnote{\cite[Chapter 3]{BCZ23}} that the usual derived Fukaya
category generated by exact, closed Lagrangians admits a TPC
refinement. In this subsection we review enough of the constructions
from \cite{BCZ23} as is necessary for the applications in
\S\ref{sec-app}. We refer to this \pbred{work} for \zjnote{more
  details}.

\subsubsection{Filtered $A_{\infty}$  categories}\label{subsubsec:ainfty-alg}
\pbred{Before we get to the filtered framework, let us briefly
  describe our general conventions regarding (unfiltered)
  $A_{\infty}$-categories.}

\pbred{For coherence with the rest of the algebraic parts of the
  paper, as well as with some of our previous work
  (e.g.~\cite{Bi-Co-Sh:LagrSh, Bi-Co:fuk-cob, Bi-Co:cob1} and, in
  part, in ~\cite{BCZ23}) we use here homological conventions for
  $A_{\infty}$-categories.}  \pbred{More specifically, this means that
  in our $A_{\infty}$-categories $\mathcal{A}$, the $\mu_d$-operation
  has degree $d-2$ (rather than $2-d$), for all $d \geq 1$. In
  particular, the hom's $\hom_{\mathcal{A}}(X,Y)$ between any two
  objects $X,Y$, endowed with $\mu_1$, are {\em chain} complexes
  (rather than cochain complexes). Moreover, whenever relevant, we
  assume the homology units (or sometimes even strict units) $e_X$ to
  be in degree $0$, $e_X \in \hom_{\mathcal{A}}(X,X)_0$. Turning to
  the {\em homological} category $H(\mathcal{A})$ of $\mathcal{A}$,
  its hom's are the {\em homologies} of the chain complexes
  $(\hom_{\mathcal{A}}(-,-), \mu_1)$ and the $\mu_2$-operation induces
  products
  $$H_i(\hom_{\mathcal{A}}(X,Y)) \otimes H_j(\hom_{\mathcal{A}}(Y,Z))
  \longrightarrow H_{i+j}(\hom_{\mathcal{A}}(X,Z)),$$ hence we obtain
  a graded category (which is unital in case $\mathcal{A}$ is
  homologically unital).} \pbred{Sometimes we will restrict to the
  degree $0$ homological category $H_0(\mathcal{A})$. In case
  $\mathcal{A}$ is homology unital, $H_0(\mathcal{A})$ becomes a
  unital category since we have assumed our homology units $e_X$ to be
  in degree $0$.}  \pbred{Later on, in~\S\ref{subsubsec:filtr-fuk},
  when dealing with Fukaya categories, we will explain how to adjust
  the standard grading setting in Floer theory to fit the conventions
  above.}

\pbred{We now turn to the filtered case. A filtered
  $A_{\infty}$-category $\mathcal{A}$ is an $A_{\infty}$-category
  (with the conventions described above)} over a given base field
$\k$, such that the $\hom$ spaces $\hom_{\mathcal{A}}(X,Y)$ between
every two objects $X,Y$ are filtered (with increasing filtrations) and
{\em all} the composition maps $\mu_d$, $d \geq 1$, respect the
filtrations. We endow $\hom_{\mathcal{A}}(X,Y)$ with the differential
$\mu_1$ and view it as a filtered chain complex. We denote by
$\hom_{\mathcal{A}}^s(X,Y)$, $s \in \mathbb{R}$, the level-$s$
filtration subcomplex of $\hom_{\mathcal{A}}(X,Y)$.  We make three
further assumptions on our filtered $A_{\infty}$-categories:
\begin{itemize}
\item[(i)] $\mathcal{A}$ is strictly unital with the units lying in persistence level $0$. 
\item[(ii)] for every two objects $X, Y \in \Ob(\mathcal{A})$, the space
$\hom_{\mathcal{A}}(X,Y)$ is finite dimensional over $\k$.
\item[(iii)] $\mathcal{A}$ is complete with respect to
persistence shifts in the sense that we have a shift ``functor''
 which consists of a family of $A_{\infty}$-functors
$\Sigma = \{\Sigma^r: \mathcal{A} \longrightarrow \mathcal{A}, r \in
\mathbb{R}\}$ whose members satisfy the following conditions:
\begin{enumerate}
\item $\Sigma^r$ is strictly unital and the higher components
  $(\Sigma^r)_d$, $d \geq 2$, of $\Sigma^r$ all vanish.
\item \zjnote{$\Sigma^0 = \mathds{1}$}, $\Sigma^s \circ \Sigma^t = \Sigma^{s+t}$.
\item We are given prescribed identifications
  $\hom_{\mathcal{A}}^s(\Sigma^r X,Y) \cong
  \hom_{\mathcal{A}}^{s+r}(X,Y)$ that are compatible with the
  inclusions
  $\hom_{\mathcal{A}}^{\alpha}(X,Y) \subset
  \hom_{\mathcal{A}}^{\beta}(X,Y)$ for $\alpha \leq \beta$. These
  identifications are considered as part of the structure of the shift
  functor $\Sigma$.
\end{enumerate}
\end{itemize}
Given a filtered $A_{\infty}$-category $\mathcal{A}$ let
$Tw \mathcal{A}$ be the category of (filtered) twisted complexes over
$\mathcal{A}$.  This is itself a filtered $A_{\infty}$-category.  The
category $\mathcal{A}$ embeds into $Tw \mathcal{A}$ in an obvious way,
the embedding being a filtered $A_{\infty}$-functor which is full and
faithful.  Further, $Tw \mathcal{A}$ is pre-triangulated in the
filtered sense (which in particular means that it is closed under
formation of filtered mapping cones). As a result the homological
category
$$\mathcal{C}'\mathcal{A}:=H_0(Tw \mathcal{A})$$
is a TPC that contains the homological persistence category
$H_0(\mathcal{A})$ of $\mathcal{A}$.

\

Another TPC associated to the $A_{\infty}$-category $\mathcal{A}$ is
provided by the category $Fmod(\mathcal{A})$ of filtered
$A_{\infty}$-modules over $\mathcal{A}$.  We will only consider
strictly unital modules here. There is an obvious shift functor on
this category $\Sigma: (\R,+) \to \mathrm{End}(Fmod(\mathcal{A}))$ and
$Fmod(\mathcal{A})$ is in fact a filtered dg-category in the sense of
\zjnote{\cite[Chapter 2]{BCZ23}} and it is pre-triangulated. Thus
\zjnote{$H_{0}(Fmod(\mathcal{A}))$} is a TPC.  A \pbred{certain}
subcategory of $Fmod(\mathcal{A})$ \pbred{will be} of particular
interest \pbred{for us}.  Because our assumption of strict unitality
the Yoneda embedding is filtered:
$$\mathcal{Y}:\mathcal{A}\to Fmod(\mathcal{A})$$
and, thus, there exists a filtered quasi-isomorphism
$$\mathcal{M}(X)\to \hom_{Fmod}(\mathcal{Y}(X),\mathcal{M})$$ for all objects
$X$ of $\mathcal{A}$ and
$\mathcal{M}\in {\rm Obj}(Fmod(\mathcal{A}))$.  Let $\mathcal{A}^{\#}$
be the triangulated closure of the Yoneda modules and their
shifts. This is a \zjnote{full subcategory} of $Fmod(\mathcal{A})$
that has as objects all the iterated cones, over filtration preserving
morphisms, of shifts of Yoneda modules (thus of modules of the form
$\Sigma^{r}\mathcal{Y}(X)$).

Finally, let $\mathcal{A}^{\nabla}$ the smallest \zjnote{full sub
  dg-category} of $Fmod(\mathcal{A})$ that contains $\mathcal{A}^{\#}$
and all the modules (and all their shifts and translates) that are
$r$-quasi-isomorphic to objects in $\mathcal{A}^{\#}$,
$r\in [0,\infty)$.  A module $\mathcal{M}$ is $r$-quasi-isomorphic to
$\mathcal{N}$ if, in $H_{0}Fmod(\A)$, there is an $r$-isomorphism
$\mathcal{M}\to \mathcal{N}$.
 
It is easy to see that $\mathcal{A}^{\nabla}$ remains
pre-triangulated, carries the shift functor induced from
$Fmod(\mathcal{A})$ and thus the homological category
\pbred{$\mathcal{C}\mathcal{A}:=H_{0}(\A^{\nabla})$} is a TPC.
 
There are filtered functors:
$$\Theta: Tw(\A)\longrightarrow Fmod[Tw(\A)] \longrightarrow Fmod(\A) $$
where the first arrow is the Yoneda embedding and the second is
pullback over the natural inclusion $\A\to Tw(\A)$. The composition
$\Theta$ is a homologically full and faithful embedding, and it
induces a full and faithful embedding of TPCs. The image of $\Theta$
lands inside $\A^{\nabla}$ (actually inside the category denoted by
$A^{\#}$) and thus we have an inclusion of TPCs:
$$\bar{\Theta}: \mathcal{C}'\A\hookrightarrow \mathcal{C}\A$$
which, in general, is not a TPC- equivalence (by contrast to the
non-filtered case).  In particular, the category $\C\A$ contains many
more objects than $\C'\A$. Each object in $\C\A$ is $r$-isomorphic to
an object in $\C'\A$ but possibly only for some $r>0$ (to have a TPC
equivalence we would need to be able to take $r=0$).

\subsubsection{Filtered Fukaya categories}\label{subsubsec:filtr-fuk}
We work in the following setting: $(X, \omega = d \lambda)$ is a
Liouville manifold, \zjnote{i.e.~an exact symplectic manifold with a
  prescribed primitive $\lambda$ of the symplectic structure $\omega$,
  and such that $X$ is symplectically convex at infinity with respect
  to these structures}. The Lagrangians used here are triples
$L = (\bar{L}, h_{L}, \theta_{L})$ consisting of a closed oriented
exact Lagrangian submanifolds $\bar{L} \subset X$ equipped with a
primitive $h_L: \bar{L} \longrightarrow \mathbb{R}$ of
$\lambda|_{\bar{L}}$ and a grading $\theta_{L}$.  We refer to such a
pair $L$ as a marked Lagrangian submanifold and to $\bar{L}$ as its
underlying Lagrangian.

\pbred{Before we go on, let us describe our algebraic conventions
  regarding Floer theory and Fukaya categories.} \pbred{First of all,
  unless otherwise stated, all Fukaya categories will be assumed from
  now on to be over the base field $\k = \Z_{2}$.}

\pbred{The standard cohomology grading convention in Floer theory is
  that the $\mu_d$-operation has degree $2-d$, which is compatible
  with standard cohomology conventions of $A_{\infty}$-categories. The
  standard homology grading convention is that the $\mu_d$-operation
  has degree $d-2+n-dn$, where $n = \tfrac{1}{2}\dim X$ is the complex
  dimension of the ambient symplectic manifold $X$. In particular this
  means that the $\mu_2$-products have degree $-n$, namely
  $\mu_2: CF_i(L_0,L_1) \otimes CF_j(L_1,L_2) \longrightarrow
  CF_{i+j-n}(L_0,L_2)$ and the (chain level) units are in degree $n$,
  $e_L \in CF_n(L,L)$. This is also the case in Morse homology theory
  (where the units are maxima of \zjred{Morse functions}, hence of degree $n$,
  and represent the fundamental class of the manifold $L$). Obviously
  this is incompatible with the conventions described
  in~\S\ref{subsubsec:ainfty-alg}. In order to adjust the Floer
  homology setting to the one from~\S\ref{subsubsec:ainfty-alg} we
  perform a simple shift in grading as follows. We define:}
\pbred{$$\CCF_j(L_0,L_1) := CF_{n+j}(L_0,L_1), \;\; \HHF_j(L_0,L_1) :=
  HF_{n+j}(L_0,L_1), \; \forall j \in \mathbb{Z}.$$} \pbred{We now
  redefine the hom's in the Fukaya category $\fuk(X)$ such that
  $$\hom_{\fuk(X)}(L_0,L_1)_j = \CCF_j(L_0,L_1),$$} \pbred{and similarly for
  Floer the homologies $\HHF$.} \pbred{With these conventions the
  Fukaya category fits the conventions described
  in~\S\ref{subsubsec:ainfty-alg}, namely the $\mu_d$ operation has
  degree $d-2$ and the (homological) units are in degree $0$. The only
  apparent downside of this grading adjustment (which is purely
  aesthetic) is that with the standard grading, $CF_*(L,L)$ (when
  endowed with suitable auxiliary structures) is the Morse complex of
  $L$ and $HF_*(L,L) \cong H_*(L; \mathbb{Z}_2)$ (recall that we work
  with exact Lagrangians). However, with the new grading $\CCF_j(L,L)$
  and $\HHF_j(L,L)$ are concentrated in the negative degrees
  $-n \leq j \leq 0$, and we have
  $\HHF_j(L,L) \cong H_{n+j}(L;\mathbb{Z}_2)$. Another thing to note
  is that with the new grading conventions Poincar\'{e} duality takes
  the following form: $\HHF_j(L,L') \cong \HHF_{-n-j}(L',L)^*$, where
  the asterisk in the superscript stands for the dual.}


\pbred{We now proceed to the case of filtered Fukaya categories.} Fix
a collection of marked Lagrangians $\mathcal{X}$ in $X$. We assume
that $\mathcal{X}$ is closed under grading translations and shifts of
the primitives in the sense that if $L = (\bar{L}, h_L)$ is in
$\mathcal{X}$, then for every $r \in \mathbb{R}$, and
$k \in \mathbb{Z}$, the marked Lagrangian
\begin{equation}\label{eq:shift-tr}
\Sigma^r L[k] := (\bar{L}, h_L+r,\theta_{L}-k)
\end{equation}
is also in $\mathcal{X}$.  We denote by $ \bar{\mathcal{X}}$ the
family of underlying Lagrangians:
$$ \bar{\mathcal{X}}=\{\bar{L}\ |\ L\in\mathcal{X}\}~.~$$
We will assume that the family $\bar{\mathcal{X}}$ is finite and that
its elements are in generic position in the sense that any two
distinct Lagrangians $L', L'' \in \bar{\mathcal{X}}$ intersect
transversely and for every three distinct Lagrangians
$\bar{L}_0, \bar{L}_1, \bar{L}_2 \in \bar{\mathcal{X}}$ we have
\zjnote{$\bar{L}_0 \cap \bar{L}_1 \cap \bar{L}_2 = \emptyset$}.

Notice that we have \zjnote{implicitly} introduced here the shift
functor $\Sigma^{r}$, and the translation functors $T^{k}L:= L[-k]$
that act on the objects of $\mathcal{X}$ as in (\ref{eq:shift-tr}).

\

In Chapter 3 of \cite{BCZ23} are constructed, in the setting above,
filtered $A_{\infty}$ categories, $\fuk(\mathcal{X};\mathscr{P})$, of
Fukaya type, with objects the elements of $\mathcal{X}$. Here
$\mathscr{P}$ is a choice of perturbation data, picked roughly by the
scheme in \cite{Se:book-fukaya-categ}, but with special care, such
that the output is indeed filtered and not only weakly filtered. This
\zjnote{construction} is delicate and the constraint that
$\bar{\mathcal{X}}$ be finite and in generic position is used
repeatedly there. The resulting category
$\fuk(\mathcal{X};\mathscr{P})$ is strictly unital, with a unit in
filtration $0$, and has also the property that
$\hom_{\fuk(\mathcal{X};\mathscr{P})}(L,L)$ is concentrated in
filtration $0$ for all $L\in \mathcal{X}$. This happens because the
moduli spaces used to defined the $A_{\infty}$ operations appeal to
``cluster'' type configurations mixing Morse flow lines and Floer type
polygons such that, as a chain complex,
$(\hom_{\fuk(\mathcal{X};\mathscr{P})}(L,L), \mu_{1})$ coincides with
the Morse complex of a Morse function $f_{L}:L\to \R$ whose choice is
part of the data $\mathscr{P}$. Moreover, for two transverse
Lagrangians $L,L'$ the Floer data has $0$-Hamiltonian term (in other
works, it satisfies the homogeneous Floer equation).

It is also shown that any two such categories, defined for two
different allowable choices of perturbation data $\mathscr{P}$ and
$\mathscr{P}'$ are filtered quasi-equivalent in the sense that there
are filtered $A_{\infty}$-functors
$\fuk(\mathcal{X};\mathscr{P})\to \fuk(\mathcal{X};\mathscr{P}')$ that
are the identity on objects and induce a (filtered) equivalence of the
homological persistence categories.

One can then apply the discussion in \S\ref{subsubsec:ainfty-alg}
resulting in two associated TPC's.

The first is $\mathcal{C}\fuk(\mathcal{X})$, the \pbred{homological}
category
$$\mathcal{C}\fuk(\mathcal{X})=
H_{0}[\fuk(\mathcal{X};\mathscr{P})^{\nabla}]$$ where
$\fuk(\mathcal{X};\mathscr{P})^{\nabla}$ is the smallest triangulated
- with respect to weight-$0$ triangles - \zjnote{full subcategory} of
$Fmod(\fuk(\mathcal{X};\mathscr{P}))$ that contains the Yoneda modules
$\mathcal{Y}(L)$ with $L\in \mathcal{X}$ and is closed to
$r$-isomorphism for all $r$ in the sense that if $j:M\to M'$ is an
$r$-isomorphism of modules, and
$M\in \fuk(\mathcal{X};\mathscr{P})^{\nabla}$, then
$M'\in \fuk(\mathcal{X};\mathscr{P})^{\nabla}$. Possibly more
concretely, each object in this category is $r$-isomorphic, for some
$r$, to a \zjred{weight-$0$} iterated cone of Yoneda modules.  Two such
triangulated persistence categories, defined for different choices of
perturbation data are equivalent (as TPCs) and thus we drop the
reference to the perturbation data from the notation.

The second type of TPC, $\mathcal{C}'\fuk(\mathcal{X})$,  is defined by
$$\mathcal{C}'\fuk(\mathcal{X})=H_{0}[Tw(\fuk(\mathcal{X};\mathscr{P}))]$$
where $Tw(\fuk(\mathcal{X};\mathscr{P}))$ is the category of filtered
twisted complexes constructed from
$\fuk(\mathcal{X};\mathscr{P})$. Any two choices of perturbation data
produce equivalent TPCs, and we again drop $\mathscr{P}$ from the
notation.  There are filtered functors:
$$\Theta: Tw(\fuk(\mathcal{X};\mathscr{P}))
\longrightarrow Fmod[Tw(\fuk(\mathcal{X};\mathscr{P}))]
\longrightarrow Fmod(\fuk(\mathcal{X};\mathscr{P})) $$ and a
comparison functor:
$$\bar{\Theta}: \mathcal{C}'\fuk(\mathcal{X})\hookrightarrow
\mathcal{C}\fuk(\mathcal{X})~.~$$

Of course, there are also unfiltered categories here as the marked
Lagrangians in $X$ are the objects of an $A_{\infty}$-category, the
Fukaya category $\fuk(X)$ of $X$, constructed as in Seidel's book
\cite{Se:book-fukaya-categ}. In case of risk of confusion we denote
this category by $\fuk_{un}(X)$ to indicate that filtrations are
neglected in this case, and we use this \zjred{convention} for the other
structures in use, such as $A_{\infty}$ modules and so
forth. Similarly, we have the category $\fuk_{un}(\mathcal{X})$ that
has as objects only the elements of $\mathcal{X}$, as well as the
respective derived versions $D\fuk(X)$, and $D\fuk_{un}(\mathcal{X})$.

The basic relation between the filtered and unfiltered versions is
that there are equivalences of triangulated categories
$$[\C'\fuk(\mathcal{X}))]_{\infty}
\cong [\C\fuk(\mathcal{X})]_{\infty}\cong D\fuk_{un}(\mathcal{X})~.~$$

\medskip

Given that $\C'\fuk(\mathcal{X})$, and $\C\fuk(\mathcal{X})$ are
TPC's, their sets of objects can be endowed with fragmentation
pseudo-metrics of the type $d^{\mathcal{F}}$, as
in~\S\ref{subsec:tr-weight}, for any choice of objects $\mathcal{F}$.

\section{General properties of the $K$-group of a TPC}
\label{sec-k-TPC}

The purpose of this section is to describe several basic properties of
the $K$-group of a TPC.  In the whole section we fix a TPC
\zjnote{denoted by $\mathcal{C}$. Its $0$-level category
  $\mathcal{C}_{0}$ has the same objects as $\mathcal{C}$ but has as
  morphisms only the $0$-level morphisms in $\mathcal{C}$,
  $\ho_{\mathcal{C}_{0}}(A,B)=\ho^{0}_{\mathcal{C}}(A,B)$}. The group
we are most interested in is $K(\mathcal{C}_{0})$, the Grothendieck
group of the category $\mathcal{C}_{0}$ which, by the definition of a
TPC, is triangulated. We refer to this group as the $K$-group of
$\mathcal{C}$ and denote it by $K(\C)$. Note that $K(\C)$ has a
natural structure of a $\Lambda_P$-module by defining
$t^r \cdot [A] = [\Sigma^r A]$ for every $r \in \mathbb{R}$.

We start in \S\ref{ssec-acyc-cat} with some remarks on
acyclics. Subsection \ref{subsubsec:TPC-funct-K} shows how TPC
functors - properly defined - induce morphisms at the $K$-group
level. In \S\ref{subsubsec:ring-K} we discuss a ring structure on the
$K$-group which is defined in case the TPC admits a tensor structure -
we call such a TPC monoidal.  In \S\ref{sec-op-tpc} we discuss the
opposite category of a TPC and its K-group. In \S
\ref{subsec:infty-lev} we consider the $\infty$-level of a TPC and the
associated $K$-group.  We package most of these properties in Theorem
\ref{thm:alg-pack1} in \S\ref{subsec:conclusion-alg}.

\subsection{The acyclic category of a TPC} \label{ssec-acyc-cat} Let
$\mathcal{C}$ be a TPC.  An important feature of a TPC is that there
is a class of objects that are ``close to $0$'' in a quantifiable
way. Explicitly, recall that an object $X \in {\rm Obj}(\C)$ is
\emph{$r$-acyclic} for $r \geq 0$, denoted by $X \simeq_r 0$, if
\zjnote{$\mathds{1}_{X}\simeq_{r} 0$} (see~\S\ref{subsec:tpc}, and for more
details~\cite{BCZ23}).

As mentioned in \S\ref{subsec:tpc}, there are 
 two (sub)categories.
\begin{itemize}
\item[(i)] The full subcategory of $\C$, denoted by $\mathcal{AC}$,
  whose objects are $r$-acyclic for any $r \geq 0$.
\item[(ii)] The full subcategory of $\C_0$, denoted by
  $\mathcal{AC}_0$, whose objects are $r$-acyclic for any $r \geq 0$.
\end{itemize}
The subcategory $\mathcal{AC}$ is obviously a persistence
category, since $\ho_{\mathcal{AC}} = \ho_{\C}$ which
admits a persistence module structure. Moreover, for
$X, Y \in {\rm Obj}((\mathcal{AC})_0) = {\rm Obj}(\mathcal{AC}) = {\rm
  Obj}(\mathcal{AC}_0)$, we have
\[ \ho_{(\mathcal{AC})_0}(X,Y) = \ho^0_{\mathcal{AC}}(X,Y)
  = \ho_{\mathcal{AC}_0}(X,Y).\] This yields
$(\mathcal{AC})_0 = \mathcal{AC}_0$. Moreover, the subcategory $\mathcal{AC}$ is a TPC and
the category $\C_{\infty}$ is the Verdier localization of $\C_{0}$
 \zjred{with} respect to $\mathcal{AC}_{0}$ \zjnote{(for more details, see \cite{BCZ23})}. Thus, the $K$- group $K(\A\C)$ of $\A\C$ is
well-defined. Note that similarly to $K(\C)$ the group $K(\A\C)$ also
has the structure of a $\Lambda_P$-module by defining
$t^r \cdot [A] = [\Sigma^r A]$.

\medskip

The two $K$-groups  $K(\C)$ and $K(\mathcal{AC})$
are closely related.  First, there is an obvious homomorphism,
\begin{equation} \label{j-map} j: K(\mathcal{AC}) \to K(\C)
  \,\,\,\,\mbox{defined by} \,\,\,\, X \mapsto X.
\end{equation}
This is well-defined since an exact triangle in $\mathcal{AC}_0$ is
also an exact triangle in $\C_0$. 

Additionally (and non-trivially), we also  have morphisms  in the opposite direction\zjnote{, called {\em acyclic truncations},} defined as follows. For any fixed $r \geq 0$ and for any
$X \in {\rm Obj}(\C)$, by \zjnote{(iv) in \S \ref{subsec:tpc}},
we have an exact triangle
\begin{equation} \label{iii-TPC} \Sigma^r X \xrightarrow{\eta_r^X} X
  \to Q_{r,X} \to T\Sigma^r X
\end{equation}
in $\C_0$ for some $r$-acyclic object
$Q_{r,X} \in {\rm Obj}(\mathcal{AC}_0)$. Consider the following map
\begin{equation} \label{Q-map} Q_r: K(\C) \to
  K(\mathcal{AC})\,\,\,\,\mbox{defined by} \,\,\,\, X \mapsto
  Q_r(X): = Q_{r,X}.
\end{equation}
Note that since $Q_{r,X}$ is $r$-acyclic and defined up to isomorphism
in $\C_0$, we know that $Q_{r,X}$ (or more precisely the class
$[Q_{r,X}]$) is a well-defined element in $K(\mathcal{AC})$.

\begin{lemma} \label{lemma-Q} The map $Q_r$ defined in (\ref{Q-map})
  is well-defined, and it is a homomorphism. \end{lemma}

\begin{proof} For the first conclusion, we need to prove that if $X$
  is replaced by $Y-Z$ in $K(\C)$, that is, $X = Y-Z$, then
  $Q_r(X) = Q_r(Y) - Q_r(Z)$ in $K(\mathcal{AC})$. By the defining
  relations in the corresponding $K$-groups, this means that the
  existence of the exact triangle $X \xrightarrow{f} Y \to Z \to TX$
  in $\C_0$ implies that \zjred{the triangle defined below}
  \begin{equation} \label{exact-tri-Q} Q_{r,X} \to Q_{r,Y} \to Q_{r,Z}
    \to TQ_{r,X}
  \end{equation}
  is also exact in $\mathcal{AC}_0$, where $Q_{r,X}, Q_{r,Y}$, and
  $Q_{r,Z}$ are the $r$-acyclic objects from the exact triangles in
  the form of (\ref{iii-TPC}). In fact, by \zjnote{(iii) in \S \ref{subsec:tpc}}, we have an exact triangle
  $\Sigma^rX \xrightarrow{\Sigma^r f} \Sigma^r Y \to \Sigma^r Z \to T
  \Sigma^r X$. Then consider the following diagram,
  \[ \xymatrix{ \Sigma^r X \ar[r]^-{\Sigma^r f} \ar[d]_-{\eta^X_r} &
      \Sigma^r Y \ar[r] \ar[d]_-{\eta^Y_r}
      & \Sigma^r Z \ar@{-->}[d] \ar[r] & T\Sigma^r X \ar[d] \\
      X \ar[r]^{f} \ar[d] & Y \ar[r] \ar[d] & Z \ar@{-->}[d] \ar[r]
      & TX \ar[d] \\
      Q_{r,X} \ar@{-->}[r] \ar[d] & Q_{r,Y} \ar@{-->}[r] \ar[d]
      & A \ar@{-->}[d] \ar@{-->}[r] & TQ_{r,X} \ar[d]\\
      T\Sigma^r X \ar[r] & T\Sigma^r Y \ar[r] & T\Sigma^r Z \ar[r] &
      T^2\Sigma^r X}
  \]
  We know that $\eta_r^Y \circ \Sigma^rf = f \circ \eta_r^X$ and the
  first two rows and first two columns are all exact triangles. Then
  by $3\times3$-lemma (see \zjred{Proposition 1.1.11 in \cite{BBD82}}), there exists
  an object $A \in {\rm Obj}(\C)$ (in fact,
  $A \in {\rm Obj}(\mathcal{AC})$) and dashed arrows such that all the
  rows and columns are exact, as well as the diagram commutes (except
  the most right-bottom square is anti-commutative). In particular, we
  have an exact triangle
  \[ Q_{r,X} \to Q_{r,Y} \to A \to T Q_{r,X}, \] as the third row in
  the diagram. Moreover, since the third column
  $\Sigma^r Z \to Z \to A \to T\Sigma^r Z$ is also exact, we have
  $A \cong Q_{r,Z}$. Hence, we obtain the exact triangle
  (\ref{exact-tri-Q}).

  To see that $Q_{r}$ is a homomorphism, we need to show that
  $Q_r(X+Y) = Q_r(X) + Q_r(Y)$.  In the $K$-group, the sum $X+Y$ is
  represented by an object $Z \in {\rm Obj}(\C)$ that fits into an
  exact triangle $X \to Z \to Y \to TX$. Then the same argument as in
  the first part shows that there is an exact triangle
  $Q_{r,X} \to Q_{r,Z} \to Q_{r,Y} \to TQ_{r,X}$. Therefore,
  \begin{align*}
    Q_r(X+Y) = Q_r(Z) & = Q_{r,Z} \\
                      & = Q_{r,X} + Q_{r,Y} = Q_r(X) + Q_r(Y) 
  \end{align*}
  in $K(\mathcal{AC}_0)$. Similarly, we have
  $Q_{r}(-X) = Q_{r, TX} = TQ_{r,X} = - Q_{r,X}$ in
  $K(\mathcal{AC}_0)$.
\end{proof}

\subsection{TPC functors and induced
  $K$-morphisms}\label{subsubsec:TPC-funct-K}
We formalize here some of the constructions before in a more abstract
fashion. We denote a TPC by a triple $(\C, T_{\C}, \Sigma_{\C})$ (or
$(\C, T, \Sigma)$) to emphasize the role of the shift functor $\Sigma$
and that of the translation functor $T$ in $\mathcal{C}$. \zjnote{Let us recall the following definition (see Definition 2.25 in \cite{BCZ23}).} 

\begin{dfn} \label{dfn-tpc-functor} Let $(\C, T_{\C}, \Sigma_{\C})$
  and $(\D, T_{\D}, \Sigma_{\D})$ be TPCs. A {\it TPC functor}
  $\F: (\C, T_{\C}, \Sigma_{\C}) \to (\D, T_{\D}, \Sigma_{\D})$ is a
  persistence functor satisfying the following conditions.
  \begin{itemize}
  \item[(i)] For any $r \in \R$, we have
    $\Sigma_{\D}^r \circ \F = \F \circ \Sigma_{\C}^r$. Moreover, for
    any $r,s \in \R$, we have the following commutative diagram,
    \[ \xymatrixcolsep{4pc} \xymatrix{ \F(\Sigma_{\C}^r X)
        \ar[r]^-{\F((\eta_{r,s})_X)} \ar[d]_-{=}
        & \F(\Sigma_{\C}^s X) \ar[d]^-{=}\\
        \Sigma_{\D}^r \F(X) \ar[r]^-{(\eta_{r,s})_{\F(X)}} &
        \Sigma_{\D}^s \F(X)} \] for any $X \in {\rm Obj}(\C)$.
  \item[(ii)] The restriction $\F|_{\C_0}: \C_0 \to \D_0$ is a
    triangulated functor (recall that by definition of a TPC, the
    $0$-level categories $\C_0$ and $\D_0$ are genuine triangulated
    categories).
  \end{itemize}
\end{dfn}

Here are some useful properties of TPC functors.

\begin{lem} \label{prop-tpc} Let
  $\F: (\C, T_{\C}, \Sigma_{\C}) \to (\D, T_{\D}, \Sigma_{\D})$ be a
  TPC functor. Then we have the following properties.
  \begin{itemize}
  \item[(i)] For any $X \in {\rm Obj}(\C)$ and $r \geq 0$, we have
    $\F(\eta_r^X) = \eta_r^{\F(X)}$.
  \item[(ii)] If $X \in {\rm Obj}(\C)$ is $r$-acyclic, then $\F(X)$ is
    also $r$-acyclic.
  \item[(iii)] If $\phi: X \to Y$ is an $r$-isomorphism, then
    $\F(\phi): \F(X) \to \F(Y)$ is an $r$-isomorphism too.
  \item[(iv)] If
    $\Delta: A \xrightarrow{\bar{u}} B \xrightarrow{\bar{v}} C
    \xrightarrow{\bar{w}} \Sigma^{-r} T A$ is a strict exact triangle
    of weight $r\geq 0$, then
    \[ \F(\Delta): \F(A) \xrightarrow{\F(\bar{u})} \F(B)
      \xrightarrow{\F(\bar{v})} \F(C) \xrightarrow{\F(\bar{w})}
      \Sigma^{-r} T \F(A)\] is also a strict exact triangle of weight
    $r$.
  \item[(v)] $\F$ descends to a group homomorphism
    $\F_{K}: K(\C) \to K(\D)$. The same conclusion holds for $K(\A\C)$
    and $K(\A\D)$ as well as for $K_{\infty}(\mathcal{C})$ and
    $K_{\infty}(\mathcal{D})$.
  \end{itemize}
\end{lem}

\begin{proof}
  The proof is rather straightforward so we only outline the main
  points. Indeed, the fact that the restriction of $\mathcal{F}$ to
  $\C_{0}$ is triangulated and that $\mathcal{F}$ is a persistence
  functor that commutes with the shift functors and the natural
  transformations $\eta$ implies easily that $\mathcal{F}$ sends
  $r$-acyclics to $r$-acyclics.  One then deduces that $\mathcal{F}$
  transforms $r$-isomorphisms into $r$-isomorphisms and the rest of
  the properties stated are easy to check.
\end{proof}

Given any TPC $(\C, T_{\C}, \Sigma_{\C})$, recall that $K(\C)$ and
$K(\A\C)$ are related by the group homomorphisms $Q_r$ defined in
(\ref{Q-map}) and $j$ defined in (\ref{j-map}). It is easy to see that
we also have:

\begin{lemma} \label{lem-q-j-fx} For any TPC
  $(\C, T_{\C}, \Sigma_{\C})$, we have the following relations,
  \[ j \circ \F_K = \F_K \circ j \,\,\,\,\mbox{and} \,\,\,\, Q_r \circ
    \F_K = \F_K \circ Q_r \] for any $r \geq 0$.
\end{lemma}

To summarize what we have shown till now, triangulated persistence
categories together with TPC functors form a category \zjnote{denoted by $\mathcal{T}PC$}.  There is an {\em acyclic} functor
$$\ac : \mathcal{T}PC \longrightarrow \mathcal{T}PC
\ \ , \ \ \ \ \C \longmapsto \ac\C ~.~$$ The morphisms $j$ and $Q_{r}$
can be viewed as appropriate natural transformations relating the
functor $\ac$ to the identity.

Recall from the beginning of~\S\ref{sec-k-TPC} that for a TPC $\C$ the
$K$-group $K(\C)$ has the structure of $\La_{P}$-module defined by
$t^{r}\cdot [A] = [\Sigma^{r} A]$.  Consider the category of modules
$\md_{\La_{P}}$ over $\La_{P}$. With these conventions we have
constructed two Grothendieck group functors:
\begin{equation}\label{eq:K-funct-1}
  K: \mathcal{T}PC \longrightarrow \md_{\La_{P}}
\end{equation}
that associates to a TPC $\C$ its $K$-group $K(\C)$ and to a TPC
functor $\mathcal{F}$ the induced morphism $\mathcal{F}_{K}$ and
\begin{equation}
  K\ac :\mathcal{T}PC \longrightarrow \md_{\La_{P}}
\end{equation}
that is defined by $K\ac =K\circ \ac$ and associates to $\C$ the
module $K(\ac \C)$.  We will often denote this module by $K\ac (\C)$
from now on.

\subsection{Monoidal structures on $\C$ and ring structures on
  $K(\mathcal{C})$}
\label{subsubsec:ring-K}
We start here with a natural generalization of the notion of tensor
category (or monoidal category, see Section 1 in Chapter VII of
\cite{Mac-cat-book}) to the TPC setting.

\begin{dfn} \label{dfn-tensor} Let $\C$ be a TPC. A {\it tensor
    structure} on $\C$ is an associative symmetric bi-operator
  $\otimes: {\rm Obj}(\C) \times {\rm Obj}(\C) \to {\rm Obj}(\C)$
  together with a unit element $e_{\C} \in {\rm Obj}(\C)$ (with
  respect to $\otimes$) such that the following axioms are satisfied.
  \begin{itemize}
  \item[(i)] For any $A \in {\rm Obj}(\C)$, we have $(-) \otimes A$ is
    a persistence functor on $\C$. Explicitly, for any
    $X,Y \in {\rm Obj}(\C)$ and $f \in \ho_{\C}^{\alpha}(X,Y)$ for
    $\alpha \in \R$,
    \[ ((-) \otimes A)(X) = X \otimes A\,\,\,\,\mbox{and}\,\,\,\,
      ((-) \otimes A)(f) = f \otimes \mathds{1}_A \in
      \ho_{\C}^{\alpha}(X \otimes A, Y \otimes A). \]
  \item[(ii)] For any $A \in {\rm Obj}(\C)$, the restriction of the
    functor $(-) \otimes A$ to $\C_0$ is triangulated, i.e., it maps
    exact triangles in $\C_0$ to exact triangles in $\C_0$. Moreover,
    $(-)\otimes A$ is additive on $\C_0$.
  \item[(iii)] For any $A \in {\rm Obj}(\C)$, the functor
    $(-) \otimes A$ is compatible with the translation
    auto-equivalence $T$ and the shift functor $\Sigma^r$ for any
    $r \in \R$ in the following sense,
    \[ TX \otimes A \cong T(X \otimes A) (\cong X \otimes TA)
      \,\,\,\,\mbox{in $\C_0$}, \] and there exists a family of {\it
      isomorphisms}
    $\{\theta_r: \Sigma^r X \otimes A \to \Sigma^r(X \otimes A)\}_{r
      \in \R}$ in $\C_0$ with $\theta_0 = \mathds{1}_{X \otimes A}$
    such that the following diagram commutes,
    \[ \xymatrixcolsep{5pc} \xymatrix{ \Sigma^r X \otimes A
        \ar[r]^-{(\eta_{r,s})_X \otimes \mathds{1}_A}
        \ar[d]_-{\theta_r} & \Sigma^s X \otimes A \ar[d]^-{\theta{s}} \\
        \Sigma^r (X \otimes A) \ar[r]^-{(\eta_{r,s})_{X \otimes A}}&
        \Sigma^s (X \otimes A)}\] for any $r, s \in \R$. A similar
    conclusion holds for the functor $A \otimes (-)$ due to the
    symmetry of $\otimes$.
  \end{itemize}
  A monoidal TPC is a category $\C$ as before together with a choice
  of tensor structure.
\end{dfn}

\begin{remark} \label{rmk-axiom1} The axiom (i) in Definition
  \ref{dfn-tensor} translates to the following commutative diagram
  \[ \xymatrixcolsep{4pc} \xymatrix{ \ho_{\C}^{\alpha}(X, Y)
      \ar[r]^-{(-) \otimes \mathds{1}_A} \ar[d]_-{i_{\alpha,
          \alpha+r}} & \ho_{\C}^{\alpha}(X\otimes A, Y\otimes A)
      \ar[d]^-{i_{\alpha, \alpha+r}}\\
      \ho_{\C}^{\alpha+r}(X, Y) \ar[r]^-{(-) \otimes \mathds{1}_A} &
      \ho_{\C}^{\alpha+r}(X\otimes A, Y\otimes A)}\] for any
  $\alpha \in \R$ and $r \in \R_{\geq 0}$. In particular, we have for
  any $f \in \ho_{\C}^{\alpha}(X, Y)$,
  $i_{\alpha, \alpha+r}(f) \otimes \mathds{1}_A = i_{\alpha,
    \alpha+r}(f \otimes \mathds{1}_A)$.
\end{remark}

\begin{rem} \label{lem-1} Recall that $\mathcal{AC}$ is the
  subcategory of $\C$ whose objects consists of only acyclic
  objects. It is a simple exercise to check that any
  $A \in {\rm Obj}(\C)$, the functor $(-) \otimes A$ restrict to an
  endofunctor of $\mathcal{AC}$.
\end{rem}

\begin{ex} \label{ex:tensor-ch} A basic example of a monoidal TPC's is
  the persistence homotopy category $\C^{\fg}=H_{0}(\fchfg)$ of
  finitely generated filtered chain complexes. The tensor structure is
  induced by tensor products of (filtered) chain complexes and the
  unit element is given by the chain complex \zjnote{$E = \k \left<x\right>$} which is
  concentrated in degree $0$, has $\partial x =0$ and whose filtration
  function $\ell$ is defined by $\ell(x)=0$.

  In the same vein, the persistence homotopy category
  $\C^{\ta}=H_{0}(\fcht)$ of tame filtered chain complexes (see
  page~\pageref{p:types-ch}) is a monoidal TPC too.
\end{ex}

Suppose that $\C$ is monoidal. For $X, Y \in \Ob(\mathcal{C})$ with 
classes $[X], [Y]\in K(\C)$ define
\begin{equation} \label{dfn-ring}
  [X] \cdot [Y] = [X \otimes Y] \in K(\C)  ~.~
\end{equation} 

\begin{lemma} \label{lem-3} If $\C$ admits a tensor structure, then
  the operation in~\eqref{dfn-ring} is well-defined and induces a
  commutative unital ring structure on $K(\C)$ and a non-unital one on
  $K(\ac \mathcal{C})$. The same operation makes $K(\ac \mathcal{C})$
  into a $K(\C)$-module, and the map~\eqref{j-map} is $K(\C)$-linear.
\end{lemma}
The proof can be found in the expanded version of the
paper~\cite{BCZ-pkt-arxiv}.


\begin{rem}
  Recall that $K(\mathcal{C})$ is also a $\Lambda_P$-module. This
  module structure is compatible with the preceding ring structure
  (defined when $\mathcal{C}$ is monoidal) and together these two
  structures make $K(\mathcal{C})$ into an algebra over $\Lambda_P$. A
  similar statement holds for $K\ac (\C)$.
\end{rem}

  

\subsection{The opposite TPC and its $K$-group} \label{sec-op-tpc}

Recall that every triangulated category $\C$ has its opposite category
denoted by $\C^{\rm op}$ which is triangulated as well. The exact
triangles in $\C^{\rm op}$ are related to those of $\mathcal{C}$ as
follows (we follow the construction given by Remark 10.1.10 (ii) in
\cite{KS-cat-sheaves}).  Denote by
${\rm op}: \C \longrightarrow \C^{\rm op}$ the canonical
(contravariant) functor that (i) is identity on
${\rm Obj}(\C) (= {\rm Obj}(\C^{\rm op}))$ and (ii) reverses the
domain and target on the morphisms in $\C$, that is,
${\rm op}(f) \in {\rm Hom}_{\C^{\rm op}}(Y,X)$ iff
$f \in {\rm Hom}_{\C}(X,Y)$. Assume that $T$ is the translation
autoequivalence of $\C$, then define
\begin{equation} \label{t-op} T^{\rm op}: = {\rm op} \circ T^{-1}
  \circ {\rm op}^{-1}: \C^{\rm op} \to \C^{\rm op}
\end{equation}
to be the translation autoequivalence in $\C^{\rm op}$. An
exact triangle in $\C^{\rm op}$ is defined by
\begin{equation} \label{ex-tri-op} Z \xrightarrow{{\rm op}(g)} Y
  \xrightarrow{{\rm op}(f)} X \xrightarrow{{\rm op}(- T^{-1}(h))}
  T^{\rm op}Z
\end{equation}
whenever $X \xrightarrow{f} Y \xrightarrow{g} Z \xrightarrow{h} TX$ is
an exact triangle in $\C$. It is easy to see that $\C^{\rm op}$ is a
triangulated category.

\begin{dfn} \label{dfn-opp-tpc} Let $(\C, T, \Sigma)$ be a TPC. Define
  its {\it opposite} by $(\C^{\rm op}, T^{\rm op}, \Sigma^{\rm op})$
  where $T^{\rm op}$ is defined in (\ref{t-op}) and
  \begin{equation} \label{shift-op} (\Sigma^{\rm op})^r: = {\rm op}
    \circ \Sigma^{-r} \circ {\rm op}^{-1}: \C^{\rm op} \to \C^{\rm op}
  \end{equation}
  for every $r \in \R$. Moreover, the persistence module structure on
  the $\ho_{\C}$'s induces in an obvious way such a structure on the
  $\ho_{\C^{\rm op}}$'s. This makes $\mathcal{C}^{\text{op}}$ into a
  persistence category and the canonical functor ${\rm op}$ is a
  (contravariant) persistence functor. Moreover, as we will see in
  Lemma~\ref{lem-opp} and~\ref{lem-kappa-opp} below,
  $(\C^{\rm op}, T^{\rm op}, \Sigma^{\rm op})$ defined as above is a
  TPC and the functor $\text{op}$ is a TPC functor.
\end{dfn}

A few remarks are in order.

\begin{remark} \label{rmk-opp}
  \begin{itemize}
  \item[(i)] The collection
    $\Sigma^{\rm op} = \{(\Sigma^{\rm op})^r\}_{r \in \R}$ forms a
    shift functor on $\C^{\rm op}$, and the natural transformation
    $\eta^{\rm op}_{r,s}: (\Sigma^{\rm op})^r \to (\Sigma^{\rm op})^s$
    is defined by
    \begin{equation} \label{op-nt} \eta^{\rm op}_{r,s}: = {\rm
        op}(\eta_{-s,-r})
    \end{equation}
    for any $r,s \in \R$, where
    $\eta_{-s,-r}: \Sigma^{-s} \to \Sigma^{-r}$ is the natural
    transformation for $\Sigma$.
  \item[(ii)] Recall that for each object $X$ in $\C$, the morphism
    $\eta^X_r = \iota_{-r,0}(\eta_{r,0}(X)) \in \ho_{\C}^0(\Sigma^r X,
    X)$ plays an important role in the definition of a TPC. Similarly
    in $\C^{\rm op}$, the counterpart is
    \begin{align*}
      (\eta^{\rm op})_r^X = \iota_{-r,0}(\eta^{\rm op}_{r,0}(X))
      & =i_{-r,0}({\rm op}(\eta_{0, -r}(X))) \\
      & = {\rm op}(\eta^X_r) \in \ho_{\C^{\rm op}}^0((\Sigma^{\rm op})^rX, X)
    \end{align*}
    where in the third equality we use the condition that ${\rm op}$
    is a persistence functor, commuting with the persistence structure
    morphisms.
  \item[(iii)] The $0$-level category of $\mathcal{C}^{\text{op}}$
    satisfies $(\C^{\rm op})_0 = (\C_0)^{\rm op}$, since the canonical
    map ${\rm op}$ is a persistence functor.
  \end{itemize}
\end{remark}

\begin{lemma} \label{lem-opp} Let $(\C, T, \Sigma)$ be a TPC, then its
  opposite $(\C^{\rm op}, T^{\rm op}, \Sigma^{\rm op})$ is also a
  TPC.
\end{lemma}

We refer the reader to the expanded version of this
paper~\cite{BCZ-pkt-arxiv} for the proof.

For brevity, denote by $\C^{\rm op}_0$ the $0$-level category of
$\C^{\rm op}$ (there is no ambiguity of the notation due to (iii) in
Remark \ref{rmk-opp}).

\begin{lemma} \label{lem-kappa-opp} The canonical functor
  ${\rm op}: (\C, T, \Sigma) \to (\C^{\rm op}, T^{\rm op}, \Sigma^{\rm
    op})$ is a TPC functor. Moreover, it induces a group {\rm
    isomorphism} ${\rm op}_K: K(\C) \to K(\C^{\rm op})$.
\end{lemma}

The proof is straightforward and can be found in~\cite{BCZ-pkt-arxiv}.


\subsection{The $\infty$-level of a TPC and an exact
  sequence}\label{subsec:infty-lev}
Recall from~\S\ref{subsec:tpc}, page~\pageref{p:C-infty}, that for
every TPC $\C$ we can associate its $\infty$-level category
$\C_{\infty}$ \zjnote{with
  $\ho_{\C_{\infty}}(A,B)=\lim_{\longrightarrow}
  \ho_{\C}^{\alpha}(A,B).$} As mentioned in \S\ref{subsec:tpc}, this
category is identified with the Verdier localization of
$\C_{0}/\mathcal{AC}_{0}$. Thus $\C_{\infty}$ is triangulated (see
\jznote{\cite[Proposition 2.38]{BCZ23}}) and its $K$-group
$K(\C_{\infty})$ is \zjnote{well-defined}. The exact triangles of
$\C_{\infty}$ contain those of $\C_{0}$ and it is easy to see that
many of the properties of $K(\C)$ remain true for $K(\C_{\infty})$:
\zjnote{\begin{itemize}
\item[(i)] $K(\C_{\infty})$ is a \zjnote{$\La_{P}$-module}. If $\C$ is monoidal, then $K(\C_{\infty})$ is \zjred{a} $\La_{P}$-algebra.
\item[(ii)] \zjnote{There is a functor
    $K_\infty : \mathcal{T}PC \longrightarrow \md_{\La_{P}}$ defined
    by $\C \longmapsto K_{\infty}(\C) :=K(\C_{\infty})$.}
\end{itemize}} \zjnote{Here} (as in~\S\ref{subsubsec:TPC-funct-K})
$\md_{\La_{P}}$ \zjnote{denotes} the category of modules over
$\La_{P}$.

\medskip

We will now relate the groups $K(\C)$ and $K(\C_{\infty})$ via an
exact sequence. Of course, this relation is to be expected given that
$\C_{\infty}=\C_{0}/\mathcal{AC}_{0}$. \zjnote{Note that} there is a surjective homomorphism $\pi: K(\C) \longrightarrow K_{\infty}(\C)$
of $\Lambda_P$-modules which is induced by the identity on objects.
Further, the inclusion $\Ob(\mathcal{AC}) \subset \Ob(\mathcal{C})$ gives
rise to a well defined map
$j: K\ac(\C)\stackrel{j}{\longrightarrow} K(\C)$. Let
$\tor K(\C) := \ker(j)$ and denote by
$\psi: \tor K(\C) \longrightarrow K\ac(\C)$ the inclusion.
Thus we have the following sequence of $\La_{P}$-modules:
\begin{equation}\label{eq:exact-seq-k}
  0\longrightarrow \tor K(\C)\stackrel{\psi}{\longrightarrow}
  K\ac(\C)\stackrel{j}{\longrightarrow}
  K(\C)\stackrel{\pi}{\longrightarrow} K_{\infty}(\C)\longrightarrow 0.
\end{equation}  
From the general properties of Verdier localization it is easy to deduce \zjnote{the following property}.

\begin{prop} \label{prop:exact-seq} The
  sequence~\eqref{eq:exact-seq-k} is exact.
\end{prop} 

\begin{rem}\label{rem:tor}
\zjnote{(a)} As mentioned before, this follows from the general properties of Verdier
localization but it is useful to have a more explicit description of $\tor K(\C)$
together with the \zjnote{exactness} argument 
and we provide them here.
  Consider the \zjnote{following} diagram of $\La_{P}$-modules.
  \begin{equation}\label{eq:big-diag}\xymatrix{
      &    &  &          0 \ar[d] &  \\
      &  0\ar[d]  & 0\ar[d] &  \widehat{\mathcal{R}}_{\C} \ar[d]^{s}  & \\
      0 \ar[r]  & \mathcal{R}_{\ac \C_{0}} \ar[d]\ar[r]^{\bar{j}'}
      & \mathcal{R}_{\C_{0}}\ar[r]^{\overline{\pi}'}\ar[d]
      & \mathcal{R}_{\infty} \ar[r]\ar[d]_{u} & 0\\
      0\ar[r] & G \langle \ac\C \rangle \ar[d]\ar[r]^{\bar{j}}
      & G \langle \C \rangle \ar[d]\ar[r]^{\overline{\pi}}
      & G_{\infty} \ar[r]\ar[d]_{v} & 0\\
      & K\ac (\C) \ar[r]^{j}\ar[d] & K(\C) \ar[r]^{\pi} \ar[d]
      & K_{\infty} (\C)\ar[r] \ar[d] & 0\\
      &  0 & 0 & 0 &
    }
\end{equation}
The diagram is constructed as follows. We denote by
$G \langle - \rangle$ the free abelian group generated by the objects
of the respective categories and by $\mathcal{R}_{-}$ the subgroup of
triangular relations associated to the triangular structure of the
$0$-level of the respective TPC, as in the definition of the
$K$-group.  Thus, by the definition of the $K$ groups in question, the
first two columns are exact. The maps $\bar{j}$ and $\bar{j}'$ are
inclusions and $j$ is induced by them. The groups $G_{\infty}$ and
$\mathcal{R}_{\infty}$ are defined as the obvious quotients making the
top two rows exact and $\overline{\pi}$, $\overline{\pi}'$ are the
projections.  The maps $u$ is induced from the left, top square in the
diagram. The group $\widehat{\mathcal{R}}_{\C}$ is the kernel of $u$.
All the groups and morphisms discussed here are also $\La_{P}$-modules
so $\widehat{\mathcal{R}}_{\C}$ is also a $\La_{P}$-module.  It is
immediate to see that $\pi\circ j=0$ because an acyclic object is, by
definition, $r$-isomorphic to $0$ for some $r\in \R$, and thus it is
isomorphic to $0$ in $\C_{\infty}$. This means that the map $v$ is
well-defined and induced by the left bottom square in the diagram. 

Diagram chasing shows that the exactness in the proposition follows by proving  that $v$ is surjective and, moreover, this implies that the kernel of $j$ is  identified to $\widehat{\mathcal{R}}_{\C}$,
$$\tor K(\C)=\ker (j)= \widehat{\mathcal{R}}_{\C}~.~$$  Given that
$K_{\infty}(\C)=G \langle \C \rangle / \mathcal{R}_{\C_{\infty}}$, the surjectivity of $v$ comes down to the
following identity
\zjred{\begin{equation}\label{eq:relations1}
  \mathcal{R}_{\C_{\infty}}= G \langle \ac \C \rangle +
  \mathcal{R}_{\C_{0}}
\end{equation}}
which follows from the description of the exact triangles in $\C_{\infty}$.

\zjnote{(b)}  It is useful to express the group
  $\tor K (\C)=\widehat{\mathcal{R}}_{\C}$ in yet a different way.  Notice
  that we have the obvious inclusion:
  $$\mathcal{R}_{\ac \C_{0}}\hookrightarrow G \langle \ac\C \rangle
  \cap \ \mathcal{R}_{\C_{0}}~.~$$ More diagram chasing shows that:
  \begin{equation}\label{eq:tor-2}
    \tor K(\C)= [G \langle \ac\C \rangle \cap \
    \mathcal{R}_{\C_{0}}]/\mathcal{R}_{\ac \C_{0}}.
  \end{equation}
\end{rem}

\zjnote{(c)} It is easy to see that $\tor K (\C)$ can be viewed as a functor:
$$\tor K : \mathcal{T}PC \longrightarrow \md_{\La_{P}}$$
and that $\psi: \tor K(\C) \longrightarrow K\ac(\C)$
from~\eqref{eq:exact-seq-k} is part of a natural transformation
between $\tor K$ and $\ac K$.

\subsection{Wrap-up of general
  properties} \label{subsec:conclusion-alg} We assemble here the
properties discussed above in one statement, by expanding on the
discussion in \S\ref{subsubsec:TPC-funct-K}.  Recall from
\S\ref{subsubsec:TPC-funct-K} that $\md_{\La_{P}}$ is the category of
modules over $\La_{P}$. Denote by
${\rm op}:\md_{\La_{P}}\to \md_{\La_{P}}^{\rm op}$ the obvious
${\rm op}$ functor with the particular feature that
${\rm op}(M)= \overline{M}$ where $\overline{M}$ is the same abelian
group as $M$ but the $\La_{P}$-module structure is
$t^{r} \cdot \overline{a} = \overline{t^{-r} \cdot a}$ where $a\in M$
and $\overline{a}$ is the corresponding element in $\overline{M}$.
Recall from \S\ref{subsubsec:TPC-funct-K} that we denote by
$\mathcal{T}PC$ the category of TPCs and TPC functors.

\begin{thm}\label{thm:alg-pack1}
  There are four functors:
  $ K, K\ac, K_{\infty}, \tor K :\mathcal{T}PC\to
  \md_{\La_{P}}$ that are defined on objects by
  $$K(\C):=K(\C_{0})\ , K\mathcal{A}(\C):=K(\mathcal{A}\C_{0})
  \ , K_{\infty}(\C):=K(\C_{\infty}) \ , \tor K(\C): = [G \langle
  \ac\C \rangle \cap \ \mathcal{R}_{\C_{0}}]/\mathcal{R}_{\ac \C_{0}}
  \ $$ (see (\ref{eq:tor-2}) for \zjnote{$\tor K$}), and correspondingly on
  morphisms.  There are natural transformations relating $K$,
  $K\mathcal{A}$, $K_{\infty}$, and $\tor K$ induced by the maps $j$,
  $Q_{r}$ from \S\ref{ssec-acyc-cat}, and $\pi$ and $\psi$ from
  \S\ref{subsec:infty-lev} and the four functors fit into exact
  sequences as in (\ref{eq:exact-seq-k}):
  \begin{equation} 0\to\tor K \stackrel{\psi}{\longrightarrow}
    K\A\stackrel{j}{\longrightarrow} K \stackrel{\pi}{\longrightarrow}
    K_{\infty}\to 0.
  \end{equation}
The four functors $\tor K$, $K\ac$, $K$, $K_{\infty}$ commute with
  the ${\rm op}$ functors in the domain and target.  Moreover, when
  restricted to the subcategory of $\mathcal{T}PC$ consisting of
  monoidal TPCs and functors, they take values in the category of
  algebras over $\La_{P}$.
\end{thm}


\section{Filtered chain complexes} \label{sec-fil-cc}

In this section, we will carry on explicit computations of the
$K$-group when $\C = H_0(\fch)$, the homotopy category of filtered
chain complexes over $\k$, and make explicit several properties from
\S \ref{sec-k-TPC} in this setting.

\subsection{Main computational result} \label{ssec-conv} We refer to
~\cite{UZ16} for the general background on the algebraic structures
related to this category.  We recall below some notation, definitions
and results.

We denote by $\fch$ the category of filtered chain complexes. 
As explained in~\cite[Section~6.2]{BCZ21}, $\fch$ has the structure
of a filtered dg-category which is pre-triangulated, therefore its
degree-0 homology category  denoted by
$H_0(\fch)$ is a TPC (see~\cite[Proposition~6.1.12]{BCZ21}). 

\medskip

Given a filtered chain complex $(C_*, \partial, \ell)$ there exists a
decomposition, unique up to filtered chain isomorphism \pbnote{(see
  the work of Barannikov~\cite{Bar94}; See also Proposition~7.4
  in~\cite{UZ16})},
\begin{equation} \label{fcc-decomp} (C_*, \partial, \ell) \cong
  \bigoplus_{k \in \Z}C_{(k)}[-k] \,\,\,\,\mbox{where}\,\, C_{(k)} =
  \bigoplus_{i=1}^{r_k} E_1(a_i) \oplus \bigoplus_{i=1}^{s_k} E_2(b_i,
  c_i).
\end{equation}
Here, $E_1(a)$ and $E_2(b,c)$ denote the following two elementary
filtered chain complexes serving as the building blocks of the
decomposition:
\begin{itemize}
\item[-]
  $E_1(a) = (( \cdots \to 0 \to \k x \to 0 \to \cdots),
  \ell(x)= a)$, $\deg(x)=0$;
\item[-]
  $E_2(b,c) = (( \cdots \to 0 \to \k y \xrightarrow{\partial} \k x \to
  0 \to \cdots), \ell(x)= b, \ell(y) = c)$, $\deg(x)=0$, $\deg(y)=1$.
\end{itemize}

\medskip

In what follows we will concentrate on two  full
subcategories of $\fch$ (see also \S\ref{subsec:pres-mod-pre}).
\begin{itemize} \label{p:fg-t}
\item[(fg)] (Finitely generated chain complexes). The full subcategory
  $\fchfg \subset \fch$ whose objects are filtered chain complexes
  $(C_*, \partial, \ell)$ for which the total space
  $\oplus_{k \in \mathbb{Z}}C_k$ is finite dimensional over $\k$.
\item[(t)] (Tame chain complexes). The full subcategory
  $\fcht \subset \fch$ whose objects are filtered chain complexes
  $(C_*, \partial, \ell)$ such that for every filtration level
  $r \in \mathbb{R}$, the truncated chain complex
  $(C_*^{\leq r}, \partial)$ is finite dimensional over $\k$ (again,
  when taking all degrees together) and, for sufficiently small $r$, we 
  have $C_{\ast}^{\leq r}=0$.
\end{itemize}
As we will see below, there are significant differences between the
computational results of $K(H_0(\fchfg))$ and $K(H_0(\fcht))$.

\medskip

The TPC $H_0(\fch)$ is monoidal in the sense of Definition
\ref{dfn-tensor}.  The bi-operator $\otimes$ is defined by the
classical tensor product of two chain complexes, that is, for
$C = (C_*, \partial_C, \ell_C)$ and $D = (D_*, \partial_D, \ell_D)$,
we have
\[ C \otimes D = \left( (C \otimes D)_*, \partial_{C \otimes D},
    \ell_{C\otimes D}\right) \] where
$(C \otimes D)_k = \bigoplus_{p+q = k} C_p \otimes D_q$, and
$(\partial_{C \otimes D})_k = \sum_{p+q=k} \left((\partial_C)_p
  \otimes \mathds{1}_D + (-1)^p (\mathds{1}_C \otimes
  (\partial_D)_q)\right)$. Moreover,
$\ell_{C \otimes D} = \ell_C + \ell_D$. Then, clearly this $\otimes$
is associative and symmetric. The unity $e$ is given by the filtered
chain complex,
\[ E_1(0) = \left( \cdots \to 0 \to \k x \to 0 \to \cdots, \,\,\,
    \partial = 0, \,\,\,\deg(x) = 0 \,\,\, \ell(x) =0 \right). \] The
fact that $e$ is a unit follows from
Example~\ref{ex-tensor-elementary} below. Most of the other axioms in
Definition~\ref{dfn-tensor} are readily verified. For instance, the
family of isomorphisms $\{\theta_r\}_{r \in \R}$ in the axiom (iii)
are in fact all identities. Indeed, by definition,
$\Sigma^r (C_*, \partial, \ell) = (C_*, \partial_C, \ell_C+r)$, so
\[ (\Sigma^r C) \otimes D = \Sigma^r(C \otimes D) = \left( (C \otimes
    D)_*, \partial_{C \otimes D}, \ell_{C\otimes D} + r \right).\] The
only non-trivial part is the verification that $(-) \otimes C_*$ is a
triangulated functor on $H_0(\fch)_0$. This follows from Lemma 15.58.2
in~\cite{Stacks}.

\begin{ex} [Example 2.8 in \cite{PSS17} or Proposition 9.2 in
  \cite{Zha19}] \label{ex-tensor-elementary} For later use, this
  example computes the tensor product of elementary filtered complexes
  that appear in the decomposition (\ref{fcc-decomp}). There are 
  three cases:
    \begin{itemize}
  \item[-] $E_1(a) \otimes E_1(b) = E_1(a+b)$.
  \item[-] $E_1(a) \otimes E_2(b,c) = E_2(a+b, a+c)$.
  \item[-]
    $E_2(a,b) \otimes E_2(c,d) = E_2(a+c, \min\{a+d, b+c\}) \oplus
    E_2(\max\{a+d, b+c\}, b+d)[-1]$.
  \end{itemize}
  We emphasize that in the last case, there is a degree shift in
  the second term.
\end{ex}

It is straightforwards to verify that the above tensor structure
restricts to tensor structures  on  $H_0(\fchfg)$ and
$H_0(\fcht)$. We denote $\C^{\ta}=H_0(\fcht)$ and $\C^{\fg}=H_0(\fchfg)$.

\medskip

The various computations of $K$-groups for these two categories are
collected in the next theorem. Recall the Novikov rings $\La$ and
$\La_{P}$ from~\S\ref{subsubsec:nov}.

\begin{thm} \label{t:computations-chains} The tensor product of
  filtered chain complexes induces ring structures on
  $\KK(\mathcal{C}^{\rm t})$ and $\KK(\mathcal{C}^{\rm fg})$. Moreover:
  \begin{enumerate}
  \item there exists a canonical isomorphism of rings
    \begin{equation} \label{eq:lambda-iso-1} \lambda:
      \KK(\mathcal{C}^{\rm t}) \longrightarrow \Lambda,
    \end{equation}
    uniquely defined by the property that $\lambda ([E_{1}(a)]) = t^a$
    for every $a \in \mathbb{R}$, where $E_{1}(a)$ is as above and
    $[-] \in \KK(\mathcal{C}^{\rm t})$ is the $\KK$-class of the
    relevant filtered chain complex.
  \item Denote by
    $i^{\text{\rm fg}} : \KK (\mathcal{C}^{\rm fg}) \longrightarrow
    \KK(\mathcal{C}^{\rm t})$ the map induced by the inclusion
    $\mathcal{C}^{\rm fg}_0 \subset \mathcal{C}^{\rm t}_0$. Then the
    composition
    $\lambda \circ i^{\text{\rm fg}} : \KK(\mathcal{C}^{\rm fg})
    \longrightarrow \Lambda$ has values in $\Lambda_P$ and is an
    isomorphism of rings $\KK(\mathcal{C}^{\rm fg}) \cong \Lambda_P$.
  \item Denote by $\ac \mathcal{C}^{\rm fg}$ the homotopy category of
    acyclic finitely generated filtered chain complexes and by
    $j_K: \KK(\ac \mathcal{C}^{\rm fg}) \longrightarrow
    \KK(\mathcal{C}^{\rm fg})$ the map induced by the inclusion
    $\ac \mathcal{C}^{\rm fg}_{0} \subset \mathcal{C}^{\rm fg}_0$.  Then the
    composition $\lambda \circ j_K$ has values in the ideal
    $\laideal = \{P(t) \in \Lambda_P \mid P(1) = 0\} \subset
    \Lambda_P$ and gives an isomorphism
    $\KK(\ac \mathcal{C}^{\rm fg}) \cong \laideal$ which is compatible
    with $\KK(\ac \mathcal{C}^{\rm fg})$ being a
    $\KK(\mathcal{C}^{\rm fg})$-module and $\laideal$ being a
    $\Lambda_P$-module.
  \end{enumerate}
\end{thm}

\begin{remark}\label{rem:various-la}
  \begin{enumerate}
  \item An older result of Kashiwara~\cite{Ka:index} establishes a
    similar isomorphism to the ones described at points~(1) and~(2) of
    Theorem~\ref{t:computations-chains} but in the framework of the
    derived category of constructible sheaves on manifolds. The
    analogy to the result of Kashiwara becomes more transparent in
    view of the identification of Novikov rings with rings of step
    functions as described in~\S\ref{sb:diff-Nov}.

    Another related result appears in~\cite{BOO:multi-param}, where
    the $K$-group of the category of multi-parameter persistence
    modules is calculated. Note that the framework in that paper is
    slightly different than ours - the category of persistence modules
    is {\em exact} while the homotopy categories $H_0(\fchfg)$,
    $H_0(\fcht)$ of filtered chain complexes are TPC's.
  \item The morphism $\la:K(\C^{\ta})\to \La$ fits with the map
    $\la_{(-)}$ defined in \S\ref{subsec:barfuc} for persistence
    modules and their barcodes. The relation is that
    $\la (C) = \la_{H(C)}$ for each filtered chain complex $C$, where
    $H(C)$ is the persistence homology of $C$.
  \end{enumerate}
\end{remark}

The proof of the theorem will occupy the next two
subsections. \zjnote{Explicitly, \S \ref{sb:K-fg}} is concerned with $\C^{\fg}$ and
\S \ref{ssec-(T)} with $\C^{\ta}$.  The combination of the results in
these two subsections cover all the claims made in Theorem
\ref{t:computations-chains}. In \S\ref{subsubsec:other-rel} we discuss
some other relations involving the acyclics for these categories. They
are relevant to understand some $K$-theoretic aspects, as in Remark
\ref{rem:var-K} \zjred{(a)}.

\subsection{Computation of the $K$-groups associated to
  $H_0(\fchfg)$} \label{sb:K-fg} Throughout this subsection we write
$\mathcal{C}^{\fg} := H_0(\fchfg)$ for the homotopy category of
finitely generated filtered chain complexes over $\k$. By a slight
abuse of notation we will denote the $K$-class of a filtered chain
complex $A$ sometimes by $[A]$ but occasionally also simply by $A$.

\begin{lemma} \label{lem-difference} In $K(\mathcal{C}^{\rm fg})$, we have
  $E_2(a,b) = E_1(a) - E_1(b)$ for any $a \leq b \in \R$.
\end{lemma}

\begin{proof} The mapping cone of the (filtration preserving) map
  $\phi: E_1(b) \to E_1(a)$ defined by $\phi(x_b) = x_a$ is
  $E_2(a,b)$, where $x_b$ and $x_a$ are the generators with the
  prescribed filtrations of $E_1(b)$ and $E_1(a)$ respectively. Thus,
  by the definition of the $K$-group, we have
  $E_1(a) = E_1(b) + E_2(a,b)$ in $K(\mathcal{C}^{\fg})$.
\end{proof}

\zjnote{Recall that $\Lambda_{P}$ denotes the ring of Novikov polynomials} with
coefficients in $\Z$, namely {\em finite} sums of monomials of the
form $nt^a$ with $n \in \mathbb{Z}$, $a \in \mathbb{R}$ and where $t$
is a formal variable. More precisely:
\begin{equation}\label{nov-poly}
  \Lambda_{P}:=
  \left\{ \sum_{k=0}^N n_k t^{a_k} \, \bigg| \, N \in \mathbb{Z}_{\geq 0},
    \, n_k \in \Z, \,\,a_k \in \R\right\}.
\end{equation}
We endow $\Lambda_P$
with the obvious ring structure whose unit is $1 = 1 \cdot t^0$. As explained at the end of~\S\ref{ssec-conv}, since $\C^{\fg}$ is
monoidal, Lemma~\ref{lem-3} implies that the $K$-group $K(\C^{\fg})$
admits a well-defined ring structure, where the product is given by
\begin{equation} \label{K-htp-product} E_1(a) \cdot E_1(b) = E_1(a+b),
\end{equation}
due to~\eqref{dfn-ring} and Example~\ref{ex-tensor-elementary}. Note
that by Lemma~\ref{lem-difference} we only need to consider classes in
$K(\C^{\fg})$ represented by the elementary filtered chain complexes
of $E_1$-type. Moreover, the unit of $K(\C^{\fg})$ is $E_1(0)$.

\medskip

\begin{prop} \label{prop-kc0} There is a ring isomorphism
  \pbnote{$\eta: \Lambda_{P} \xrightarrow{\; \; \cong \; \; }
    K(\C^{\rm fg})$} that sends \zjred{$t^a \in \Lambda_P$ to the
    $K$-class $[E_1(a)] \in K(\mathcal{C})$ for every
    $a \in \mathbb{R}$}.
\end{prop}

\begin{proof} According to~\eqref{fcc-decomp}, any object in $\C^{\fg}_0$ is
  isomorphic to a chain complex of the following form
  \begin{equation} \label{eq:can-form-ch}
    C=\bigoplus_{k = k_0}^{k_1} C_{(k)}[-k] \,\,\,\,\mbox{where}\,\,
    C_{(k)} = \bigoplus_{i=1}^{r_k} E_1(a_i) \oplus
    \bigoplus_{i=1}^{s_k} E_2(b_i, c_i).
  \end{equation}
  Since the objects of $\mathcal{C}^{\fg}$ are all finitely generated
  chain complexes, we have that $k_0$, $k_1$, $r_k$ and $s_k$ are all
  finite. Therefore, the class $[C] \in K(\C^{\fg})$ of such a complex
  is of the following form (recall that the degree shift $[-1]$ of a
  chain complex corresponds to multiplying by $-1$ in the $K$-group)
  \[ [C] = \sum_{k=k_0}^{k_1} (-1)^{-k} \left( \sum_{i=1}^{r_k}
      [E_1(a_i)] + \sum_{i=1}^{s_k} [E_2(b_i, c_i)] \right).\] By
  Lemma~\ref{lem-difference}, and by grouping terms, the expression
  above can be further simplified to:
  \begin{equation} \label{express-g2} [C] = \sum_{i=1}^{k} n_i
    [E_1(x_i)] \,\,\,\mbox{for some $x_i \in \R$ and $n_i \in \Z$},
  \end{equation} 
  with $k$ finite and such that each $x_{i}$ only appears once in the
  sum. We may also assume that the $x_{i}$'s are written in increasing
  filtration order.

  In summary, each element in $K(\C^{\fg})$ can be written in the
  form~\eqref{express-g2}.  We claim that the representation of $[C]$
  as written in~\eqref{express-g2} is unique. This is equivalent to
  showing that if for a filtered complex $C$ we have
  $[C]=0\in K(\C^{\fg})$, then the expression~\eqref{express-g2}
  associated to $C$ does not contain any terms on the right hand side.
  The argument is based on the filtered Euler characteristics defined
  by
  \begin{equation}\label{eq:Euler-ch0}
    \chi_{\alpha}(C)=\chi(C^{\leq \alpha}) \ ,  \ \  \alpha \in \R~.~
  \end{equation}
  These Euler characteristics $\chi_{\alpha}$ descend to morphisms
  $\chi_{\alpha}: K(\C^{\fg})\to \Z$ due to the standard behavior of
  the Euler characteristic on short exact sequences of
  complexes. Additionally, it is immediate to see that
  $$\chi_{\alpha}\Bigl(\sum_{i=1}^{k} n_i
  [E_1(x_i)] \Bigr)=0 \ , \ \forall \alpha \in \R, \ \ \text{iff} \ \
  n_{i}=0 \ \ \forall \ i, $$ which shows our claim.

  What we have proved so far apriori holds only for $K$ classes of the
  form $[C]$ with $C \in \Ob(\mathcal{C}^{\fg})$. However, every
  element $X \in K(\mathcal{C}^{\fg})$ is of that form, i.e.~$X=[C]$
  for some finitely generated filtered chain complex $C$. This is
  because $\Ob(\mathcal{C}^{\fg})$ is closed under direct sums (and
  shifts in grading). Therefore the preceding statements in fact hold
  for all elements of $K(\mathcal{C}^{\fg})$. In the rest of the proof
  we continue to represent every element of $K(\mathcal{C}^{\fg})$ by
  an object of $\mathcal{C}^{\fg}$.
  
  Consider the map
  \begin{equation} \label{eta-map} \eta: \Lambda_{P} \to
    K(\C^{\fg}) \,\,\,\,\,\mbox{defined by} \,\,\,\,\, t^a \mapsto
    E_1(a)
  \end{equation}
  and extend it linearly (over $\mathbb{Z}$) to $\Lambda_{P}$. We
  claim that $\eta$ is surjective. Indeed, for any $g \in K(\C^{\fg})$
  written as on the right side of (\ref{express-g2}), consider the
  Novikov polynomial $p = \sum_{i=1}^{k} n_i t^{x_i}$. We then have
  $\eta(p) = g$. Next, we claim that $\eta$ is injective. This is an
  immediate consequence of the fact that the writing
  (\ref{express-g2}) of an element in $K(\C^{\fg})$ is unique.
  
  Finally, $\eta$ is a ring homomorphism. Indeed,
  $\eta(1) = \eta(t^0) = E_1(0) = e$, the unit, and
  \begin{align*}
    & \eta\left(\sum_{k} n_k t^{a_k} \right) \cdot
      \eta\left( \sum_{k} m_k t^{b_k}\right)
      = \left(\sum_{k} n_k E_1(a_k)\right) \cdot
      \left(\sum_{k} m_k E_1(b_k)\right) \\
    & = \sum_{k} \sum_{i+j = k} n_i m_j E_1(a_i +b_j)
    = \eta \left(\sum_{k} \sum_{i+j = k} n_i m_j t^{a_i + b_j}
    \right)
    = \eta \left(\left(\sum_{k} n_k t^{a_k}\right) \cdot
    \left(\sum_{k} m_k t^{b_k}\right) \right),
  \end{align*}
  where the second equality comes from the product structure on
  $K(\C^{\fg})$ defined in (\ref{K-htp-product}).
\end{proof}

We denote by
 \begin{equation}\label{eq:lambda-def}
 \lambda: = \eta^{-1}: K(\C^{\fg}) \to \Lambda_{P}
\end{equation} the inverse of the isomorphism $\eta$
from  Proposition \ref{prop-kc0}. We have
$\lambda (E_{1}(a))=t^{a}$  for all $a\in \R$.

 We now turn to the category of acyclics $\ac \mathcal{C^{\fg}}$ and
 its associated $K$-group $K\ac(\C^{\fg})=K(\ac \C^{\fg}_0)$ that was
 introduced in~\S\ref{ssec-acyc-cat}. Note that, unlike for
 $K(\C^{\fg})$, the generators of $K\ac(\C^{\fg})$ are all of the form
 of $E_2(a,b)$ and the relation in Lemma~\ref{lem-difference} does not
 hold anymore in $K\ac(\C^{\fg})$. Still, we have the following
 analogous useful result.

\begin{lemma} \label{lem-E_2-sum} For any $a \leq b \leq c \in \R$, we
  have $E_2(a,b) + E_2(b,c) = E_2(a,c)$ in
  $K\ac (\C^{\rm fg})$.
\end{lemma}

\begin{proof} Suppose that the generators of $E_2(a,c)$ are $y_0$ and
  $y_1$ satisfying $\partial(y_1) = y_0$ and filtrations
  $\ell(y_0)= a$, $\ell(y_1) = c$. Similarly, suppose that the
  generators of $E_2(b,c)$ are $x_0$ and $x_1$ satisfying
  $\partial(x_1)= x_0$ and filtrations $\ell(x_0) = b$,
  $\ell(x_1) = c$. Then consider a filtration preserving map
  $\phi: E_2(b,c) \to E_2(a,c)$ defined by
  \[ \phi(x_0) = y_0 \,\,\,\,\mbox{and}\,\,\,\, \phi(x_1) = y_1. \] 
 \zjred{We} have
  ${\rm Cone}(\phi) = E_2(a,b)$. In other words, we have an exact
  triangle
  \[ E_2(b,c) \xrightarrow{\phi} E_2(a,c) \to E_2(a,b) \to
    E_2(b,c)[-1] \]  and thus
  $E_2(a,c) = E_2(a,b) + E_2(b,c)$ in $K\ac(\C^{\fg})$.
\end{proof}

Similarly to (\ref{nov-poly}), we introduce the following (non-unital)
ring $\Lambda'_{\text P}$ of ``double-exponent'' Novikov polynomials
(in the formal variable $s$). Additively it is defined as:
\begin{equation} \label{double-exp} \Lambda'_{P} : = \left\{
    \sum_k n_k s^{a_k, b_k} \, \bigg| \, n_k \in \Z, \,\, a_k< b_k <
    \infty \right\}\bigg/ \mathcal R \, , 
\end{equation}
where the sums above are assumed to be finite. The relation subgroup
$\mathcal R$ is generated by the elements
$\{s^{a,b} + s^{b,c} - s^{a,c}\}$ for any $a < b < c \in \R$. The
product on \zjnote{$\Lambda'_{P}$} goes as follows,
\begin{equation} \label{double-Novikov-product} s^{a,b} \cdot s^{c,d}
  : = s^{a+c, a+d} - s^{b+c, b+d} = s^{a+c, b+c} - s^{a+d, b+d}
\end{equation}
and extended by linearity over $\mathbb{Z}$. Note that this is
well-defined and the second equality comes from the relation
$\mathcal R$. It is easy to verify that the product defined
in~\eqref{double-Novikov-product} is associative and distributive.

The tensor structure on $\C^{\fg}$ descends to a well-defined
(non-unital) tensor structure on $\ac \C^{\fg}$, therefore,
Lemma~\ref{lem-3} implies that the $K$-group $K\ac (\C^{\fg})$ admits a
well-defined (non-unital) ring structure, where the product is given
by
\begin{align*} \label{ac-product} E_2(a,b) \cdot E_2(c,d)
  & =  E_2(a+c, \min\{a+d, b+c\})- E_2(\max\{a+d, b+c\}, b+d)\\
  &= E_2(a+c, a+d) - E_2(b+c, b+d).
\end{align*}
Here the first equality comes from Example~\ref{ex-tensor-elementary}
and the second equality comes from Lemma~\ref{lem-E_2-sum}. Indeed, if
$b+c < a+d$, then
\begin{align*}
  E_2(a+c, b+c) - E_2(a+d, b+d)
  & = \left(E_2(a+c, b+c) + E_2(b+c, a+d)\right) \\
  & \,\,\,\,\,- \left(E_2(b+c, a+d) + E_2(a+d, b+d)\right) \\
  & = E_2(a+c, a+d) - E_2(b+c, b+d).
\end{align*}
Similarly, one can also write
\begin{equation} \label{product-use} E_2(a,b) \cdot E_2(c,d) =
  E_2(a+c, b+c) - E_2(a+d, b+d).
\end{equation}
The expression in (\ref{product-use}) explains the relations in
$\mathcal R$ listed above.

\begin{prop} \label{prop-kac0-2} There is a {\em multiplicative}
  isomorphism $K\ac(\mathcal{C}^{fg}) \cong \Lambda_{P}'$ that
  sends the $K$-class $E_2(a,b) \in K\ac (\mathcal{C}^{fg})$ to
  $s^{a,b} \in \Lambda_{\text P}'$, for every $a < b$.
\end{prop}
For the proof of Proposition \ref{prop-kac0-2} it is useful to relate
\zjnote{$\Lambda_{P}$ to $\Lambda'_{P}$}.  Consider the following ideal
of $\Lambda_{\rm P}$,
\begin{equation} \label{euler-0-poly} \Lambda_{P}^{(0)} : =
  \left\{p(t) \in \Lambda_{P} \,| \, p(1) = 0\right\}.
\end{equation}
The condition $p(1) = 0$ is equivalent to the condition that the sum
of all coefficients of the Novikov polynomial $p$ is equal to $0$.

\begin{lemma} \label{lem-lambda-2} Consider the map
  \zjnote{  \begin{equation} \label{sigma-map}  \sigma: \Lambda'_{P} \to
\Lambda_{P}^{(0)} \,\,\,\,\,\mbox{defined by}
    \,\,\,\,\,s^{a,b} \mapsto t^a - t^b
  \end{equation}}
  and extend it to \zjnote{$\Lambda'_{P}$} by linearity over
  $\mathbb{Z}$. Then $\sigma$ is well-defined and is a non-unital ring
  isomorphism.
\end{lemma}

The proof of the Lemma can be found in the expanded version of the
paper~\cite{BCZ-pkt-arxiv}. We are now ready to prove Proposition
\ref{prop-kac0-2}.

\begin{proof} [Proof of Proposition \ref{prop-kac0-2}] Consider the
  following ideal of $K(\C^{\fg})$,
  \begin{equation} \label{euler-0} K(\C^{\fg})^{(0)} : = \left\{A \in
      K(\C^{\fg}) \, | \, \chi(A) = 0 \right\}
  \end{equation}
  where $\chi: K(\C^{\fg}) \to \Z$ is the Euler characteristic. Note
  that $\chi$ is well defined for all objects of $\mathcal{C}^{\fg}$
  since they are assumed to be finitely generated.  By the same
  argument as in Proposition~\ref{prop-kc0}, the ring isomorphism
  $\eta: \Lambda_{P} \to K(\C^{\fg})$ restricts to a ring isomorphism
  $\eta|_{\Lambda^{(0)}_{P}}: \Lambda^{(0)}_{P} \to
  K(\C^{\fg})^{(0)}$. Then consider the following diagram,
  \begin{equation} \label{lam-K-com-2} \xymatrix{ \Lambda'_{P}
      \ar[r]^-{\sigma}_-{\cong} \ar[d]_-{\eta'} & \Lambda^{(0)}_{
        P}
      \ar[d]_-{\cong}^-{\eta|_{\Lambda^{(0)}_{P}}}\\
      K\ac(\mathcal{C}^{\fg}) \ar[r]^-{j} & K(\mathcal C^{\fg})^{(0)}}
  \end{equation}
  where the top horizontal map $\sigma$ is a ring isomorphism by Lemma
  \ref{lem-lambda-2}. Here,
  \begin{equation} \label{ac-eta} \eta': \Lambda'_{P} \to
    K\ac(\mathcal{C}^{\fg})\,\,\,\,\,\mbox{is defined by} \,\,\,\,\,
    \sum_{k} n_k s^{a_k, b_k} \mapsto \bigoplus_{k} n_k E_2(a_k, b_k).
  \end{equation}
  Lemma \ref{lem-E_2-sum} and the relations in $\mathcal R$ of
  $\Lambda'_{\rm P}$ in (\ref{double-exp}) yields that $\eta'$ is
  well-defined. Also, it is readily checked that $\eta'$ is
  surjective. Meanwhile, $j$ is the ring homomorphism defined in
  (\ref{j-map}). We claim that the diagram (\ref{lam-K-com-2})
  commutes. Indeed, for any $s^{a,b}$, we have
  \begin{equation} \label{com-j} \xymatrix{
      s^{a,b} \ar[r]^-{\sigma}  \ar[d]_-{\eta'} & t^a - t^b \ar[d]^-{\eta}\\
      E_2(a,b) \ar[r]^-{j} & E_1(a) - E_1(b)}
  \end{equation}
  which verifies the commutativity. Since
  $j \circ \eta' = \eta|_{\Lambda^{(0)}_{P}} \circ \sigma$ which
  is an isomorphism, we deduce that $\eta'$ is injective.
\end{proof}

\subsection{Computation of the $K$-groups associated to
  $H_0(\fcht)$} \label{ssec-(T)} We are now interested in the
triangulated persistence category $\C^{\ta}:=H_0(\fcht)$ of tamed
filtered chain complexes over $\k$ (see page~\pageref{p:fg-t}) and its
$K$-theory.

The computations of $K(\C^{\ta})$ and $K\ac(\C^{\ta})$ partially overlap with
the ones for $H_0(\fchfg)$, but they differ at some points. Instead of
\zjnote{$\Lambda_{P}$ and $\Lambda'_{P}$}, we will now need to work
with the full Novikov ring 
\begin{equation} \label{nov-ser} \Lambda := \left\{
    \sum_{k=0}^{\infty} n_k t^{a_k} \, \bigg| \, n_k \in \Z, \,\,a_k
    \to \infty\right\}
\end{equation}
as well as with
\begin{equation} \label{nov-ser-double} \Lambda' : = \left\{
    \sum_{k=0}^{\infty} n_k s^{a_k, b_k} \,\bigg| \, n_k \in \Z,
    \,\,a_k < b_k < \infty, \,\, a_k \to \infty\right\}\bigg/
  \mathcal{R}
\end{equation} 
where the relation subgroup $\mathcal R$ is defined as in~\eqref{double-exp}.

\medskip

Here is the analogue of Proposition~\ref{prop-kc0} in this setting.
\begin{prop} \label{prop-kc0-T} There is a ring
  isomorphism $K(\C^{\rm t}) \cong \Lambda$ that sends the $K$-class
  $E_1(a) \in K(\mathcal{C}^{\rm t})$ to $t^a \in \Lambda$ for every
  $a \in \mathbb{R}$.
\end{prop}

The proof follows from arguments almost identical to those in
Proposition~\ref{prop-kc0}. Additional care is required when comparing
two elements in $\Lambda$ that possibly involve infinitely many
monomials. We can reduce this to the polynomial situation via a
truncation at a filtration level, due to the definition of tame chain
complexes.

\medskip

We also have the following result, analogous to
Proposition~\ref{prop-kac0-2}. We omit its proof which can be found
in~\cite{BCZ-pkt-arxiv}.
\begin{prop} \label{prop-kac0} There is a multiplicative isomorphism
  $K\ac(\mathcal{C}^{\rm t}) \cong \Lambda'$ that sends the $K$-class
  $E_2(a,b) \in K\ac (\mathcal{C}^{\rm t})$ to $s^{a,b} \in \Lambda'$,
  for every $a < b$. Moreover, there is a multiplicative isomorphism
  $\La'\cong \La$.
\end{prop}

We continue to denote by $\lambda$ the inverses of the maps $\eta$ and
$\eta'$ before, as in~\eqref{eq:lambda-def}.

\begin{rem} The map
  $$j: K\ac (\C^{\ta})\to K (\C^{\ta})$$ is an isomorphism.
  This follows from diagram~(63) in~\cite{BCZ-pkt-arxiv}. It can also
  be seen in a more direct way by first noticing that the unit
  $e = E_1(0)$ of $K(\mathcal{C}^{\ta})$ belongs to
  $K\ac(\mathcal{C}^{\ta})$.  Indeed, consider the element
  $$e' = \sum_{k=0}^{\infty} E_2(k, k+1) \in K\ac (\C^{\ta}).$$  This
  element is a unit in $K\ac(\C^{\ta})$ as shown by the next
  computation.  Fix some $E_{2}(a,b)$ and write:
  \begin{align*}
    E_2(a,b) \cdot e'
    & = E_2(a,b) \cdot \left(E_2(0,1) \oplus E_2(1,2) \oplus \cdots \right) \\
    & = (E_2(a, b) - E_2(a+1, b+1)) + (E_2(a+1, b+1) - E_2(a+2, b+2)) + \cdots \\
    & = E_2(a,b)
  \end{align*}
  where the second equality comes from the product formula in
  (\ref{product-use}) and the last equality comes from successive
  cancellations with only the first term left. As $e'=e'\cdot e =e$ we
  deduce that
  $E_1(0) \in \textnormal{image\,} (K\ac (\C^{\ta}) \xrightarrow{\; j
    \;} K(\mathcal{C}^{\ta}))$. Now use the fact that the map $j$ is
  linear over $\Lambda$ to conclude that $j$ is an isomorphism.
\end{rem}

%

For further use, we summarize two direct consequences of
Proposition~\ref{prop-kac0-2} and Proposition \ref{prop-kac0} (see
also Proposition \ref{prop:exact-seq}).

\begin{cor} \label{c:K-TorK-ch} We have the ring isomorphisms
  \zjnote{$j: K\ac(\C^{\rm fg}) \to K(\C^{\rm fg})^{(0)}$ (where $K(\C^{\rm fg})^{(0)}$}
  is defined as in~\eqref{euler-0}) and
  \zjnote{$j: K\ac(\C^{\rm t}) \to K(\C^{\rm t})$}.  In particular, these maps are
  injective and thus $\tor K (\C^{\rm fg})=0$, $\tor K
  (\C^{\rm t})=0$. Further, $K_{\infty}(\C^{\rm fg})=\Z$,
  $K_{\infty}(\C^{\rm t})=0$.
\end{cor}

\begin{rem} \label{r:ses-ch-splits} Corollary~\ref{c:K-TorK-ch} says
  that in case $\mathcal{C} = \mathcal{C}^{\text{fg}}$ the exact
  sequence~\eqref{eq:exact-seq-k} becomes a short exact sequence of
  the form
  \begin{equation} \label{eq:ses-K-ch} 0\longrightarrow
    K\ac(\mathcal{C}^{\text{fg}})\stackrel{j}{\longrightarrow}
    K(\mathcal{C}^{\text{fg}})\stackrel{\pi}{\longrightarrow}
    K_{\infty}(\mathcal{C}^{\text{fg}})\longrightarrow 0.
  \end{equation}
  Identifying the first two terms, via the isomorphism $\lambda$, with
  $\Lambda_P^{(0)}$ and $\Lambda_P$ respectively, and the third term
  with $\mathbb{Z}$ via the Euler characteristic, we can identify the
  above sequence with the following one:
  \begin{equation} \label{eq:ses-K-Nov-ch} 0\longrightarrow
    \Lambda^{(0)}_P \longrightarrow \Lambda_P \longrightarrow
    \mathbb{Z},
  \end{equation}
  where the first map is the inclusion and the second map is given by
  $\Lambda_P \ni P(t) \longmapsto P(1) \in \mathbb{Z}$.  The
  sequence~\eqref{eq:ses-K-ch} obviously splits. A preferred splitting
  seems to come from the following right and left inverses to $\pi$
  and $j$ respectively:
  \begin{equation} \label{eq:split-ses-ch}
    \begin{aligned}
      & K_{\infty}(\mathcal{C}^{\text{fg}}) \longrightarrow
      K(\mathcal{C}^{\text{fg}}), \quad [C] \longmapsto \chi(C)[E_1(0)], \\
      & K(\mathcal{C}^{\text{fg}}) \longrightarrow
      K\ac(\mathcal{C}^{\text{fg}}), \quad [C] \longmapsto [C] -
      \chi(C)[E_1(0)].
    \end{aligned}
  \end{equation}
  In terms of the sequence~\eqref{eq:ses-K-Nov-ch} the first map
  in~\eqref{eq:split-ses-ch} is just the inclusion
  $\mathbb{Z} \longrightarrow \Lambda_P$ and the second map is
  $P(t) \longmapsto P(t)-P(1)$.
\end{rem}

\subsection{More relations} \label{subsubsec:other-rel} In this
subsection, we investigate the acyclic truncation maps $Q_r$ defined
in~\eqref{Q-map}.

Notice  that $Q_r$ is only a
homomorphism of abelian groups, since
\[ Q_r(E_1(a) \cdot E_1(b)) = Q_r(E_1(a+b)) = E_2(a+b, a+b+r) \] while
\begin{align*}
  Q_r(E_1(a)) \cdot Q_r(E_1(b))
  & = E_2(a,a+r) \cdot E_2(b,b+r) \\
  & = E_2(a+b, a+b+r) - E_2(a+b+r, a+b + 2r). 
\end{align*}

The proof of the following result can be found
in~\cite{BCZ-pkt-arxiv}.
\begin{prop} \label{prop:qqq} For any $r \geq 0$, the map
  $Q_r: K(\C^{\ta}) \to K\ac(\mathcal{C}^{\ta})$ is an isomorphism of
  abelian groups.
\end{prop}

\begin{proof} Consider the following diagram
  \begin{equation} \label{K-lambda-diag} \xymatrixcolsep{5pc}
    \xymatrix{ K(\C^{\ta}) \ar[r]^-{Q_r} \ar[d]_-{\lambda}^-{\cong}
      & K\ac(\mathcal{C}^{\ta}) \ar[d]_-{\lambda}^-{\cong} & \\
      \Lambda \ar[r]^-{\tilde{Q}_r} & \Lambda'
      \ar[r]^-{\sigma}_{\cong} & \Lambda}
  \end{equation}
  where the two vertical isomorphisms $\lambda$ are the inverses of
  the isomorphisms $\eta$ and $\eta'$ from Proposition \ref{prop-kc0}
  and Proposition \ref{prop-kac0} respectively, the right-bottom
  isomorphism $\sigma$ is from Proposition \ref{prop-kac0-2}, and the
  map $\tilde{Q}_r: \Lambda \to \Lambda'$ is defined by
  $t^a \mapsto s^{a, a+r}$ and extended to $\Lambda$ by linearity over
  $\mathbb{Z}$. Notice that $\tilde{Q}_r$ makes the square
  commutative. Moreover the composition
  \[ (\sigma \circ \tilde{Q}_r)(t^a) = \sigma(s^{a, a+r}) = t^a -
    t^{a+r} \] is injective, thus $\tilde{Q}_r$ is injective. Next we
  will show that $\tilde{Q}_r$ is also surjective.
  
  The surjectivity of $\tilde{Q}_r$ is proved in detail
  in~\cite{BCZ-pkt-arxiv}, and here we only outline the main idea.
  First note that $\Lambda$ is a $\Lambda'$-module with structure
  defined by $s^{a,b} \cdot t^c: = t^{a+c} - t^{b+c}$. It is
  straightforward to check that $\tilde{Q}_r$ is $\Lambda'$-linear
  with respect to this structure.
  Now, note that
  $\tilde{Q}_r(1 + t^r + t^{2r} + \cdots) = s^{0,r} + s^{r, 2r} +
  \cdots = u$, the unit in $\Lambda'$. This shows that $\tilde{Q}_r$
  is surjective. It now follows that $\tilde{Q}_r$ is an isomorphism,
  which implies that $Q_{r}$ is an isomorphism too.
\end{proof}

We now investigate the map $Q_r$ in the case of
$\mathcal{C}^{\fg}$. We have the following commutative diagram
that is analogue of diagram~\eqref{K-lambda-diag} (we keep denoting by $\la$ the inverses
of the corresponding maps $\eta$, as before).
\begin{equation} \label{K-lambda-diag-2} \xymatrixcolsep{5pc}
  \xymatrix{ K(\C^{\fg} ) \ar[r]^-{Q_r} \ar[d]_-{\lambda}^-{\cong}
    & K\ac(\mathcal{C}^{\fg} ) \ar[d]_-{\lambda}^-{\cong} & \\
    \Lambda_{P} \ar[r]^-{\tilde{Q}_r} & \Lambda'_{P}
    \ar[r]^-{\sigma}_{\cong} & \Lambda^{(0)}_{P}}
\end{equation}
Again, since $\sigma \circ \tilde{Q}_r$ is injective, $\tilde{Q}_r$
and $Q_r$ are injective. However, now $Q_r$ is not necessarily
surjective, again due to the lack of a unit in $K\ac(\mathcal{C}^{\fg} )$
(and also in $\Lambda'_{P}$). 

\begin{prop} \label{prop-frac-Euler} For any $r>0$, we have
  \[ {\rm Im}(\sigma \circ \tilde{Q}_r) = \left\{\sum_{k} n_k
      t^{a_k} \in \Lambda^{(0)}_{P} \,\bigg|\,
      \mbox{$\sum_{j \in J_c}n_j = 0$ for any $c \in \R/r \Z$}
    \right\} \] where
  $J_c = \{k \in \Z \, |\, a_k \equiv c \, \mbox{{\rm (mod $r$)}}\}$.
\end{prop}


\begin{remark} The ring $\Lambda^{(0)}_{P}$ can be regarded as the
  set of elements satisfying the condition in Proposition
  \ref{prop-frac-Euler} for $r = \infty$, where
  $a_k \equiv c \, \mbox{(mod $r$)}$ simply means
  $a_k=c$. \end{remark}

\begin{proof} [Proof of Proposition \ref{prop-frac-Euler}]
  \pbnote{Denote by $\Lambda_P^{(0), r}$ the ring obtained from
    $\Lambda^{(0)}_{P}$ by reducing \mbox{mod $r$} the exponents that
    appear in the variable $t$. (Namely, we take the elements
    $\sum_{k} n_k t^{a_k} \in \Lambda^{(0)}_{P}$ and reduce \mbox{mod
      $r$} all the exponents $a_k$.) The outcome $\Lambda_P^{(0), r}$
    inherits from $\Lambda^{(0)}_{P}$ the structure of a non-unital
    ring (recall that $\Lambda^{(0)}_{P}$ is itself non-unital). We
    have a projection
    $\pi: \Lambda^{(0)}_{P} \longrightarrow \Lambda_P^{(0), r}$}
  defined by
  \[ \sum_k n_k t^{a_k} \longmapsto \sum_k n_k t^{[a_k]}
    \,\,\,\,\mbox{where $[a_k]$ is the projection of $a_k$ on
      $\R/r\Z$}. \] We need to show that
  ${\rm Im}(\sigma \circ \tilde{Q}_r) = \ker(\pi)$.  If
  $p \in {\rm Im}(\sigma \circ \tilde{Q}_r)$, then each term $t^{a}$
  of $p$ is mapped by $\sigma \circ \tilde{Q}_r$ to $t^{a} -
  t^{a+r}$. Applying the map $\pi$, we have
  $t^{a} - t^{a+r} \mapsto t^{[a]} - t^{[a]} = 0$. We get $\pi(p) =0$.
  This shows ${\rm Im}(\sigma \circ \tilde{Q}_r) \subset \ker(\pi)$.

  Fix $p =\sum_{n_k} n_k t^{a_k} \in \ker(\pi)$ and consider the
  sub-polynomial $p_c: = \sum_{n_k} n_k t^{a_k}$ with $[a_k] = c$ for
  each $c \in \R/r \Z$ . We will show that any sub-polynomial $p_c$ is
  in the image of the map $\sigma \circ \tilde{Q}_r$ - this is enough
  to show that $p$ is in this image.  Write
  \begin{equation} \label{pc} p_c = n_1 t^{c+m_1r} + \cdots + n_L
    t^{c+m_L r} \,\,\,\,\mbox{where
      $m_1 \leq \cdots \leq m_L \in \Z$},
  \end{equation}
  and $\sum_{k=1}^L n_k = 0$. At this point the argument is similar to
  the proof of Lemma \ref{lem-lambda-2}, by induction on
  $L$. The first non-trivial case is when $L = 2$, that is
  $p_c = n t^{c+m_1r} - n t^{c+m_2 r}$. Then consider
  $x = n (t^{c+m_1r} + t^{c + (m_1+1) r} + \cdots + t^{c + (m_2-1)r})
  \in \Lambda_{\rm P}$, and we have
  \begin{align*}
    (\sigma \circ \tilde{Q}_r)(x)
    & = n (t^{c+m_1r} -t^{c+(m_1+1)r} +t^{c+(m_1+1)r} - t^{c+(m_1+2)r} \\
    \,\,\,\,\,\,\, & + \cdots +  t^{c + (m_2-1)r} - t^{c + m_2r}) = p_c. 
  \end{align*}
  Now, assume that any $p_c$ in (\ref{pc}) with $L$-many non-zero
  terms is in the image. We rewrite a $p_c$ with $(L+1)$-many non-zero terms,
   as follows,
  \begin{align*}
    p_c & = n_1 t^{c+m_1r}  + \cdots + n_L t^{c+m_L r} +
          \left(-\sum_{k=1}^L n_k\right) t^{c + m_{L+1}r} \\
        & = \left(n_1 t^{c+m_1 r} + \cdots n_{L-1} t^{c+m_{L-1} r} -
          \left(\sum_{k=1}^{L-1} n_k\right) t^{c+m_L r}\right)  \\
    \,\,\,\,\, &+\left(\sum_{k=1}^L n_k\right) \left(t^{c+m_L r} -
                 t^{c + m_{L+1}r}\right).
  \end{align*}
  Each of the two terms above is in the image by our inductive
  hypothesis. Therefore, $p_c$ is in the image. \end{proof}

\begin{remark}\label{rem:var-K} (a) Let $\C$ be a TPC and fix
  $r\in [0,\infty)$.  Recall the notations from
  \S\ref{subsec:infty-lev}. In particular, denote by
  $G \langle \C \rangle$ the free abelian group generated by the
  objects of $\C$. A natural question is whether there are $K$-type
  groups of interest obtained by dividing $G\langle \C \rangle$ by the
  subgroup $\mathcal{R}_{r}$ generated by the relations of the form
  $B=A-C$ for each \zjred{strict exact} triangle
  $$A\to B\to C\to \Sigma^{-s}TA$$ \zjred{of weight
  $s\leq r$}. In general, this does not lead to something useful.  For
  instance, it is easy to see that if acyclics of any weight belong to
  the subgroup of $K(\C)$ generated by $s$-acyclics, then
  $G \langle \C \rangle / \mathcal{R}_{r}=\C_{\infty}$.  This \zjred{is}
  what happens in $\C^{\ta}$, as it follows from Proposition
  \ref{prop:qqq}.
 
 (b) Another natural way to map $K(\C^{\fg} )$ to $K\ac(\C^{\fg} )$ is
  as explained in Remark~\ref{r:ses-ch-splits}. Namely:
  $$K(\C^{\fg} ) \longrightarrow  K\ac(\C^{\fg} ) \ \ , \ \  [C]
  \longmapsto [C] - \chi(C) E_1(0).$$ This definition parallels (via the
  isomorphisms defined by $\lambda$) the map
  $\Lambda_P \longrightarrow \Lambda_P^{(0)}$, defined by
  $P(t) \longmapsto P(t)-P(1)$.
\end{remark}


\section{Barcodes, step functions and $K$-theory for persistence
  modules} \label{sb:diff-Nov}

It is natural to associate to an ungraded persistence module
$\{M^{\alpha}\}_{\alpha \in \mathbb{R}}$ the real function given by
its dimension $\alpha \longmapsto \dim_{\k}( M^{\alpha})$. By taking
into account also graded modules - as in \S\ref{subsubsec:iso-rings} -
this provides a way to associate to each persistence module $M$, or
graded barcode $\mathcal{B}$, a function $\bar{\sigma}_{\mathcal{B}}$
belonging to a ring $LC_{B}[x,x^{-1}]$ where $x$ is a formal variable
keeping track of the grading and $LC_{B}$ is a ring of bounded step
functions $\R\to \R$ with a convolution type product.  The value of
this representation for $x=-1$, $\bar{\chi}=\bar{\sigma}|_{-1}$, is
interpreted in \S\ref{subsubsec:Novikov-etc} as the application that
associates to a persistence module its \zjnote{$K$-class}, in an
appropriately defined \zjnote{$K$-group}. This is a consequence of the fact
that $LC_{B}$ is isomorphic as a ring with the universal Novikov
polynomial ring, $\La_{P}$, and fits with the description
in~\S\ref{sec-fil-cc} of the $K$-group of the homotopy category of
filtered chain complexes. Finally, in~\S\ref{subsubsec:morse} we use
the theory developed in~\S\ref{subsubsec:iso-rings}
and~\S\ref{subsubsec:Novikov-etc} to formulate a general form of Morse
inequalities for persistence homology.

Before we go on we remark that some of the material presented
in~\S\ref{subsubsec:iso-rings} and~\S\ref{subsubsec:Novikov-etc} below
is analogous to older work of Schapira~\cite{Scha:operations} done in
the framework of the category of constructible sheaves on manifolds.

\subsection{The rings $LC_{B}$ and $LC$} \label{subsubsec:iso-rings}
It is very natural to represent barcodes in a ring of functions that
we now describe.

For a subset $S \subset \mathbb{R}$ denote by
$\mathds{1}_S: \mathbb{R} \longrightarrow \mathbb{R}$ the indicator
function on $S$:
\[ \mathds{1}_S(x) = \left\{
    \begin{array}{lcl} 1 & \mbox{if} & x \in S \\ 0 & \mbox{if} & x
      \in \R \backslash S
    \end{array}  \right..
\]

Denote by $\mathcal{I}$ the collection of intervals
$I \subset \mathbb{R}$ that are left-closed and right-open, i.e.~of
the type $I = [a,b)$ with $a \in \mathbb{R}$, $a < b \leq \infty$.
Consider the abelian group $LC_B$ generated by the indicator functions
$\sigma_I$ on intervals $I \in \mathcal{I}$, namely:
\begin{equation} \label{LC-1} LC_{B} = \Bigl\{\sigma: \mathbb{R}
  \longrightarrow \mathbb{R} \, \Big| \, \sigma = \sum_{\text{finite}}
  n_I \mathds{1}_I, \; I \in \mathcal{I}, \, n_I \in \Z \Bigr\}.
\end{equation}
Note that the functions $\sigma_I$, $I \in \mathcal{I}$, are not
linear independent. For instance we have
$\mathds{1}_{[a, \infty)} - \mathds{1}_{[b, \infty)} =
\mathds{1}_{[a,b)}$. If we endow
$\mathbb{R}$ with the lower limit topology (i.e.~the topology
generated by $\mathcal{I}$), then the elements of $LC_B$ are precisely
the locally constant integer valued functions (with respect to this
topology) that are bounded and additionally vanish at $-\infty$
(i.e.~at small enough values in $\mathbb{R}$), hence the notation
$LC_B$.

In a similar way, we also have:
\zjnote{\begin{equation} \label{LC-2} LC_B' = \left\{\sum_{\text{finite}} n_I
  \mathds{1}_I \, \Big| \, I \in \mathcal{I} \, \text{is an interval
    of finite length}, \, n_I \in \Z \right\} \subset LC_B.
\end{equation}}
Note that in contrast to $LC_B$, all the functions in $LC'_B$ must
have bounded support, in particular they vanish also at
$+\infty$. The relation
$\mathds{1}_{[a, \infty)} - \mathds{1}_{[b, \infty)} =
\mathds{1}_{[a,b)}$ from $LC_B$ does not hold anymore in $LC'_{B}$ since
$\mathds{1}_{[a, \infty)}$ and $\mathds{1}_{[b, \infty)}$ do not exist
in $LC_B'$. Still, we have the relation
$\mathds{1}_{[a,b)} + \mathds{1}_{[b, c)} = \mathds{1}_{[a,c)}$ for
any $a < b < c< +\infty$. There is yet another variant of the above which we denote by $LC$. The
functions $\sigma \in LC$ are those functions
$\sigma: \mathbb{R} \longrightarrow \mathbb{R}$ that can be written as
\begin{equation} \label{eq:LC-infinite}
  \sigma = \sum_{I \in \mathcal{I}_0} n_I \mathds{1}_{I},
\end{equation}
where the sum is taken only over collections of intervals
$\mathcal{I}_0 \subset \mathcal{I}$ that have the following property:
either $\mathcal{I}_0$ is finite, or
$\mathcal{I}_0 = \{[a_j, b_j)\}_{j \in \mathbb{N}}$ with
\zjnote{$a_j \rightarrow \infty$ as $j \rightarrow \infty$}. Clearly we
have $LC'_B \subset LC_B \subset LC$.
In~\S\ref{subsubsec:Novikov-etc} below we will endow $LC$, $LC_{B}$
with ring structures such that $LC'_{B}$ will become an ideal of
$LC_{B}$. We remark already here that these ring structures are not
the ``standard'' ones in the sense that the multiplication of two
functions from $LC$ (or $LC_{B}$) is not the function obtained by
pointwise multiplication of the two functions, but quite a different
operation.

\

We use here the barcode conventions from~\S\ref{subsec:barfuc}.  Let
$\mathcal{B} = \{(I_j, m_j)\}_{j \in \mathcal{J}}$ be an ungraded
barcode.  There is a natural function $\sigma_{\mathcal{B}}\in LC$
associated to $\mathcal{B}$ by the formula:
\begin{equation} \label{eq:sig-B} \sigma_{\mathcal{B}} = \sum_{j \in
    \mathcal{J}} n_j \mathds{1}_{I_j}.
\end{equation}
One reason this definition is natural is that if $M$ is an (ungraded)
persistence module with barcode $\mathcal{B}_{M}$,
then $$\sigma_{\mathcal B_{M}}(\alpha)=\dim_{\k}(M^{\alpha})~.~$$

The following are straightforward to verify for every ungraded barcode
  $\mathcal{B}$:
  \begin{itemize}
  \item[(i)] $\sigma_{\mathcal{B}} \geq 0$. \label{i:sig-pos}
  \item[(ii)] $\sigma_{\mathcal{B}} \in LC_B$ iff $\mathcal{B}$ has finitely
    many bars.
  \item[(iii)] Assume $\sigma_{\mathcal{B}} \in LC_B$ (i.e.~$\mathcal{B}$ has
    finitely many bars). Then $\sigma_{\mathcal{B}} \in LC'_B$ iff all
    the bars in $\mathcal{B}$ have finite length.
  \item[(iv)] If $\mathcal{B}$ contains only bars of infinite length then
    $\sigma_{\mathcal{B}}:\mathbb{R} \longrightarrow \mathbb{R}$ is a
    non-decreasing function. 
    \item[(v)]  For
  every non-negative and non-decreasing function $\sigma \in LC$ there
  exists a barcode $\mathcal{B}$ all of whose bars are of infinite
  length, such that $\sigma_{\mathcal{B}} = \sigma$.
  \end{itemize}

We now extend the definition of $\sigma_{\mathcal{B}}$ to graded
barcodes. For this we extend the rings $LC$ and $LC_{B}$ by adding formal polynomial variables $x, x^{-1}$
 (that take into account the degree) thus obtaining rings
$LC[x,x^{-1}]$ and $LC_{B}[x,x^{-1}]$.

Let $\mathcal{B} = \{\mathcal{B}_i\}_{i \in \mathbb{Z}}$ be
a graded barcode, with
$\mathcal{B}_i = \{(I_j^{(i)}, m^{(i)}_j)\}_{j \in
  \mathcal{J}^{(i)}}$.  We define
\begin{equation} \label{eq:sig-B-gr} \bar{\sigma}_{\mathcal{B}} := \sum_{i
    \in \mathbb{Z}}  \sigma_{\mathcal{B}_i}x^i = \sum_{i \in
    \mathbb{Z}} \sum_{j \in \mathcal{J}^{(i)}}  m^{(i)}_j
  \mathds{1}_{I^{(i)}_j}x^i .
\end{equation}

We will also need the following expression:

\begin{equation}\label{eq:sigma2}
  \bar{\chi}_{\mathcal{B}}=\bar{\sigma}_{\mathcal{B}}|_{x=-1}=\sum_{i
    \in \mathbb{Z}} (-1)^i \sigma_{\mathcal{B}_i}
\end{equation}

\begin{rem}\label{rem:Euler-ch1}
  Let $M_{\bullet}=\{M_{i}\}_{i\in\Z}$ be a graded persistence module
  whose barcode satisfies the conditions at the beginning of
  \S\ref{sb:diff-Nov}.  Recall that we denote by
  $\mathcal{B}_{M_{\bullet}}$ the graded barcode of $M_{\bullet}$.
  Thus we have:
\zjnote{\begin{equation}\label{eq:sigam-mod}
    \bar{\sigma}_{\mathcal{B}_{M_{\bullet}}}(\alpha)=
    \sum_{i} \dim_{\k} (M^{\alpha}_{i}) x^{i} \left( = \sum_i \sigma_{\mathcal B_{M_i}}(\alpha) x^i\right)
  \end{equation}}
which shows that $\bar{\sigma}_{\mathcal{B}_{M_{\bullet}}}$ simply
tracks
\pbred{the sum of the dimension of the modules $M^{\alpha}_{i}$ over all $i$}
when $\alpha$ varies.  Similarly,
$\bar{\chi}_{\mathcal{B}_{M_{\bullet}}}(\alpha)=\chi_{\alpha}(M_{\bullet})$
where $\chi_{\alpha}(M_{\bullet})$ is the Euler characteristic of
the graded vector space $M^{\alpha}_{\bullet}$ that we have already
seen in~\eqref{eq:Euler-ch0}.
\end{rem}

\subsection{Ring structures on $LC$, $LC_{B}$ and $K$-groups}
\label{subsubsec:Novikov-etc}

We will now endow $LC$ and $LC_B$ with ring structures. Obviously $LC$
and $LC_B$ are closed under standard (point-wise) multiplication of
functions (because
$\mathds{1}_{I} \cdot \mathds{1}_{I'} = \mathds{1}_{I \cap I'}$;
though note that $1_{\mathbb{R}} \notin LC$). However the ring
structure relevant for our purposes will be completely different. \zjnote{To this end we need to use the notions of distributions and
convolutions from classical analysis which we briefly recall
next. }

\medskip

Given a locally integrable function $f$ on $\R$, we identify it
with the regular distribution $T_f$, a linear functional on
$C^{\infty}(\R)$ defined by $T_f(\phi) = \int_{\R} f(x)
\phi(x)dx$. There are two important operations on distributions.
\begin{itemize}
\item[(i)] The derivative $D$ of a distribution, defined by
  $D(T_f) (\phi) : = - T_f(\phi')$.
\item[(ii)] For any two locally integrable functions $f$ and $g$ on
  $\R$, define the convolution $T_f \ast T_g = T_{f \ast g}$, where
  $(f \ast g)(x) := \int_{\R} f(x-y) g(y) dy$.
\end{itemize}
Recall also that we have:
\begin{equation} \label{eq:der-conv} D(T_f \ast T_g) = D(T_f) \ast T_g =
  T_f \ast D(T_g).
\end{equation}
We define a product on $LC$, which we denote by $\circ$, as follows:
\begin{equation} \label{eq:product-LC} \sigma_1 \circ \sigma_2 :=
  D(\sigma_1 \ast \sigma_2).
\end{equation}
Note that by~\eqref{eq:der-conv} we have
$\sigma_1 \circ \sigma_2 = D(\sigma_1) * \sigma_2 = \sigma_1 *
D(\sigma_2)$. The right-hand side of~\eqref{eq:product-LC} is a priori
only a distribution, but as we will see below it is a regular
distribution that moreover corresponds to a function in $LC$. We view
this function as the result of the operation
$\sigma_1 \circ \sigma_2$.  Furthermore, we will also see that if
$\sigma_1, \sigma_2 \in LC_B$ then $\sigma_1 \circ \sigma_2 \in LC_B$.

\medskip

The rings $LC_{B}$, $LC$ are in fact isomorphic to the Novikov rings
that we have seen earlier in the paper. Before we state this result,
let us recall some relevant notation from~\S\ref{subsubsec:nov}
and~\S\ref{sec-fil-cc}: $\La$ is the universal Novikov ring and
$\La_{P}$ and $\La_{P}^{(0)}$ are respectively the Novikov polynomials
and the ideal corresponding to those polynomials that vanish for
$t=1$.

\begin{prop} \label{p:LC-isos} The operation $\circ$ defined
  in~\eqref{eq:product-LC} turns $LC$ and $LC_B$ into commutative
  unital rings where the unity is given by $\mathds{1}_{[0,\infty)}$.
  Moreover, the subgroup $LC'_B \subset LC_B$ is an ideal with respect
  to this ring structure. Furthermore, there is a ring isomorphism
  $\theta: \Lambda \xrightarrow{\; \cong \;} LC$, uniquely defined by
  the property that $\theta(t^a) = \mathds{1}_{[a,\infty)}$. The
  restriction of $\theta$ to $\Lambda_P \subset \Lambda$ gives a ring
  isomorphism $\Lambda_P \xrightarrow{\; \cong \;} LC_B$, and $\theta$
  sends the ideal $\Lambda_{P}^{(0)}$ to $LC'_B$.
\end{prop}
The proof of this proposition is quite elementary, based on direct
calculations and identities involving convolutions. A detailed proof
can be found in~\cite{BCZ-pkt-arxiv}.

\begin{rem} \label{r:conv-sheaves} \zjred{In terms of the
    ``convolution product'' introduced in
    \cite[Section~4]{Scha:operations}, the product $\circ$ on $LC$
    defined in (\ref{eq:product-LC}) is the convolution product inside
    the space of constructible functions. As a concrete example, for
    finite intervals $[a,b)$ and $[c,d)$ in $\R$, one computes that
    $\mathds{1}_{[a,b)} \circ \mathds{1}_{[c,d)} = \mathds{1}_{[a+c,
      \min\{b+c, a+d\})} - \mathds{1}_{[\max\{b+c, a+d\}, b+d)}.$
    Interestingly, this shares certain common features with
    Example~\ref{ex-tensor-elementary} on the tensor product of
    elementary filtered chain complexes.  From a different
    perspective, if we identify $\mathds{1}_{[a,b)}$ with the locally
    constant sheaf $\k_{[a,b)}$ over $\R$, then the computational
    result on the product
    $\mathds{1}_{[a,b)} \circ \mathds{1}_{[c,d)}$ also coincides with
    the ``sheaf convolution'', denoted by $\ast$, of $\k_{[a,b)}$ and
    $\k_{[c,d)}$ (see~\cite[Example~3.6]{Zha20}). The study of the
    sheaf convolution $\ast$ was initiated by Tamarkin in
    \cite[Section~3.1.2]{Tam18} with further applications, for
    instance, to persistence module theory (see
    \cite[Section~2.1]{KS18}).}  \pbnote{Our
    Proposition~\ref{p:LC-isos} above and \cite[Theorem
    3.4]{Scha:operations} identify (via the Euler characteristic) the
    Novikov ring $\Lambda$, the function space $LC$, and the $K$-group
    of the derived category of constructible sheaves on $\R$.}
\end{rem}

In view of Proposition~\ref{p:LC-isos} we have the following relation
for all graded barcodes $\mathcal{B}$:
\begin{equation} \label{eq:theta-lambda-1} \bar{\chi}_{\mathcal{B}} =
  \theta(\lambda_{\mathcal{B}})
\end{equation}
where $\theta$ is the ring isomorphism in Proposition~\ref{p:LC-isos},
$\bar{\chi}_{\mathcal{B}}$ appears in~\eqref{eq:sigma2} and
$\la_{(-)}$ is the map from barcodes to $\La$ as defined in
\S\ref{subsec:barfuc}, see also Remark \ref{rem:various-la}.  This
identity has a $K$-theoretic interpretation that we outline before
proceeding to the proof of Proposition \ref{p:LC-isos}. It is best
understood through the following diagram.

\begin{equation} \label{K-pers-mod-diag} \xymatrix{ \mathcal{PM}\ar[r]
    & Tw(\mathcal{PM}) \ar[r] & H_{0}(Tw(\mathcal{PM}))
    \ar@{-->}[rr]^{\bar{\sigma}}
    & & LC_{B}[x,x^{-1}]\ar[d]_{ (-)|_{x=-1}} \\
    \fchfg \ar[u]_{H} \ar[rr] & & H_{0}\fchfg
    \ar[u]_{E}\ar@{-->}[r]^{\lambda}& \La_{P} \ar[r]^{\theta} & LC_{B}
  }
\end{equation}
Here $\mathcal{PM}$ is the category of persistence modules that have
barcodes with only finitely many bars. This category can be viewed as
a filtered dg-category (with trivial differential) and thus can be
completed, by using twisted complexes, to a pre-triangulated filtered
dg-category $Tw(\mathcal{PM})$ as in~\cite[\S6.1]{BCZ23}. Therefore
the associated homotopy category $H_{0}Tw(\mathcal{PM})$ is a TPC. The
two left top horizontal arrows are the obvious functors. The category
of filtered, finitely generated chain complexes maps to $\mathcal{PM}$
through the persistence homology functor, denoted $H$.  The functor
$E: H_{0}\fchfg \to H_{0}(Tw(\mathcal{PM}))$ is a TPC equivalence
defined as a composition
$$H_{0}\fchfg \stackrel{j}{\longrightarrow}
H_{0}(Tw (\fchfg))\stackrel{H'}{\longrightarrow}
H_{0}(Tw(\mathcal{PM}))$$ where $H'$ is induced by $H$ and $j$ is the
obvious inclusion. It is easy to see that they are both equivalences,
in the case of $j$ because $H_{0}\fchfg $ is already pre-triangulated.
The dashed arrows $\bar{\sigma}$, $\lambda$ are defined only on the
objects of the respective categories ($\bar{\sigma}$ from
(\ref{eq:sig-B-gr}) admits an obvious extension to twisted modules).
In view of Proposition \ref{prop-kc0}, $\lambda$ is the map taking a
chain complex $C$ to its class $[C]\in K(H_{0}\fchfg)\cong \La_{P}$.
The two squares in the diagram commute and $\theta$ is a ring
isomorphism.  In summary, we have
$$\La_{P}\cong K(H_{0} \fchfg )\cong K(H_{0}Tw(\mathcal{PM})) \cong LC_{B}$$
and $\bar{\chi}$ can be viewed as associating to a persistence module
its $K$-class in $LC_{B}\cong K(Tw(\mathcal{PM}))$.

\subsection{Persistence Morse inequalities}\label{subsubsec:morse}
The representation of persistence modules in the ring $LC[x,x^{-1}]$
leads to a natural version of the Morse inequalities in this context.

Let $C = (C_{\bullet}, \partial, \ell)$ be a tame filtered chain
complex - see \S\ref{subsec:pres-mod-pre}. Forgetting the
differential, we can \zjnote{regard} $C$ as a graded persistence module
$\bar{C}_{\bullet}$ with $\bar{C}_k^{\alpha}= C_k^{\leq \alpha}$ and
persistence structural maps given by the inclusions
$C_k^{\leq \alpha} \subset C_k^{\leq \beta}$,
$\forall \alpha \leq \beta$. Its barcode has obviously only infinite
bars.  We also consider the graded persistence module
$H_{\bullet}(\bar{C})$ obtained by taking the persistence homology of
$C$:
$$H_k(\bar{C}^{\alpha})= H_k(C^{\leq \alpha})~.~$$

We put
\zjnote{$$\mathbb P_C := \bar{\sigma}_{\mathcal{B}_{\bar{C}}} \ \mathrm{and} \
\ \mathbb H_{C}:= \bar{\sigma}_{\mathcal{B}_{H(\bar{C})}}$$} where
$\bar{\sigma}$ appears in (\ref{eq:sigam-mod}).  These two functions
belong to the ring $LC[x,x^{-1}]$. They track the dimensions of the
respective persistence modules, in each degree, when $\alpha$ varies.

\begin{prop} \label{perp-Morse-inequality} Assume that the filtered
  chain complex $C$ is bounded. There exists a polynomial
  $\mathbb Q_C(x) \in LC_B[x,x^{-1}]$ such that
  \[ \mathbb P_C(x) - \mathbb H_C(x) = (1+x) \mathbb Q_C(x),\] and
  such that each coefficient of $\mathbb Q_C(x)$ is a non-negative and
  non-decreasing function in $LC_B$.
\end{prop}

Note that when evaluating the coefficients of $\mathbb P_C(x)$,
$\mathbb H_C(x)$ and $\mathbb Q_C(x)$ at $\infty$ (recall that every
function in $LC_B$ is constant near $\infty$), the conclusion in
Proposition~\ref{perp-Morse-inequality} recovers the classical Morse
inequalities expressed in terms of the so-called Poincar\'e
polynomials.

\begin{ex} Here is a simple example reflecting the
  conclusion in Proposition \ref{perp-Morse-inequality}. Assume
  $C = E_2(a,b)$ with the generator in filtration $a$ of degree $0$
  and  the generator in filtration $b$ of degree $1$. Then
  \[ \mathbb P_C(x) = \mathds{1}_{[a, \infty)} + \mathds{1}_{[b,
      \infty)}x \,\,\,\,\mbox{and}\,\,\,\, \mathbb H_C(x) =
    \mathds{1}_{[a,b)}. \] Then we have
  \begin{align*}
    \mathbb P_C(x) - \mathbb H_C(x)
    & =\left( \mathds{1}_{[a, \infty)} +
      \mathds{1}_{[b, \infty)}x\right) - \mathds{1}_{[a,b)}\\
    & = \mathds{1}_{[b, \infty)} +
      \mathds{1}_{[b, \infty)}x = (1+x)\cdot \mathds{1}_{[b, \infty)}.
  \end{align*}
  The function $\mathds{1}_{[b, \infty)} \in LC_B$ is the function
  $\mathbb Q_C(x)$ claimed by Proposition~\ref{perp-Morse-inequality}.
\end{ex}

\begin{proof} [Proof of Proposition~\ref{perp-Morse-inequality}]

We will use the following lemma.

\begin{lemma} \label{sigma-ses-sum} Let $U, V, W$ be ungraded
  persistence $\k$-modules, and $0 \to U \to V \to W \to 0$ a short
  exact sequence of persistence $\k$-modules (in particular, the
  arrows are all persistence morphisms), then we have
  $\sigma_{\mathcal B_V} = \sigma_{\mathcal B_U} + \sigma_{\mathcal
    B_W}$.
\end{lemma}
The proof of the lemma is straightforward since we work over a field
and $\sigma_{\mathcal B_V}(\alpha)=\dim_{\k}(V^{\alpha})$ and
similarly for $U$ and $W$.

\medskip

We now get back to the proof of Proposition~\ref{perp-Morse-inequality}. For any degree $k \in \N$, there are two persistence $\k$-modules
associated to the filtered vector space $C_k$ defined as follows
\[ Z_k = \left\{\left(Z_k^{s} : = \ker(\partial_k: C^{\leq s}_k \to
      C^{\leq s}_{k-1}) \right)_{s \in \R}, \,\, \{i^{Z_k}_{s,t} =
    \mbox{inclusion}\}_{s\leq t \in \R} \right\}, \]
\[ B_k = \left\{\left(B_k^{s} : = {\rm Im}(\partial_k: C^{\leq
        s}_{k+1} \to C^{\leq s}_{k}) \right)_{s \in \R}, \,\,
    \{i^{B_k}_{s,t} = \mbox{inclusion}\}_{s\leq t \in \R} \right\}.\]
We have the following two exact sequences of persistence modules,
 \zjnote{ \[ 0 \to Z_k \xrightarrow{i} \bar{C}_k \xrightarrow{\partial_k}
    B_{k-1} \to 0, \,\,\,\,\mbox{and}\,\,\,\, 0 \to B_k
    \xrightarrow{i} Z_k \xrightarrow{\pi} H_k(\bar{C}) \to 0 \] }where
  \zjnote{$H_k(\bar{C})$} is the degree $k$ part of \zjnote{$H_*(\bar{C})$}. 
  By Lemma
  \ref{sigma-ses-sum}, we have
  \[ \sigma_{\mathcal B_{\bar{C}_k}} = \sigma_{\mathcal B_{Z_k}} +
    \sigma_{\mathcal B_{B_{k-1}}} = \sigma_{\mathcal B_{B_{k}}} +
    \sigma_{\mathcal B_{B_{k-1}}} + \sigma_{\mathcal B_{H_k(\bar{C})}}. \]
  Summing over all the degrees $k$, \zjnote{by definition in (\ref{eq:sigam-mod})}, we get
  \begin{align*}
    \mathbb P_C(x)
    & = \sum_{k} \sigma_{\mathcal B_{\bar{C}_k}} x^k \\
    & = \sum_{k} \left(\sigma_{\mathcal B_{B_{k}}} +
      \sigma_{\mathcal B_{B_{k-1}}} \right) x^k +
      \sum_{k} \sigma_{\mathcal B_{H_{k}\bar{C}}}x^k \\
    & = \left( \sum_{k} \sigma_{\mathcal B_{B_{k}}}x^k\right) +
      \left(\sum_k \sigma_{\mathcal B_{B_{k-1}}} x^{k-1}\right) x  +
      \mathbb H_C(x) \\
    & = (1+x) \left( \sum_{k} \sigma_{\mathcal B_{B_{k}}}x^k\right) +
      \mathbb H_C(x).
  \end{align*}
  The desired polynomial $\mathbb Q_C(x)$ is equal to
  $\sum_{k} \sigma_{\mathcal B_{B_{k}}}x^k$. Note that only
  infinite-length bars appear in each of the persistence $\k$-modules
  $B_k$, therefore each coefficient
  $\sigma_{\mathcal{B}_{B_{k}}}\in LC_B$ of the polynomial
  $\mathbb{Q}_C(x)$ is a non-negative and non-decreasing function.
\end{proof}

A useful application of Proposition \ref{perp-Morse-inequality} is \zjnote{given} next.

\begin{cor}\label{cor:count-bars} Let $\phi : G\to G'$ be a filtered 
  chain \zjnote{map} of two finitely generated chain complexes and let
  $G''= {\rm Cone\/}(\phi)$ be the filtered cone of $\phi$.  We have
  the inequality:
  \zjnote{$$\# (\mathcal{B}_{H(G'')})\leq \# (\mathcal{B}_{H(G)})+\#
    (\mathcal{B}_{H(G')})$$} \zjnote{where $\#(\mathcal B_{\bullet})$
    denotes the cardinality of the intervals} \pbred{in the
    corresponding barcode.}
\end{cor} 

\begin{proof}
  The proof is based on a simple interpretation of the polynomial
  $\mathbb Q_C$ from Proposition
  \ref{perp-Morse-inequality}. \zjnote{Let $C$ be a filtered chain
    complex.} By rewriting the complex $C$ in normal form as in
  (\ref{fcc-decomp})
  $$ C\cong \bigoplus E_{1}(a_{i})[-k_{i}]\oplus
  \bigoplus E_{2}(b_{j},c_{j})[-k_{j}]$$ we see that the value of
  $\mathbb Q_C(1)\in LC_{B}$ for $\alpha=\infty$, \zjnote{denoted by}
  $\mathbb Q_{C}(1)|_{\infty}$, gives the number of the terms
  \zjnote{$E_2(b_{j},c_{j})$} in this decomposition.  This number
  satisfies:
  \zjnote{$$\mathbb Q_{C}(1)|_{\infty}=\# (\mathcal{B}_{H(C)}) +
    N^{0}_{C} - \mathbb H_{C}(1)|_{\infty}$$} where $N^{0}_{C}$ is the
  number of terms \zjnote{$E_2(b_{j},c_{j})$} with
  \zjnote{$b_{j}=c_{j}$} (these terms correspond to so-called ghost
  bars) and clearly $\mathbb H_{C}(1)|_{\infty}$ is the number of
  infinite bars in \zjnote{$\mathcal{B}_{H(C)}$}. From Proposition
  \ref{perp-Morse-inequality} we have that
  $\mathbb Q_{C}(1)= \frac{\mathbb P_{C}(1) - \mathbb H_{C}(1)}{2}$
  and thus
  $$\#(\mathcal{B}_{H(C)})+N^{0}_{C}=\frac{[\mathbb P_{C}(1) - \mathbb
    H_{C}(1)]|_{\infty}}{2} + \mathbb H_{C}(1)|_{\infty}=
  \frac{[\mathbb P_{C}(1) + \mathbb H_{C}(1)]|_{\infty}}{2}$$

  Returning to the statement, assume for a moment, to simplify the
  argument, that $N^{0}_{G}=0=N^{0}_{G'}$. We then have  
 \zjred{ \zjnote{\begin{align*}
 2 \# (\mathcal{B}_{H(G'')})& \leq  \mathbb  P_{G''}(1)|_{\infty} +
\mathbb H_{G''}(1)|_{\infty} \\
& =  \mathbb P_{G}(1)|_{\infty}+\mathbb P_{G'}(1)|_{\infty} + \mathbb H_{G''}(1)|_{\infty} \\
& =  2 \# (\mathcal{B}_{H(G)}) + 2\#(\mathcal{B}_{H(G')}) - \mathbb H_{G}(1)|_{\infty} - \mathbb H_{G'}(1)|_{\infty}+\mathbb H_{G''}(1)|_{\infty}\\
& \leq 2 \# (\mathcal{B}_{H(G)}) + 2\#(\mathcal{B}_{H(G')}).
\end{align*}}}
\zjnote{We have used here the obvious relations
  $\mathbb P_{G''}(1)=\mathbb P_{G}(1)+\mathbb P_{G'}(1)$.}

\medskip

We are left to showing the general case, when $N_{G}^{0}$ and
$N_{G'}^{0}$ do not necessarily vanish. \zjnote{By \cite[Theorem
  B]{UZ16}, for these $G$ and $G'$, there exist filtered chain
  complexes $G_1, G'_1$ and filtered chain homotopies $\psi, \psi'$ as
  follows,
  \[ \psi: G \simeq G_1 \,\,\,\,\mbox{and}\,\,\,\, \psi': G_1 \simeq
    G'_1 \] where $G_1$ and $G'_1$ do not have any ghost bars. In
  particular, $\mathcal{B}_{H(G)}=\mathcal{B}_{H(G_1)}$ and
  $\mathcal{B}_{H(G')}=\mathcal{B}_{H(G'_1)}$. Replace
  $\phi: G \to G'$ in our assumption with a filtered chain map
  $\phi_1: G_1 \to G'_1$, constructed from $\phi$ by (pre)-compositing
  the homotopy inverse of $\psi$ and map $\psi'$, then the resulting
  filtered cone $G''_{1}={\rm Cone}(\phi')$ is filtered
  quasi-isomorphic to $G''$. Thus
  $\mathcal{B}_{H(G'')} = \mathcal{B}_{H(G''_{1})}$. Therefore, we can
  conclude the proof from the previous discussion where no ghost bars
  exist.}
\end{proof}

\begin{rem} (a) Related, and somewhat more refined, bar \zjnote{counting} inequalities 
appear in \cite{BPPPSS22} where they are shown by different methods.

(b) It is easy to see that the Corollary above implies that if 
$$0\to M\to N\to Z\to 0$$
is a short exact sequence of persistence modules with $M$ and $Z$
having \zjnote{barcodes} with finite numbers of bars, then the number
of bars in the \zjnote{barcode} of $N$ is at most the sum of the
numbers of bars of $\mathcal{B}_{M}$ and $\mathcal{B}_{Z}$, a result
that again appears in \cite{BPPPSS22}.

\zjnote{(c) By rotating an exact triangle $G \to G' \to G'' \to TG$ of
  filtered chain complexes in the category $\mathcal C^{\rm fg}$ (see
  \cite[Remark 2.46]{BCZ23}), the Corollary above also implies that
  $\# (\mathcal{B}_{H(G)})\leq \# (\mathcal{B}_{H(G')})+\#
  (\mathcal{B}_{H(G'')})$.}
\end{rem}


\section{A $K$-theoretic pairing} \label{sec-pairing}

In this section, we will prove one of the main results in this paper,
that is, there exists a bilinear pairing on the $K$-group of a TPC
$\C$ under the assumption that $\C$ is tame and bounded.

\subsection{The pairing}\label{subsec:pairing}
We will work here with a tame and bounded triangulated persistence
category $\mathcal{C}$ as in Definition~\ref{def:tpc-tame}, with shift
functor $\Sigma$ and translation functor $T$.

We will construct a group homomorphism
$q_X: K(\C) \longrightarrow \Lambda$ for every $X \in {\rm Obj}(\C)$
as follows.  Let $X \in {\rm Obj}(\C)$. Given $N \in {\rm Obj}(\C)$,
recall the graded persistence module
$$\ho_{\C}(X,T^{\bullet}N) := \{\ho_{\C}(X,T^iN)\}_{i \in
  \mathbb{Z}}$$ from (\ref{eq:pers-grad-mor}).  Define
$\tilde{q}_X: {\rm Obj}(\C) \longrightarrow \Lambda$ by
\begin{equation} \label{tq-X} \tilde{q}_X(N) :=
  \lambda_{\ho_{\C}(X,T^{\bullet}N)},
\end{equation}
where $\lambda_{(-)}$ is as in~\eqref{eq:mor-pers-la} (see
also~\eqref{eq:lambda-M-blt}). This is well-defined because $\C$ is
tame and bounded.
We can extend $\tilde{q}_X$ linearly on the free abelian group
generated by ${\rm Obj}(\C)$.
  
\begin{thm} \label{thm-pairing} The linear map $\tilde{q}_X$ defined
  above descends to a homomorphism
  $$q_X: K(\C) \longrightarrow \Lambda~.~$$ 
  
  Moreover, the map
  $\kappa : {\rm Obj} (\C)\times {\rm Obj} (\C)\to \La_P$ defined by
  $\kappa (X,Y)= q_{X}(Y)$ induces a bi-linear map of
  $\Lambda_{P}$-modules
  $$\kappa: \overline{K(\C)} \otimes_{\Lambda_{P}} K(\C)
  \longrightarrow \Lambda ~.~$$
\end{thm}
Here $\overline{K(\C)}$ is the same as $K(\C)$ as a group but the
$\La_{P}$-structure is reversed in the sense that $t^{a}$ acts as
$\Sigma^{-a}$, as in \S \ref{subsec:conclusion-alg}, where
$\Sigma^{s}$ is the shift functor in $\C$.  We denote the action of
$\La_{P}$ on $K(\C)$ by $\cdot$ and the action on $\overline{K(\C)}$
by $\overline{\cdot}$.
 
\begin{rem}\label{rem:skew-symm}
  \begin{enumerate}
  \item Given that $\kappa$ is bilinear with the $\La_{P}$-module
    $\overline{K(\C)}$ on the first term, it follows that $\kappa$ is
    not symmetric. Indeed, assuming $\kappa(X,Y)=\kappa(Y,X)\not=0$
    for two objects $X$, $Y$, then, for $r\not=0$, we have
    $\kappa (\Sigma^{r}X,Y) =t^{-r}\kappa (X,Y)=
    t^{-r}\kappa(Y,X)=t^{-2r}\kappa(Y, \Sigma^{r}X)$ and thus
    $\kappa(\Sigma^{r}X,Y)\not=\kappa (Y,\Sigma^{r}X)$. We will see
    that in many cases of interest the form $\kappa$ is {\em
      skew-symmetric} in the sense that we have
    \begin{equation}\label{eq:skew-symm}
      \kappa(X,Y)(t)=(-1)^{\epsilon}\kappa(Y,X)(t^{-1})
    \end{equation}
    where $\epsilon=\epsilon_{\C}\in \Z/2$ is a global
    constant associated to the category $\C$.
  \item The self-pairing $\kappa(X,X)$, applied to objects
    $X \in \Ob(\mathcal{C})$, is invariant under shifts and
    translations in the sense that
    $\kappa(\Sigma^r T^k X, \Sigma^r T^kX) = \kappa(X,X)$ for every
    $r \in \mathbb{R}$, $k \in \mathbb{Z}$, where $\Sigma^r$ and $T$
    are the shift and translation functors, respectively. This easily
    follows from the definitions.
  \end{enumerate}
\end{rem}

\begin{ex} \label{ex-pairing-kf} Consider the category
  $\C^{\fg} = H_0(\fchfg)$. We will describe $\kappa$ in this case by
  using the identification $K(\C^{\fg}) \cong \Lambda_{P}$ as in
  Proposition \ref{prop-kc0}.  It is enough to calculate
  $\kappa(E_1(a), E_1(c))$ and extend this by linearity. We claim:
  \begin{equation} \label{pairing-easy-rel} \kappa(E_1(a), E_1(c)) =
    t^{c-a}.
  \end{equation}
  Indeed, we have $\ho_{\C}(E_1(a), T^i E_1(c)) = E_1(c-a)$ when
  $i = 0$ and equal to $0$ otherwise (here $T^i$ is the shift of the
  degree by $i$), so
  $$\kappa(E_1(a), E_1(c)) = q_{E_1(a)}(E_1(c))
  = \lambda_{\ho_{\C}(E_1(a), E_1(c))} = t^{c-a}.$$
  Therefore, using the identification of $K(\C^{\fg})$ with $\La_{P}$,
  $\kappa$ can be rewritten on elements of $\La_{P}$ in the form
  $$\kappa (P,Q)(t)=P(t^{-1})Q(t)~.~$$
  Thus, for $\mathcal{C}^{\fg}$ the pairing $\kappa$ is non-degenerate
  and skew-symmetric (with $\epsilon=0$).
\end{ex}

\begin{proof}[Proof of Theorem \ref{thm-pairing}.]
  We start with the statement concerning $q_{X}$.  We need to show
  that for each exact triangle $A \to B \to C \to TA$ in $\C_0$ we
  have
  \begin{equation}\label{eq:exact3}
    \tilde{q}_X(A) - \tilde{q}_X(B) + \tilde{q}_X(C) = 0.
  \end{equation}
  To prove this relation in $\La$ we first remark that a Novikov
  polynomial $P \in \La$ vanishes if and only if we have
  $P_{\alpha}(1)=0$ for all $\alpha\in\R$. Here and below, we denote
  by $P_{\alpha}$ the $\alpha$-truncation of $P$, namely the sum of
  the monomials in $P$ having $t$ raised at exponents $\leq \alpha$.

  Given a graded persistence module $M_{\bullet}$ (of finite total
  rank below each filtration level) we can define its $\alpha$-Euler \zjred{characteristic} as:
  \begin{equation} \label{eq:Euler-pers} \chi_{\alpha}(M_{\bullet}):=
    \chi(M^{\alpha}_{\bullet})~.~
  \end{equation}

  Using the preceding notation we have:
  \begin{equation} [ \lambda_{M_{\bullet}}]_{\alpha}(1)
    =\chi_{\alpha}(M_{\bullet})
  \end{equation}
  This follows directly from formula~\eqref{eq:lambda-M-blt}. See also
  Remark~\ref{rem:Euler-ch1}.
  
  In summary, to prove~\eqref{eq:exact3} we need to show that for each
  $\alpha$ we have:
  \begin{equation}\label{eq:Euler-char-eq2}
    \chi_{\alpha} (\ho_{\C}(X,T^{\bullet}A)) -
    \chi_{\alpha}(\ho_{\C}(X,T^{\bullet}B))+
    \chi_{\alpha}(\ho_{\C}(X,T^{\bullet}C ))=0.
  \end{equation}
  Notice that $\ho^{\alpha}_{\C}(X, \cdot)$ is a cohomological functor
  on $\C_0$. We apply it to the exact triangle $A\to B\to C\to TA$ and
  get a long exact sequence
  \[\cdots \to \ho^{\alpha}_{\C}(X, T^iA) \to
    \ho^{\alpha}_{\C}(X, T^iB) \to \ho^{\alpha}_{\C}(X, T^iC) \to
    \ho^{\alpha}_{\C}(X, T^{i+1}A) \to \cdots.\] By the boundedness
  assumption, this long exact sequence has only finitely many terms.
  We view this long exact sequence as a cochain complex $S^{\bullet}$
  (with trivial homology) and regrade its terms so that the degree of
  $\ho^{\alpha}_{\C}(X, T^iA)$ is $3i$, the degree of
  $\ho^{\alpha}_{\C}(X, T^iB)$ equals $3i+1$ and finally the degree of
  $\ho^{\alpha}_{\C}(X, T^iC)$ is $3i+2$.  The Euler characteristic of
  $S^{\bullet}$ vanishes, hence:
  \begin{align*}
    0 = \chi (S^{\bullet}) & = \sum_{i\in \N} (-1)^{3i} {\rm dim\,}
                             \ho^{\alpha}_{\C}(X, T^iA)  +
                             \sum_{i\in \N} (-1)^{3i+1} {\rm dim\,}
                             \ho^{\alpha}_{\C}(X, T^iB) \\
                           & \,\,\,\,\,\,+ \sum_{i\in \N} (-1)^{3i+2}
                             {\rm dim\,} \ho^{\alpha}_{\C}(X, T^iC).
  \end{align*}
  The parities of $3i$ and $3i+2$ are the same as that of $i$, while
  the parities of $3i+1$ and $i$ are opposite. Therefore,
  \begin{align*} \label{euler-char}
    0 &= \sum_{i\in \N} (-1)^{i} {\rm dim\,} \ho^{\alpha}_{\C}(X, T^iA)
        - \sum_{i\in \N} (-1)^{i} {\rm dim\,} \ho^{\alpha}_{\C}(X, T^iB) \\
      & \,\,\,\,\,\,+ \sum_{i\in \N} (-1)^{i} {\rm dim\,}
        \ho^{\alpha}_{\C}(X, T^iC)
  \end{align*}
  which shows relation (\ref{eq:Euler-char-eq2}) and concludes the
  proof.

  We now turn to the part of the statement concerned with the map
  $\kappa$.  Given an exact triangle $X\to Y\to Z\to TX$ in $\C_{0}$,
  we need to show that for each $N$ we have:
  \begin{equation}\label{eq:exact-3}
    \tilde{q}_X(N) - \tilde{q}_Y(N) + \tilde{q}_Z(N) = 0 ~.~
  \end{equation}
  Given that the functor $\ho^{\alpha}_{\C}(-, N)$ is homological on
  $\C_{0}$ and that
 \begin{equation} \label{equ-hom-duality}
 \ho_{\C}(T^{-i}X, N) \simeq \ho_{\C}(X, T^iN) \; \; \forall i,
 \end{equation}
  the same argument as in the first part shows~\eqref{eq:exact-3}. We
  deduce that $\kappa$ descends to the tensor product of groups
  $K(\C)\otimes K(C)$.  The $\La_{P}$-module part of the statement now
  follows because:
  \[ t^a \cdot \kappa([X], [Y]) = \kappa(t^a \ \overline{\cdot}\ [X],
    [Y]) = \kappa([X], t^a \cdot [Y]) ~.~\] This is due to the fact
  that
  $$\ho^{\alpha}_{\C}(\Sigma^{-r}X,Y)=\ho^{\alpha}_{\C}(X,\Sigma^r(Y))$$
  for all $\alpha, r\in \R$.
\end{proof}

\begin{rem}
  The pairing $\kappa$ above is natural from the point
  of view of a more general formalism of \zjnote{a filtered version of} derived functors.
  \zjnote{More details will appear in \cite{BCZ-derived}}.
\end{rem}

\subsection{The exact sequence~\eqref{eq:exact-seq-k} and the pairing
  $\kappa$} \label{subsec:ex-pair} Assume that $\C$ is a TPC and
moreover tame and bounded (see Definition~\ref{def:tpc-tame}).  Recall
that we have the exact sequence of $\La_{P}$ modules:
\begin{equation} \label{eq:ES-again}
  0\longrightarrow
  \tor K(\C)\longrightarrow K\ac (\C)\stackrel{j}{\longrightarrow}
  K(\C)\longrightarrow K_{\infty}(\C)\longrightarrow 0~.~
\end{equation}
\zjnote{Also r}ecall from Theorem~\ref{thm-pairing} that we have a pairing
$\kappa: \overline{K(\C)} \otimes_{\Lambda_{P}} K(\C)
\longrightarrow \Lambda$. We will discuss here some relations between
this pairing and the exact sequence~\eqref{eq:ES-again}.

\begin{cor}\label{cor:pairing-ex} There exists a pairing:
  $$\kappa' : \overline{K\ac (\C)}\otimes_{\La_{P}} K (\C)
  \longrightarrow \La'$$ that vanishes on
  $\overline{\mathrm{Tor} K(\C)}\otimes_{\La_{P}} K(\C) $ where $\La'$
  is defined in~\eqref{nov-ser-double}. The pairings $\kappa'$ and
  $\kappa$ are related by:
  $\sigma \circ \kappa' = \kappa \circ (\bar{j} \otimes \id)$, where
  $\bar{j}: \overline{K\ac (\C)} \longrightarrow \overline{K(\C)}$ is
  the same as $j$ from~\eqref{eq:ES-again} but viewed here as a map
  between modules with the reversed $\Lambda_P$-structure.  The map
  $\sigma: \Lambda' \longrightarrow \Lambda$ is defined
  in~\eqref{sigma-map} and in~\S\ref{ssec-(T)}.

  In case $\C$ is of finite type, then the pairing $\kappa$ has values
  in $\La_{P}$ and the pairing $\kappa'$ has values in
  $\La'_{P}\simeq\La_{P}^{(0)}$. The pairing $\kappa$ induces the
  usual Euler pairing
  $$\kappa'':K_{\infty}(\C)\otimes K_{\infty}(\C) \longrightarrow \Z$$
  through the formula $\kappa''=\kappa|_{t=1}$.
\end{cor}

\begin{proof}
  The construction of $\kappa'$ is just like that of $\kappa$ but is
  based on a map $q'_{X}$ that replaces $q_{X}$, for
  $X\in \mathrm{Obj}(\ac \C)$. In this case:
  $\tilde{q}'_{X}(N)=\lambda'_{\ho_{\C}(X,T^{\bullet}N)}$ and
  $\lambda'_{(-)}$ is defined for persistence modules with barcodes
  with only finite bars through the formula $\lambda'_{I}= s^{a,b}$
  for $I=[a,b)$, extended by linearity. Given that $X$ is acyclic,
  $\ho_{\C}(X,T^{\bullet}N)$ has a barcode with only finite bars and
  thus this $q'_{X}(-)$ is well defined.  Recall from Proposition
  \ref{prop-kac0-2} that there is a map $\sigma : \La'\to \La$ defined
  by $\sigma (s^{a,b})=t^{a}-t^{b}$ which is an isomorphism. It is
  clear that
  $\sigma \circ \kappa'=\kappa \circ (\bar{j} \otimes \id|_{K(\C)})$
  where $\bar{j}: \bar{K\ac(\C)} \longrightarrow \bar{K(\C)}$ is the
  same map as $j$, but with domain and target being $K\ac(\C)$ and
  $K(\C)$ with the reversed $\Lambda_P$-module structures.

  Out of this we deduce that $\tilde{q}'_{X}$ descends to a map
  $q'_{X}$ defined on the \zjnote{$K$-group} and that $\kappa'$ is well defined
  through $\kappa'(X, Y)=q'_{X}(Y)$. By using the exact sequence
  (\ref{eq:exact-seq-k}) it also follows that $\kappa'$ vanishes on
  $\overline{\mathrm{Tor} K(\C)}\otimes_{\La_{P}} K(\C) $.

  In case $\C$ is of finite type, then, by Definition
  \ref{def:tpc-tame} all the persistence modules
  $\ho (X,T^{\bullet}Y)$ have only finitely many bars. Thus the image
  of $q_{X}$ is in $\La_{P}$ and so is that of $\kappa$. The classical
  Euler characteristic of the graded vector space
  $\ho_{\C_{\infty}} (X,T^{\bullet}Y)$ is well defined as an integer
  and it is equal to $\kappa (X,Y)(1)$ because the finite bars have a
  trivial contribution to this sum. Finally, in this case $\kappa'$
  lands in $\La'_{P}$ which is mapped by $\sigma$ isomorphically to
  $\La_{P}^{(0)}$ \zjnote{(see Lemma \ref{lem-lambda-2})}.
\end{proof}

\subsection{Skew-symmetry of the bilinear form $\kappa$}
\label{sb:skew-sym-kappa}
We discuss here conditions \zjnote{which} ensure that the bilinear form $\kappa$
is skew-symmetric in the sense of~\eqref{eq:skew-symm}.

\begin{ex}\label{ex:duality-chain}
  As a motivating example, we return to the category $\fchfg$. We
  already know from Example \ref{ex-pairing-kf} that $\kappa$ is
  skew-symmetric in this case, but regardless, it is worth noticing
  that this fact is a reflection of the Calabi-Yau property of
  $\fchfg$ (see e.g.~\cite{Co:CY, KoSo:ai, Ga:symp-duality, Camp}).

  Let $C = (C_*, \partial_*, \ell) \in \Ob(\fchfg)$, namely a finitely
  generated filtered chain complex. Following~\cite[Section~2.4]{UZ16}
  one can construct an algebraic dual of $C$ (with degree shifted by a
  fixed number $m_0 \in \Z$): this is a filtered chain complex,
  denoted by
  \begin{equation} \label{dual-ch} C^{\vee} = \left(C_*^{\vee}: = {\rm
        Hom}(C_{m_0-*}, \k), \delta_*, \ell^{\vee} \right)
  \end{equation}
  where $\delta_*: C_*^{\vee} \to C_{*-1}^{\vee}$ is defined by
  $\phi \in {\rm Hom}(C_{m_0-*}, \k) \,\mapsto \, \delta_*(\phi) :=
  (-1)^{*} \cdot \phi \circ \partial_{m_0-*} \in {\rm
    Hom}(C_{m_0-*+1}, \k)$ and the dual filtration function is defined
  by
  \begin{equation} \label{dual-fil}
    \ell^{\vee}(x^{\vee}_i) = - \ell(x_i)  
  \end{equation}
  if $\{x_i\}_{i}$ is a basis for $C_*$ and $\{x^{\vee}_i\}_i$ is the
  corresponding dual basis for $C^{\vee}_*$.

  We have the following relation:
  \begin{equation}\label{eq:symm-ch}
    \lambda(C)(t) = (-1)^{m_0} \lambda(C^{\vee})(t^{-1}) ~.~
  \end{equation}
  Indeed, Proposition 6.7 in \cite{UZ16} implies that the
  decompositions (\ref{fcc-decomp}) for $C$ and $C^{\vee}$ are related
  via the following correspondence:
  \begin{itemize}
  \item[(i)] $E_1(a)$ is a summand of $C$ if and only if
    $E_1(-a)[-m_0]$ is a summand of $C^{\vee}$;
  \item[(ii)] $E_2(a,b)$ is a summand of $C$ if and only if
    $E_2(-b,-a)[-m_0+1]$ is a summand of $C^{\vee}$
  \end{itemize}
  and this implies (\ref{eq:symm-ch}). The category $\fchfg$ satisfies
  a form of Calabi-Yau duality in the sense \zjred{that} for each two
  objects $X,Y$ there is a filtered quasi-isomorphism
  \begin{equation}\label{eq:dual-ch2} \varphi_{X,Y}: \bigl(\ho (Y, X)
    \bigr)^{\vee} \longrightarrow \ho (X,Y)
  \end{equation} 
  Using now~\eqref{eq:symm-ch} with $m_{0}=0$ we see that $\kappa$ is
  skew-symmetric on $H_{0}(\fchfg)$.
\end{ex}
 
Example~\ref{ex:duality-chain} admits obvious generalizations. For
instance, suppose $\mathcal{A}$ is a pre-triangulated filtered
dg-category such that each filtered vector space
$\ho_{\mathcal{A}} (X,Y)$ for $X,Y\in \mathrm{Obj}(\mathcal{A})$ is
finite dimensional and assume further that $\mathcal{A}$ satisfies a
filtered form of the Calabi-Yau property such that, in particular,
that for each two objects $X,Y$ there is a filtered quasi-isomorphism
as in~\eqref{eq:dual-ch2}). We then again deduce that $\kappa$ is
skew-symmetric on the \zjnote{$K$-group of $H_{0}(\mathcal{A})$}. Of course, a
similar result remains true for filtered $A_{\infty}$-categories as
long as the filtered quasi-isomorphisms~\eqref{eq:dual-ch2} exist.
  
\

It is natural to wonder what is the analogue of the Calabi-Yau
property for TPCs. One possible definition - that will not be
developed in full here - is roughly as follows. As a preliminary,
consider a persistence module $M$ that is assumed to have only
finitely many bars.  Define a new persistence module $M^{\vee}$ as
follows $$(M^{\vee})^{\alpha}= (M^{-\alpha})^{\ast}, $$ where
$V^* := \hom_{\k}(V,\k)$ is the dual vector space of $V$. The
persistence structural maps on $M^{\vee}$ are induced in the obvious
way and possible degree shifts are neglected.
 
\begin{rem} Notice that if $M$ is lower semi-continuous and lower
  bounded as in the definition of triangulated persistence categories
  used in this paper \zjnote{(see \S\ref{subsec:pres-mod-pre})}, then
  $M^{\vee}$ is upper semi-continuous and upper bounded.  In terms of
  barcodes, the bars of $M$ are of the form $[a,+\infty)$ and $[c,d)$
  and the corresponding ones for $M^{\vee}$ are $(-\infty, -a]$ and
  $(-d,-c]$, respectively.
\end{rem}
 
Assume now that $\C$ is a TPC. We will say that $\C$ is Calabi-Yau if
the triangulated category category $\C_{\infty}$ is (strict)
Calabi-Yau and moreover the following holds: for any two objects $X$,
$Y$ of $\C$ there is an isomorphism
$\Phi= \Phi_{X,Y}: [\ho_{\C_{\infty}}(Y,X)]^{\ast}\to
\ho_{\C_{\infty}}(X,Y)$ and there are persistence morphisms:
$$\varphi=\varphi_{X,Y}: [\ho_{\C}(Y, X) ]^{\vee}\to \ho_{\C}(X,Y)$$
such that for each $\alpha\in \R$ the natural map
$\ker (\varphi ^{-\infty}) \to \ker(\varphi^{\alpha})$ is surjective
and the map $\coker (\varphi ^{\alpha}) \to \coker(\varphi^{\infty})$
is injective. Moreover, the isomorphism $\Phi$ is required to be tied
to $\varphi$ as follows.  First notice that
$\ker (\varphi^{-\infty})=[\ho_{\C_{\infty}}(Y,X)]^{\ast}$ and
\zjnote{$\ho_{\C_{\infty}}(X,Y)=\coker (\phi_{\infty})$} so we can view
\zjnote{$\Phi: \ker (\varphi^{-\infty}) \to \coker (\varphi^ {\infty})$}. We
now require for each $x\in \ker (\varphi^{-\infty})$,
$$\sup \left\{\alpha \in \R \ | \ x\in \ker (\varphi^{\alpha})\right\}=
\inf \left\{\beta\in \R\ | \ \Phi(x)\in\coker (\varphi^{\beta})
\right\}~.~$$ These conditions ensure, in particular, that if the
\zjnote{barcode} in the domain of $\varphi$ contains, respectively,
$(-\infty, a]$ and $(c,d]$ then the image of $\varphi$ contains
$[a,\infty)$ and $[c,d)$. As a result, if $\C$ is Calabi-Yau in this
sense, then the bilinear form $\kappa$ is skew-symmetric on $K(\C)$.
\begin{rem} \label{rmk-usher-new} \pbred{The Poincar\'e duality
    property for persistence modules, or more essentially for filtered
    chain complexes, has been developed in detail in the recent
    work~\cite{Ush23}, where it is called the chain-level
    Poincar\'e-Novikov structure (see~\cite[Section~6]{Ush23}). This
    could provide a rigorous foundation for the Calabi-Yau property of
    our TPCs.}
\end{rem}

\section{Measurements on the $K$-group and on
  barcodes.}\label{S:bar-m-p}
In contrast to usual triangulated categories, triangulated persistence
categories are endowed with exact triangles of arbitrary non-negative
weights. As a result, there is a notion of ``energy'' required to
decompose objects through (iterated) exact triangles and, in
particular, the set of objects of a TPC, $\C$, carries a class of
pseudo-metrics called fragmentation pseudo-metrics as recalled in
\S\ref{subsec:tr-weight}. In this section we first see in
\S\ref{subsubsec:metrics-K-gen} that these metric structures descend
to \zjnote{$K(\C)$} and they sometimes give non-trivial structures
there. On the other hand, in \S\ref{s:num-inv} we discuss numerical
invariants associated to persistence modules and barcodes that are
natural from the point of view of the representation in $LC$ from
\S\ref{subsubsec:iso-rings}.

\subsection{Group semi-norms on $K(\C)$}\label{subsubsec:metrics-K-gen}  
Let $\C$ be a TPC. Let $d$ be a pseudo-metric on $\mathrm{Obj}(\C)$
that is subadditive in the sense that
\begin{equation}\label{eq:add}
  d(X \oplus X', Y \oplus Y') \leq d(X, Y) + d(X', Y') ~.~
\end{equation}
Such pseudo-metrics exist in abundance. Indeed, as shown
in~\cite[\S5]{BCZ21} and recalled in~\S\ref{subsec:tpc}, the category
$\C_{\infty}$ admits a persistence triangular weight $\overline{w}$
that is subadditive.  As a result, ${\rm Obj}(\C)$ carries a class of
pseudo-metrics $d^{\mathcal{F}}$ (where $\mathcal{F}$ is an auxiliary
family of objects of $\C$), called fragmentation pseudo-metrics,
associated to $\overline{w}$. Their construction is recalled in
\S\ref{subsec:tr-weight}.  The subadditivity of $\overline{w}$ implies
(\ref{eq:add}) for all such $d^{\mathcal{F}}$.

The main result in this section deals with ``measurements'' induced by
these pseudo-metrics on $K(\C)$.  There are two natural notions that
we will discuss.

First, define a map $\bar{d}: K(\C) \times K(\C) \to \R_{\geq 0}$ as
follows,
\begin{equation} \label{induce-metric} \bar{d}(x, y) = \inf\{d(X, Y)
  \,| \, X, Y \in \mathrm{Obj}(\C) \,\,\mbox{and}\,\, [X]=x, \, [Y] =
  y \}.
\end{equation}
We will see below that this definition leads to a pseudo-metric on
$K(\C)$ that is invariant by translation and, thus, the map
$x \longmapsto \bar{d}(x,0)$ provides a so-called group semi-norm
(possibly infinite) on $K(\C)$ that determines $\bar{d}$. However, we
will see that this pseudo-metric (and the associated semi-norm), is
trivial in some of our main examples.

Another possibility is to consider a sub-class
\zjnote{$\mathcal{S}\subset \mathrm{Obj}(\C)$} such that $\mathcal{S}$ is
closed with respect to the direct sum and, moreover, each element of
$K(\C)$ has a representative in $\mathcal{S}$ and define
$||-||_{d, \mathcal{S}}: K(\C) \longrightarrow \R_{\geq 0}$ by
\begin{equation}\label{eq:semi}
  ||x||_{d,\mathcal{S}}= \inf\{ d(X,0) \, |
  \, X\in\mathcal{S} \,\,\mbox{and}\,\, [X]=x\}~.~
\end{equation}
This definition also gives a group semi-norm on $K(\C)$, that we will
call the {\em strong semi-norm induced by} $d$ (relative to
$\mathcal{S}$). In contrast to $\bar{d}(-,0)$, there are simple
examples when $d$ is a fragmentation pseudo-metric and
$||-||_{d, \mathcal{S}}$ is finite and not trivial. This will be
shortly shown in Proposition~\ref{thm:ind-metr} below.

\begin{rem} \label{r:seminorms} Before we go on let us explain in more
  detail the essential difference between $\bar{d}$ and $\|-\|_{d,S}$.
  By definition $\bar{d}(x,0)$ is the infimum of the ``distances''
  $d(X,Y)$ where $X$ and $Y$ go over all elements of
  $\Ob(\mathcal{C})$ whose classes in $K(\mathcal{C})$ are $[X]=x$ and
  $[Y]=0$. (Thus the second entry in $\bar{d}(x,0)$ stands for
  $0 \in K(\mathcal{C})$.)  By contrast, in the definition of
  $\|x\|_{d,S}$ we infimize $d(X,0)$ only over $X$, so that the second
  entry in $d(X,0)$ is taken to be the $0$ object rather than
  only a representative of $0 \in K(\mathcal{C})$. In addition, we vary
  $X$ only over the subset $\mathcal{S} \subset \Ob(\mathcal{C})$.
\end{rem}

\begin{prop}\label{thm:ind-metr} Let $\mathcal{C}$ be a triangulated
  persistence category.
  \begin{itemize}
  \item[(i)] Any pseudo-metric $d$ on ${\rm Obj}(\C)$ that is
    \zjnote{subadditive} as in (\ref{eq:add}) induces a pseudo-metric
    $\bar{d}$ given by (\ref{induce-metric}) that is translation
    invariant in the sense that $$\bar{d}(x+z,y+z)=\bar{d}(x,y)$$ for
    all $x,y,z\in K(\C)$.
  \item[(ii)] \zjnote{The bottleneck metric $d_{\rm bot}$ on the
      objects of
      $\C^{\rm fg}= H_{0}(\mathcal F{\rm {\bf Ch}}^{\rm fg})$ descends
      to a pseudo-metric on $K(\C^{\rm fg})\cong \La_{P}$ that only
      takes the values $0$ and $\infty$.} We refer the reader
    to~\cite{PRSZ20} for the definition and more details on the
    bottleneck metric.
  \item[(iii)] Let $d^{\mathcal{F}}$ be the fragmentation
    pseudo-metric on $\C^{\rm fg}$ associated to the persistence
    weight $\bar{w}$ on $\C^{\rm fg}$ with
    $\mathcal{F}=\{E_{1}(0)[s]\}_{s\in\Z}$. (See~\S\ref{subsec:tr-weight}
    for the definitions.) The induced pseudo-metric
    $\bar{d}^{\mathcal{F}}$ on $K(\C^{\rm fg})\cong \La_{P}$ vanishes.
  \item[(iv)] The application $||- ||_{d,\mathcal{S}}$ from
    (\ref{eq:semi}) is a group semi-norm on $K(\C)$ and, for
    $\C=\C^{\rm fg}$, $d=d^{\mathcal{F}}$ as at point (iii), and
    $\mathcal{S}$ the subclass of filtered chain complexes with
    trivial differential, $|| -||_{d,\mathcal{S}}$ is finite and we
    have $$||t^{a}||_{d, \mathcal{S}}=|a|$$ for all $a<0$. Thus
    $||-||_{d,\mathcal{S}}$ is not trivial in this case.
  \end{itemize}
\end{prop}

A detailed proof of the proposition can be found
in~\cite{BCZ-pkt-arxiv}.

\begin{rem}
  \begin{enumerate}
  \item A similar result to point (iii) of
    Proposition~\ref{thm:ind-metr} has been independently established
    by Berkouk~\cite{Ber:non-stab} in the more general context of the
    $K$-group of the category of constructible sheaves over vector
    spaces of arbitrary finite dimension.
  \item If the family $\mathcal{F}$ used before is replaced with
    $\mathcal{F}=\{E_{1}(r)[s]\}_{s\in\Z}$, for some fixed $r > 0$,
    the resulting strong fragmentation semi-norm
    $||-||_{d,\mathcal{S}}$ will not vanish on polynomials with
    $a_{0} <r$ (here $\mathcal{S}$ is as before the class of filtered
    complexes with trivial differential). On the other hand, if the
    family $\mathcal{F}$ is as in the proposition but
    \zjnote{$\mathcal{S}=\mathrm{Obj}(\C^{\fg})$}, then the resulting
    \zjnote{pseudo-norm} $|| - ||_{d,\mathcal{S}}$ is trivial.
  \item The strong fragmentation pseudo-norm associated to the
    bottleneck metric on $\C^{\fg}$ is infinite for all complexes $A$
    such that $H(A)\not=0$.
  \end{enumerate}
\end{rem}


\subsection{$K$-theoretic numerical invariants}
\label{s:num-inv}

In this section we will construct several numerical invariants of
filtered chain complexes that can be defined at the level of the
$K$-theories of $H_0(\fchfg)$, $\mathcal{A} H_0(\fchfg)$ and
$H_0(\fcht)$. This means that the values of our invariants on a given
object $C$ (i.e.~a filtered chain complex) \zjred{depend} only on its
$K$-class.

Since (graded) barcodes are in $1-1$ correspondence with isomorphism
types of filtered chain complexes in the categories $H_0(\fchfg)$ and
$H_0(\fcht)$, we will view the numerical invariants constructed below
as invariants associated to barcodes.

Before we begin recall that we have the following isomorphisms
\begin{equation} \label{eq:K-lambda}
  K \mathcal{A} (H_0(\fchfg)) \xrightarrow{\; \cong \; }
  \Lambda_P^{(0)}, \quad K(H_0(\fchfg)) \xrightarrow{\; \cong \; }
  \Lambda_P, \quad K(H_0( \fcht)) \xrightarrow{\; \cong \; } \Lambda,
\end{equation}
which assign to the $K$-class $[C]$ of a filtered chain complex $C$
the Novikov polynomial (resp. series)
$\lambda([C]) = \lambda_{\mathcal{B}(H_{\bullet}(C))}$, where
$H_{\bullet}(C)$ is the persistence homology of $C$,
$\mathcal{B}(H_{\bullet}(C))$ is its (graded) barcode and
$\lambda_{\mathcal{B}(-)}$ is defined in~\S\ref{subsec:barfuc}. By
abuse of notation we denote all three isomorphisms above by $\lambda$
since they are compatible one with the other with respect to \zjnote{the following maps}
$$K \mathcal{A} (H_0(\fchfg)) \longrightarrow K(H_0(\fchfg))
\longrightarrow K(H_0( \fcht))$$ induced from the inclusions
$\mathcal{A} \fchfg \subset \fchfg \subset \fcht$ (see
Theorem~\ref{t:computations-chains}).

Recall also from Proposition~\ref{p:LC-isos} that we have an
isomorphism of rings $\theta: \Lambda \xrightarrow{\; \cong \;} LC$
that also sends $\Lambda_P^{(0)} \subset \Lambda_P \subset \Lambda$
isomorphically to $LC'_B \subset LC_B \subset LC$. Thus, by composing
$\lambda$ with $\theta$ we obtain an $LC_B$-function
(resp. $LC$-function) $\bar{\chi}([C])$ associated to every class
$[C]$ in $ K(H_0( \fchfg))$ and $K(H_0( \fcht))$ respectively, namely
$\bar{\chi}([C]) := \theta \circ \lambda ([C])$.

As mentioned above one can define the functions $\bar{\chi}$ directly
on the level of barcodes. \zjred{Recall} from~\S\ref{subsubsec:iso-rings} that
given a graded barcode
$\mathcal{B} = \{\mathcal{B}_i\}_{i \in \mathbb{Z}}$ we define
\begin{equation} \label{eq:chi-bar-B} \bar{\chi}_{\mathcal{B}} =
  \sum_{i \in \mathbb{Z}} (-1)^i \sigma_{\mathcal{B}_i} \in LC.
\end{equation}
If $\mathcal{B}^{\text{tot}}$ contains only a finite number of bars,
then $\bar{\chi}_{\mathcal{B}} \in LC_B$ and if moreover each of these
bars is of finite-length then $\bar{\chi}_{\mathcal{B}} \in LC'_B$.
Below we will construct several numerical invariants of barcodes (or
$K$-classes of filtered chain complexes) using the functions
$\bar{\chi}$.

Some of the material
in~\S\ref{sbsb:euler-char}-~\S\ref{sbsb:gnrl-length} bears similarity
with older work of Schapira~\cite{Scha:operations} done in the
framework of the category of constructible sheaves on manifolds.
Other related, more recent, work can be found in~\cite{BB:euler,
  GH:magnitude, Leb:hybrid}.

\subsubsection{Euler characteristic} \label{sbsb:euler-char} The
classical Euler characteristic of finitely generated chain complexes
induces a well defined homomorphism
$\chi : K(H_0(\fchfg)) \longrightarrow \mathbb{Z}$ that factors
through the obvious map
$K(H_0(\fchfg)) \longrightarrow K(H_0(\ch^{\fg}))$ and the classical
Euler characteristic. Here we have denoted by $K(H_0(\ch^{\fg}))$ the
$K$-group of the homotopy category of finitely generated chain
complexes (endowed with no filtrations).

In terms of barcodes this has the following description. Suppose
$\mathcal{B}$ is stable at infinity, that is, the associated graded
persistence module has its persistence structure maps being
isomorphisms when the filtration parameter is sufficiently
large. Assume in addition that $\mathcal{B}^{\text{tot}}$ has only \zjred{
finitely many} bars of infinite length. Then the {\it Euler
  characteristic} of $\mathcal B$, denoted by $\chi(\mathcal B)$, is
equal to $\bar{\chi}_{\mathcal B}(\infty)$.

If we use the $\Lambda_P$-model for $K(H_0(\fchfg))$ then
$\chi(\mathcal B) = \lambda_{\mathcal{B}}|_{t=1}$. This coincides also
with the total Euler characteristic of the underlying chain complex of
a filtered chain complex $C$ in the sense that
$\chi(C) = \bar{\chi}([C])(\infty) = \lambda([C])|_{t=1}$.

\subsubsection{Length} \label{sbsb:length} Here we define a
homomorphism $\ell: K(H_0(\fchfg)) \longrightarrow \mathbb{R}$ which
we call length as follows. \zjred{Identify} $K(H_0(\fchfg)) \cong LC_B$ (using
$\theta \circ \lambda$ as explained at the beginning
of~\S\ref{s:num-inv}. Recall that all functions $\sigma \in LC_B$ are
bounded and moreover constant at infinity. We define:
\begin{equation} \label{eq:length-1} \ell(\sigma) = \int_{\mathbb{R}}
  \bigl( \sigma - \sigma(\infty)\mathds{1}_{[0,\infty)} \bigr) d\mu,
  \quad \forall \; \sigma \in LC_B,
\end{equation}
where $d\mu$ is the Lebesgue measure on $\mathbb{R}$. Note that the
function $\sigma - \sigma(\infty)\mathds{1}_{[0,\infty)}$ has compact
support hence its integral over $\mathbb{R}$ is finite.

If we restrict $\ell$ to $K \mathcal{A} (H_0(\fchfg)) \cong LC'_B$,
then
$$\ell(\sigma) = \int_{\mathbb{R}} \sigma d\mu,
\quad \forall \; \sigma \in LC'_B,$$ since all \zjred{functions}
$\sigma \in LC'_B$ satisfy $\sigma(\infty)=0$.

The analogous definition of length for barcodes goes as follows. Let
$\mathcal{B}$ be a graded barcode with $\mathcal{B}^{\text{tot}}$
consisting of finitely many bars. Let
$\bar{\chi}_{\mathcal{B}} \in LC_B$ be its associated $LC$-function as
defined in~\eqref{eq:chi-bar-B}.  We define the length of
$\mathcal{B}$ to be
$\ell(\mathcal{B}):= \ell(\bar{\chi}_{\mathcal{B}})$, where the latter
is defined using~\eqref{eq:length-1}. More explicitly:
\begin{equation} \label{eq:length-2} \ell(\mathcal{B}) =
  \int_{\mathbb{R}} \bigl(\bar{\chi}_{\mathcal{B}} -
  \bar{\chi}_{\mathcal{B}}(\infty) \mathds{1}_{[0,\infty)} \bigr)
  d\mu.
\end{equation}

\begin{rem} \label{r:length} (a) The reason for the name ``length'' is
  the following.  In case
  $\mathcal{B} = \{(I_j, m_j)\}_{j \in \mathcal{J}}$ is an {\em
    ungraded} barcode with finitely many bars all of which of {\em
    finite length} then $\ell(\mathcal{B})$ is the sum of the lengths
  of the bars in $\mathcal{B}$ (taking multiplicities into
  account). This is so because $\bar{\chi}_{\mathcal{B}}(\infty) = 0$,
  hence
  $$\ell(\mathcal{B}) = \int_{\mathbb{R}} \bar{\chi}_{\mathcal{B}} d\mu =
  \sum_{j \in \mathcal{J}} m_j \int_\mathbb{R} \mathds{1}_{I_j} d\mu =
  \sum_{j \in \mathcal{J}} m_j \, \text{length}(I_j).$$ Note that in
  this case $\ell(\mathcal{B}) \geq 0$ with equality if and only if
  $\mathcal{B} = \emptyset$.
    
  In case $\mathcal{B} = \{\mathcal{B}_i\}_{i \in \mathbb{Z}}$ is a
  graded barcode (with $\mathcal{B}^{\text{tot}}$ having only finitely
  many bars and all of finite length) then
  $$\ell(\mathcal{B}) = \sum_{i \in \mathbb{Z}}(-1)^i \ell(\mathcal{B}_i)$$
  is just the signed sum of the lengths of the $\mathcal{B}_i$'s. In
  contrast to the case of ungraded bars, here $\ell(\mathcal{B})$ can
  assume also negative values and it may also happen that
  $\ell(\mathcal{B})=0$ even if $\mathcal{B} \neq \emptyset$.

  \smallskip
 (b) The integrand in the \zjnote{definition} of length in~\eqref{eq:length-1}
  is very much related to \zjnote{the K-group of the acyclics in $\mathcal F{\rm {\bf Ch}}^{\rm fg}$}. Indeed, the map
  $$LC_B \longrightarrow LC'_B, \quad
  \sigma \longmapsto \sigma - \sigma(\infty) \mathds{1}_{[0,\infty)}$$
  corresponds under the isomorphisms $K(H_0(\fchfg)) \cong LC_B$,
  $K \mathcal{A}(H_0(\fchfg)) \cong LC'_B$, defined by
  $\theta \circ \lambda$, to the map
  $$K(H_0(\fchfg)) \longrightarrow K \mathcal{A}(H_0(\fchfg)),
  \quad [C] \longmapsto [C] - \chi(C)[E_1(0)]$$
  from~\eqref{eq:split-ses-ch} in Remark~\ref{r:ses-ch-splits}.  \Qed
\end{rem}

Next we give an alternative expression for $\ell$ which is useful when
working with the identification $K(H_0(\fchfg)) \cong \Lambda_P$ (via
the isomorphisms $\lambda$ as in~\eqref{eq:K-lambda}). More precisely,
we consider here
$\ell \circ \theta: \Lambda_P \longrightarrow \mathbb{R}$, where
$\theta: \Lambda_P \longrightarrow LC_B$ is the isomorphisms from
Proposition~\ref{p:LC-isos}. By abuse of notation we denote
$\ell \circ \theta$ by $\ell$.

Let $P(t) \in \Lambda_P$ be a Novikov polynomial. A simple calculation
shows that
\begin{equation} \label{eq:ell-P} \ell(P(t)) = -P'(1),
\end{equation}
where $P'(1)$ stands for the derivative of $P(t)$ with respect to the
formal variable $t$, calculated at $t=1$.

The last formula has the following consequence. Recall that
$K(H_0(\fchfg))$ has a product induced from tensor products of
filtered chain complexes $[C] \cdot [D] = [C \otimes D]$.
\begin{cor} \label{c:ell-prod} For every $X, Y \in K(H_0(\fchfg))$ and
  $n \geq 1$ we have
  \begin{equation} \label{eq:ell-prod}
    \begin{aligned}
      \ell(X \cdot Y) & = \ell(X) \chi(Y) + \chi(X) \ell(Y), \\
      \ell(X^n)  & = n\chi(X)^{n-1} \ell(X).
    \end{aligned}
  \end{equation}
\end{cor}
This follows easily from~\eqref{eq:ell-P}, the fact that $\lambda$ is
an isomorphism of rings and the chain rule for derivative of products
of functions.

Recall that $K(H_0(\fchfg))$ and $K \mathcal{A}(H_0(\fchfg))$ are both
$\Lambda_P$-modules. With respect to these structures we have the
following formula
\begin{cor} \label{c:ell-mod} For all $X \in K(H_0(\fchfg))$,
  $P(t) \in \Lambda_P$ we have:
  $$\ell(P(t)X) = P(1) \ell(X) + P'(1) \chi(X).$$
  In particular, if $X \in K \mathcal{A}(H_0(\fchfg))$ we have
  $\ell(P(t)X)) = P(1) \ell(X)$.
\end{cor}
The very last identity means that the restriction
$\ell|_{K \mathcal{A} (H_0(\fchfg))}: K \mathcal{A} (H_0(\fchfg))
\longrightarrow \mathbb{R}$ is a map of $\Lambda_P$-modules, if we
endow its target $\mathbb{R}$ with the $\Lambda_P$-module structure
$P(t) \cdot x = P(1)x$, $\forall x \in \mathbb{R}$.

\subsubsection{Generalized length} \label{sbsb:gnrl-length} For
simplicity we assume here that the graded barcodes $\mathcal{B}$ under
consideration satisfy that $\mathcal{B}^{\text{tot}}$ has only
finitely many bars of infinite length.

Given a functional $T: LC_B \longrightarrow \mathbb{R}$ we can define
an associated functional $\widetilde{\ell}_T$ on graded barcodes by
\begin{equation} \label{functional-barcode-1}
  \widetilde{\ell}_T(\mathcal{B}) = T(\bar{\chi}_{\mathcal{B}}).
\end{equation}

Similarly, given a functional $T: LC'_B \longrightarrow \mathbb{R}$ we
can define an associated functional $\ell_T$ on graded barcodes by
\begin{equation} \label{functional-barcode-2} \ell_T(\mathcal{B}) = T
  \bigl( \bar{\chi}_{\mathcal{B}} - \bar{\chi}_{\mathcal{B}}(\infty)
  \mathds{1}_{[0,\infty)} \bigr).
\end{equation}

If $h: \mathbb{R} \longrightarrow \mathbb{R}$ is locally integrable
function we can apply the above to the functional $T_h$, i.e.~the
distribution associated to $h$. We denote by
\begin{equation} \label{functional-barcode-3} \ell_h(\mathcal{B}) :=
  \int_{\R} h \cdot \bigl( \bar{\chi}_{\mathcal{B}} -
  \bar{\chi}_{\mathcal{B}}(\infty) \mathds{1}_{[0,\infty)} \bigr) d\mu
\end{equation}
the corresponding functional on barcodes. In case $h$ satisfies in
addition that $\int_{a}^{\infty} |h| d \mu < \infty$ for every
$a \in \mathbb{R}$, then
$\widetilde{\ell}_{h} := \widetilde{\ell}_{T_h}$ is also well defined.

Both the Euler characteristic and the length considered
in~\S\ref{sbsb:euler-char} and~\S\ref{sbsb:length} are special cases
of $\widetilde{\ell}_T$ and $\ell_T$ for different $T$'s. More
precisely,
$\chi(\mathcal B) = \widetilde{\ell}_{\delta_{\infty}} (\mathcal B)$
where $\delta_{\infty}$ is the Dirac delta function at $\infty$.
Similarly, the length $\ell(\mathcal{B})$ coincides with
$\ell_h(\mathcal{B})$ for the function $h \equiv 1$.

\subsection{Other measurements} \label{sb:ot-mes} Here we briefly
describe some other numerical invariants of barcodes that are not
$K$-theoretic (in the sense that the value of the invariants on a
barcode $\mathcal{B}$ does not depend solely on
$\bar{\chi}_{\mathcal{B}}$). For simplicity we assume below that our
barcodes $\mathcal{B}$ are such that $\mathcal{B}^{\text{tot}}$
contains only a finite number of bars.

Given a graded barcode
$\mathcal{B} = \{\mathcal{B}_i\}_{i \in \mathbb{Z}}$ as above,
consider the following $LC_B$ function:
\begin{equation} \label{abs-indicator-barcode} |\sigma|_{\mathcal{B}}
  := \sigma_{\mathcal{B}^{\text{tot}}} = \sum_{i \in \mathbb{Z}}
  \sigma_{\mathcal{B}_i} \in LC_B.
\end{equation}
The difference between $\bar{\chi}_{\mathcal{B}}$ and
$|\sigma|_{\mathcal{B}}$ is that in the latter we take an unsigned sum
of the functions $\sigma_{\mathcal{B}_i}$, hence we always have
$|\sigma|_{\mathcal{B}} \geq 0$ with equality if and only if
$\mathcal{B} = \emptyset$.

One can now apply different functionals, as
in~\S\ref{sbsb:gnrl-length} above to the functions
$|\sigma|_{\mathcal{B}}$ to obtain various measurements on barcodes.

For example, in case $\mathcal{B}^{\text{tot}}$ has only bars of
finite length we can define its {\em absolute total length}:
$$|\ell|(\mathcal{B}) := \int_{\mathbb{R}} |\sigma|_{\mathcal{B}} d \mu.$$
When $\mathcal{B}$ comes from the filtered Morse persistent homology,
$|\ell|(\mathcal{B})$ has been studied and applied in analysis and
spectral theory~\cite[Chapter~6]{PRSZ20}.

Another measurement is the following (here we do {\em not} assume that
$\mathcal{B}$ has no bars of infinite length).
\begin{equation} \label{bar-length-chi} \overline{\ell}(\mathcal{B}) =
  \ell(\mathcal B_{\rm finite}) + \chi(\mathcal{B}),
\end{equation}
where $\mathcal B_{\rm finite}$ is the subset of $\mathcal{B}$
consisting of all its finite-length bars, and $\chi(\mathcal{B})$ is
the Euler characteristic of $\mathcal{B}$ defined
in~\S\ref{sbsb:euler-char}.

The preceding invariant $\bar{\ell}$ behaves nicely with respect to
tensor products.  For any two graded barcodes $\mathcal B'$ and
$\mathcal B''$ one can define a product
$\mathcal B' \otimes \mathcal B''$ which is the barcode of the tensor
product of the corresponding persistence modules or filtered chain
complexes (cf.~Example \ref{ex-tensor-elementary}). We have the
following identity (see~\cite{BCZ-pkt-arxiv} for the proof):
\begin{equation} \label{bar-length-tensor} \overline{\ell}(\mathcal B'
  \otimes \mathcal B'') = \overline{\ell}(\mathcal B') \chi(\mathcal
  B'') + \chi(\mathcal B') \overline{\ell}(\mathcal B'') -
  \chi(\mathcal B') \chi(\mathcal B'').
\end{equation}

\section{Symplectic Applications} \label{sec-app}

\zjnote{The purpose of this} section is to consider the persistence
$K$-theory developed earlier in the paper for the
\zjnote{triangulated} persistence category $\C= \C\fuk(\mathcal{X})$
associated to a filtered Fukaya category,
$\fuk(\mathcal{X},\mathscr{P})$, as in \S\ref{subsec:Fuk-TPC}, and
deduce a couple of symplectic applications.

Both applications identify algebraic properties that differentiate
embedded Lagrangians from immersed ones. The first one,
in~\S\ref{sb:k-lag}, Corollary \ref{thm-emb-char} (which is a direct
consequence of the $K$-theoretic constructions in the paper), shows
that the Euler type pairing $\kappa$ from \S\ref{sec-pairing}, has a
special behavior on classes $A\in K(\C)$ that admit embedded
representatives in the sense that $\kappa (A,A)$ is a constant
polynomial, equal to the Euler \zjnote{characteristic} of the
representative. \zjnote{As it turns out this fails for some classes
  that represent Lagrangian immersions}. The second
\zjnote{application}, in Theorem \ref{thm:dist-immersed}, provides a
lower bound \zjnote{on} the ``energy'' required to resolve the
singularities of \zjnote{Lagrangian immersions} that represent classes
$A$ with $\kappa(A,A)\in\La_P$ non-constant.  In
\S\ref{subsec:immersed} \zjnote{we discuss} some relations
\zjnote{between} the constructions above and immersed Floer
\jznote{theory}. In particular, how the expression of $\kappa(A,A)$ is
expected to provide a lower bound for the number of self-intersection
points of an immersed representative of $A$. This \zjnote{subsection}
is somewhat speculative in the sense that we do not fully develop in
this paper the \zjnote{theory of} filtered Fukaya \zjnote{category} in
the immersed case. In \S\ref{sb:exp} we provide an example to
illustrate \pbred{our} results.

\subsection{The pairing $\kappa$ and Lagrangian
  submanifolds} \label{sb:k-lag}

We assume here the setting and notation
from~\S\ref{subsubsec:filtr-fuk}. In particular, we work over the
\pbred{base} field $\Z_{2}$, \pbred{and we use the homological and
  grading conventions for Fukaya categories as described at the
  beginning of~\S\ref{subsubsec:filtr-fuk}.}  Thus \pbred{our ambient
  manifold} \zjnote{$(X^{2n},\omega = d\lambda)$ is a Liouville
  manifold.} The objects of the category
$\fuk(\mathcal{X},\mathscr{P})$ are a family $\mathcal{X}$ of marked
exact Lagrangians in $X$ - these are triples
$(\bar{L},h_{L},\theta_{L})$ where $L\subset X$ is a closed (embedded)
Lagrangian submanifold, $h_{L}$ is a primitive of $\la|_{L}$,
$\theta_{L}$ is a choice of grading. We assume $\mathcal{X}$ to be
closed under shifts of primitives and translation of grading. Finally,
$\mathscr{P}$ is a choice of perturbation \zjnote{data} such that the
resulting \zjnote{$A_{\infty}$-category}
$\fuk(\mathcal{X},\mathscr{P})$ is filtered.  The triangulated
persistence category $\C\fuk(\mathcal{X})$ is associated to
$\fuk(\mathcal{X};\mathscr{P})$ by an algebraic construction described
in \S\ref{subsubsec:ainfty-alg} and \S\ref{subsubsec:filtr-fuk}.  In
brief, this TPC is the homological category
$$\C\fuk(\mathcal{X}) = H_{0}[\fuk(\mathcal{X};\mathscr{P})^{\nabla}]$$
where
$$\fuk(\mathcal{X};\mathscr{P})^{\nabla} \subset
Fmod(\fuk(\mathcal{X};\mathscr{P}))$$ is \zjred{the} closure inside the
filtered dg-category $Fmod(\fuk(\mathcal{X};\mathscr{P}))$ of filtered
$\fuk(\mathcal{X};\mathscr{P})$- modules of the (filtered) Yoneda
modules $\mathcal{Y}(L)$, $L\in \mathcal{X}$, with respect to iterated
cones of weight $0$ and with respect to $r$-isomorphism, for all
$r\geq 0$. More details appear in \S\ref{subsec:Fuk-TPC} and we refer
to \cite{BCZ23} for the full construction.


Before we go on, let us mention the following two points that will be
useful later. Firstly, if $C, C'$ are two general objects in
$\mathcal{C}\fuk(\mathcal{X})$, then it might be difficult
to explicitly write out $\hom(C, C')$ in this category. However, when
$C = L$ and $C' = L'$ are two objects from $\mathcal{X}$ (i.e.~two
marked Lagrangians from the family $\mathcal{X}$), then we have:
\pbred{
\begin{equation} \label{hom-emb-mod}
  \hom_{\mathcal{C}\fuk(\mathcal{X})}(L, L') = \HHF_0(L, L')
\end{equation}
}
as persistence modules. More generally, again as persistence 
modules we have:
\pbred{$$ \hom_{\mathcal{C}\fuk(\mathcal{X})}(L, L'[k])= \HHF_{k}(L,L')~.~$$}
Recall that, by definition,
$T^{k}L'=L'[-k]= (\bar{L}',h_{L'}, \theta_{L'}- k)$.

\pbred{
Recall also from \S\ref{subsubsec:filtr-fuk} the second TPC,
$\C'\fuk(\mathcal X)=H_{0}[Tw(\fuk(\mathcal{X};\mathscr{P})]$, which
is the homological category of the filtered twisted complexes over
$\fuk(\mathcal{X};\mathscr{P})$ as well as the embedding:
$$\bar{\Theta}: \mathcal{C}'\fuk(\mathcal{X})
\hookrightarrow \mathcal{C}\fuk(\mathcal{X})~.~$$}

The second remark is that for any
\pbred{$C, C' \in \Ob(\mathcal{C}'\fuk(\mathcal{X}))$} there exists
$k_0 \in \mathbb{N}$ such that
\pbred{
\begin{equation} \label{symp-bd}
  \hom_{\mathcal{C}'\fuk(\mathcal{X})}(C, T^iC') = 0, \;
  \text{whenever } |i| \geq k_0~.~
\end{equation}
}

Finally, we denote by $\fuk(X)$ the (unfiltered) Fukaya category of
all \zjnote{exact} marked Lagrangians in $X$ and by $D\fuk(X)$ the
associated derived category.  In both cases these notions neglect
filtrations and if this aspect needs to be emphasized we write instead
$\fuk_{un}(X)$, $D\fuk_{un}(X)$ and similarly for other relevant
algebraic structures.

A last notation to be recalled here is the Novikov polynomial
$\lambda_M \in \Lambda_P$ associated to a persistence module $M$
(whose barcode has finitely many bars) as defined
in~\S\ref{subsec:barfuc}, see formula~\eqref{eq:lambda-M} there.

\begin{lemma} \label{lem-geo-dual} Let $L, L' \in \mathcal{X}$, viewed
  as objects of
  $\C=\mathcal{C}\fuk(\mathcal{X})$. Then
  \begin{equation} \label{eq:lambda-LL'}
    \lambda_{\hom_{\mathcal{C}}(L, L')}(t) =
    \lambda_{\hom_{\mathcal{C}}(L', L[-n])}(t^{-1})
  \end{equation}
  in $\Lambda_P$. Moreover, for every $i \in \mathbb{Z}$, the Novikov
  polynomial~\eqref{eq:lambda-LL'} for $L'=L[i]$ is constant and
  equals to:
  $$\lambda_{\hom_{\mathcal{C}}(L, L[-i])}(t) = \dim_{\mathbb{Z}_2}
  H_{n-i}(L; \mathbb{Z}_2).$$
\end{lemma}

\begin{proof}
  The first formula follows from duality in Floer homology. Namely,
  $HF_k(L,L') \cong HF_{n-k}(L',L)^*$ as persistence modules,
  \pbred{which with our grading conventions
    from~\S\ref{subsubsec:filtr-fuk} reads
    $\HHF_j(L,L') \cong \HHF_{-n-j}(L',L)^*$.}  This together with the
  discussion in~\S\ref{sb:skew-sym-kappa} relating the barcodes of a
  persistence module and its dual implies~\eqref{eq:lambda-LL'}.

  The second statement follows from the way we defined the filtrations
  on the Floer complexes in the category
  $\fuk(\mathcal{X};\mathscr{P})$, as recalled
  in~\S\ref{subsubsec:filtr-fuk}, \pbred{namely $CF_*(L,L)$ coincides
    with the Morse complex $CM_{*}(f_L)$ of a Morse function
    $f_{L}: L \to \R$, with all generators in action level $0$. With
    our grading conventions we now have
    $$\CCF_j(L,L[-i]) = CF_{n+j-i}(L,L) = CM_{n+j-i}(f_L),$$ hence
    $\HHF_0(L,L[-i]) \cong H_{n-i}(L; \mathbb{Z}_2)$.}
\end{proof}

\pbred{Denote $\mathcal{C}' := \mathcal{C}' \fuk(\mathcal{X})$.}
\pbred{
Recall from~\S\ref{sec-pairing} that we have two pairings:
\begin{equation*}
  \begin{aligned}
    & \kappa': \overline{K(\C')} \otimes_{\Lambda_{P}} K(\C')
    \longrightarrow \Lambda_P,  \\
    & \kappa: \overline{K(\C)} \otimes_{\Lambda_{P}} K(\C)
    \longrightarrow \Lambda,
  \end{aligned}
\end{equation*}
and $\kappa'$ is compatible with $\kappa$ via the map
$K(\mathcal{C}') \longrightarrow K(\mathcal{C})$, induced by the
inclusion $\mathcal{C}' \longrightarrow \mathcal{C}$. Therefore we
will denote below both pairings by $\kappa$.}

Due to the boundedness
property~\eqref{symp-bd} of \pbred{the category $\mathcal{C}'$} the
values of the pairing $\kappa$ are in $\Lambda_P$ and not just
$\Lambda$.  Although $\kappa$ is defined on the tensor product
\pbred{$\overline{K(\C')} \otimes_{\Lambda_{P}} K(\C')$} we will often
write $\kappa(-,-)$ viewing it as a bilinear map.  This pairing can be
explicitly written out as follows
\begin{equation} \label{eq:pairing-symp} \kappa([C], [C']) = \sum_{i}
  (-1)^i \lambda_{\hom_{\mathcal{C}}(C, T^iC')}(t),
\end{equation}
for any two elements \pbred{$[C], [C'] \in K(\mathcal{C}')$.} Notice
also that, according to our conventions, $T^iC' = C'[-i]$. We now
have.

\begin{cor} \label{thm-emb-char} For every two objects $C, C'$ of
  $\mathcal{C}'\fuk(\mathcal{X})$ we have
  $$\kappa([C], [C'])(t) = (-1)^n \kappa([C'], [C])(t^{-1}).$$
  In particular,
  $\kappa([C], [C])(t) = (-1)^n \kappa([C],[C])(t^{-1})$.  Moreover if
  a class $A \in K[\C\fuk(\mathcal{X})]$ can be represented by an embedded
  Lagrangian $L \in \mathcal{X}$ then
  $$\kappa(A,A) = (-1)^n \chi(\bar{L}),$$
  where $\chi(\bar{L})$ is the Euler characteristic of the underlying
  Lagrangian $\bar{L}$ of $L$.
\end{cor}

\begin{proof} Recall that the objects of $\C'\fuk(\mathcal{X})$ can be
  viewed as (filtered) twisted complexes of the Lagrangians in
  $\mathcal{X}$.  Equivalently, they are iterated cones of weight $0$
  of the Yoneda modules $\mathcal{Y}(L)$, $L\in \mathcal{X}$. The
  proof now follows by linearity from Lemma~\ref{lem-geo-dual}
  and~\eqref{eq:pairing-symp}.
\end{proof}

\begin{rem}
  (i) It is not expected for the relation in Corollary
  \ref{thm-emb-char} to hold for all the objects in
  $\C\fuk(\mathcal{X})$.  Each such object is $r$-isomorphic to an
  object of $\C'\fuk(\mathcal{X})$ for some $r\geq 0$ but if $r>0$
  this is insufficient to deduce the relation in the statement.
 
  (ii) Since the Euler characteristic of closed odd-dimensional
  manifolds vanishes, the last formula in Corollary~\ref{thm-emb-char}
  could just as well be written as $\kappa(A,A) = \chi(\bar{L})$.
    
  (iii) The self $\kappa$-pairing $\kappa([C],[C])$ of an element
  $[C] \in K(\mathcal{C})$, where $C \in \Ob(\C'\fuk(\mathcal{X}))$,
  is invariant both under action shifts as well as translations of
  $C$, in the sense that
  $$\kappa([\Sigma^r T^k C], [\Sigma^r T^k C]) = \kappa([C],[C])~,~$$
  for all $r \in \mathbb{R}$, $k \in \mathbb{Z}$, where $\Sigma^r$ and
  $T$ are the shift and translation functors, respectively.
\end{rem}

\medskip

The core of Corollary \ref{thm-emb-char} is that if a class
$A\in K(\C)$ can be represented by an \zjnote{embedded Lagrangian},
then the Euler type Novikov polynomial $\kappa(A,A)$ is
constant. \zjnote{An interesting question is to consider cases when}
$\kappa(A,A)$ is not constant and try to estimate how far from
embeddings are representatives of $A$.

We will discuss below two answers to this question, the first
interprets ``how far'' in the sense of the measurements introduced in
relation to the fragmentation pseudo-metrics discussed in
\S\ref{subsec:tr-weight}.  The second is concerned with lower bounds
for the minimal number of self-intersection points of
\zjnote{possible} immersed representatives of $A$, and is only briefly
covered in Remark \ref{rem:imme-no}.

To formulate the first result we need some further notation. The first
is the notion of {\em gap}, $\mathrm{gap}(Q)$, of a polynomial
$Q\in \La_{P}$.  Assuming that
$Q=\sum_{i=1}^{k} a_{i} t^{\alpha_{i}}$, $\alpha_{i}<\alpha_{i+1}$
\zjnote{and $a_i \neq 0$ for all $i$}, \pbred{we define}
$$\mathrm{gap}(Q)=\min \{\alpha_{i+1}-\max\{\alpha_{i},0\}
\ |\ 1\leq i \leq k-1, \ \alpha_{i+1}\geq 0 \} \in \R_{\geq 0}$$
\zjnote{while if \pbnote{$\alpha_j < 0$} for all $j$, we simply}
\pbred{set $\mathrm{gap}(Q) := 0$.}

The second notation applies to any persistence category $\C$ 
with the property that the persistence modules $\hom_{C}(A,B)$
have barcodes with finitely many bars, for all objects $A, B$. In this case,
assume that $\mathcal{K}\subset \mathrm{Obj}(\C)$ is a finite family
of objects. We let
$$b^{\#}(\mathcal{K})= \max \{ \#(\mathcal{B}_{\hom_{\C}(A,B)})
\ | \ A, B\in \mathcal{K} \}~.~$$ Here $\#(\mathcal{B}_{M})$ is the
number of bars of the barcode $\mathcal{B}_{M}$ associated to the
persistence module $M$. We call this the {\em max bar count} of
$\mathcal{K}$.

In case $\C$ is a TPC, recall that we have the graded persistence
modules $\ho(X,T^{\bullet}Y)$ defined in (\ref{eq:pers-grad-mor}) by
$\ho(X,T^{\bullet}Y)_{i} = \ho_{\mathcal{C}}(X,T^iY), \; i \in
\mathbb{Z}$. We define the corresponding {\em graded} max bar count of
$\mathcal{K}$ by \pbred{
\begin{equation*}
  b^{\#}_{\max}(\mathcal{K}) =
  \max\{ \#(\mathcal{B}_{\hom_{\C}(A,T^{i}B)})
  \ | \ A, B\in \mathcal{K}, \, i \in \Z\}.
\end{equation*}
}


We will need now a measurement related to reduced
decompositions. Recall from~\S\ref{subsec:tr-weight} that such
decompositions $D = (\bar{D},\phi)$ have a linearization $\ell(D)$ and
a weight which we denote here by $\bar{w}(D)$. Let
$\mathcal{F} \subset \Ob(\C)$ invariant \zjnote{under} shift and
translation. Assume that $D$ is a reduced decomposition in $\C$ of
$N\in \mathrm{Obj}(\C)$, with linearization
$\ell({{D}})=(X_{1},X_{2},\ldots, X_{n})$, where each
$X_{i}\in \mathcal{F}$.  We fix a bit more notation that only depends
on $\ell(D)$:
\begin{itemize}
\item[-] \pbred{$[D] = [X_1] + \cdots + [X_n] \in K(\C)$.}
\item[-] $\{\ell(D)\}=\{X_{1},X_{2},\ldots, X_{n}\}$, this is the {\em
    set} of objects in the linearization of $D$.
\item[-] $n_{D}=n$, the number of terms in the linearization of $D$.
\end{itemize}

We then let 
\zjnote{\begin{equation} \label{eq:weigthed-R}
\mathscr{R}^{\F}(N) =\inf \left\{n_{D}^{2}\cdot \frac{\bar{w}(D)\, b^{\#}_{\max}( \{\ell(D)\})}{\mathrm{gap}(\kappa([D],[D]))}\  \bigg| \, \begin{array}{l} \mbox{$D$ is a reduced decomposition of $N$} \\ \mbox{such that ${\rm gap}(\kappa([D],[D]))\neq 0$} \end{array} \right\}~.~
\end{equation}}
In case no reduced decompositions $D$ as above exist we put
$\mathscr{R}^{\mathcal{F}}(N)=\infty$.

\

This measurement can be defined for any TPC $\C$ but we will apply it
to the triangulated persistence category $\C=\C\fuk(\mathcal{X})$ and
for $\mathcal{F}\subset \mathcal{X}$.  \

We first state and prove the result and we will discuss its meaning in
\S\ref{subsec:immersed}.
\begin{thm}\label{thm:dist-immersed} \zjnote{Let} $\C=\C\fuk(\mathcal{X})$ and $\mathcal{F}\subset \mathcal{X}$ \zjnote{be} as above. If $N\in \mathrm{Obj}(\C)$
  is the Yoneda module of an embedded Lagrangian $\in \mathcal{X}$,
  then:
  $$\mathscr{R}^{\F}(N)\geq \frac{1}{4}.$$
\end{thm}

We will give an example illustrating the theorem in
\S\ref{sb:exp}. This result should be interpreted in the sense that,
if $N$ admits a reduced decomposition $D=(\bar{D},\phi)$, as in the
definition of $\mathscr{R}^{\mathcal{F}}$, then the $s$-isomorphism
$\phi$ is required to satisfy
\zjnote{$$s\geq \frac{1}{4n^{2}_{D}}\cdot
  \frac{\mathrm{gap}(\kappa([D],[D]))}{b_{\max}^{\#}(\{\ell(D)\})}$$}
for $N$ to represent an embedded Lagrangian. The reason why this is of
use is that the right hand side of this inequality only depends on the
linearization of $D$. It can be computed easily from knowledge of
$HF(F,F')$ as a persistence module for each $F,F'\in \mathcal{F}$. In
this case the max bar count is obvious and the $\mathrm{gap}$ term is
given by easy linear algebra calculations in $K(\C\fuk(\mathcal{X}))$.
In many examples, $\mathcal{F}$ forms a set of generators of
$[\C\fuk(\mathcal{X})]_{\infty}$ and the groups $HF(F,F')$ are
understood.
 


\begin{proof}[Proof of Theorem \ref{thm:dist-immersed}] Assume that
  $r > \mathscr{R}^{\mathcal{F}}(N)$. By inspecting the definition of
  reduced decompositions \zjnote{in \S\ref{subsec:tr-weight}} we
  deduce that there are exact triangles $\Delta_{i}$ in $\C_{0}$,
  $1\leq i \leq k$:
  \begin{equation}\label{eq:iterated-cones}
    \Delta_{i} \ :\ \ F_{i}\longrightarrow Y_{i-1}
    \longrightarrow Y_{i}\longrightarrow TF_{i}
  \end{equation}
  with $Y_{0}=0$, and an $r'$-isomorphism $\phi :Y_{k}\to N$ such that
  \zjnote{$$r> k^{2}\cdot \frac{r' \
      b_{\max}^{\#}(\{\ell(D)\})}{\mathrm{gap}(\kappa([D],[D]))}$$}
  where \zjnote{$[D]=\sum_i [F_{i}]$ and
    $b_{\max}^{\#}(\{\ell(D)\})=b_{\max}^{\#}(\{F_{1},F_{2}, \ldots,
    F_{k}\})$}. We denote \zjnote{$b=b_{\max}^{\#}\{\ell(D)\}$} and
  $g=\mathrm{gap}(\kappa([D],[D]))$.  We intend to show that, because
  $N$ is embedded, we have
  \begin{equation}\label{eq:ineq-to} k^{2}br'\geq \frac{1}{4}\  g 
  \end{equation}
  which implies the statement.

  The first step is to consider the two graded persistence modules
  \pbred{$$\HHF_{\bullet}(N,N)=\hom_{\C}(N,T^{\bullet}N)\ \
    \mathrm{and}\ \ \hom_{\C}(Y_{k}, T^{\bullet}Y_{k})$$} and notice
  that they are $2r'$-interleaved. This is because, by \cite[Lemma
  2.83]{BCZ23}, the two modules $N$ and $Y_{k}$ are $r'$-interleaved
  as objects of the triangulated persistence category $\C$ and this
  implies that $\hom_{\C}(N,T^{\bullet}Y_{k})$ is $r'$-interleaved
  with respect to both \pbred{$\HHF_{\bullet}(N,N)$} and
  $\hom_{\C}(Y_{k},T^{\bullet}Y_{k})$ (as persistence modules) which
  implies our claim.  This means, by the isometry Theorem in
  \cite{PRSZ20}, that the bottleneck distance between the two barcodes
  \pbred{${\mathcal B_{1}:=\mathcal{B}_{\HHF_{\bullet}(N,N)}}$ and
    ${\mathcal
      B_{2}:=\mathcal{B}_{\hom_{\C}(Y_{k},T^{\bullet}Y_{k})}}$} is at
  most $2r'$.  Recall from Lemma \ref{lem-geo-dual} that
  \zjnote{$\mathcal B_{1}$} consists of only infinite intervals of the
  form $[0,\infty)$. This implies that all the finite intervals in
  \zjnote{$\mathcal B_{2}$} are of length at most $r'$ and all the
  infinite intervals in \zjnote{$\mathcal B_{2}$} are of the form
  $[x,\infty)$ with $x\in[-2r',2r']$.

  The second step is to notice that
  $-[Y_k]=[F_{1}]+\ldots + [F_{k}]=[D]$.  Thus
  \zjnote{$$\kappa([D],[D])=
    \la_{\hom_{\C}(Y_{k},T^{\bullet}Y_{k})}~.~$$} In particular,
  $\mathrm{gap}(\kappa([D],[D])) =
  \mathrm{gap}(\la_{\hom_{\C}(Y_{k},T^{\bullet}Y_{k})})$. For further
  use, let
  \zjnote{$$P= \la_{\hom_{\C}(Y_{k},T^{\bullet}Y_{k})}
    =\sum_{i=1}^{m}a_{i}t^{\alpha_{i}}$$} and fix
  $\alpha'_{i} <\alpha_{i+1}$ such that
  $\mathrm{gap}(P)=\alpha_{i+1}-\alpha'_{i}$ with $\alpha_{i+1}\geq 0$
  and $\alpha'_{i}=\max\{\alpha_{i},0\}$.

  The third step is an application of the persistence Morse
  inequalities from \S\ref{subsubsec:morse} and, in particular, the
  bar counting estimate from Corollary \ref{cor:count-bars}. That
  estimate implies that:
  \begin{lem}\label{lem:bars-cts}  We have 
    $$\#(B_{2})\leq k^{2}b$$ (where, as above,
    $B_{2}=\mathcal{B}_{\hom_{\C}(Y_{k},T^{\bullet}Y_{k})}$).
  \end{lem}

\begin{proof}[Proof of Lemma \ref{lem:bars-cts}] Recall that
  \zjnote{$b=b_{\max}^{\#}(\{F_{1},F_{2}, \ldots, F_{k}\})$}. By the
  definition of \zjnote{$b^{\#}_{\max}$}, we have that $b$ is the
  maximal number of bars in the \zjnote{barcodes} of the graded
  persistence modules \zjnote{$HF_*(F_{i}, F_{j})$} for
  $1\leq i,j\leq k$ \zjnote{and $\ast \in \Z$}. To prove the statement
  we will use Corollary \ref{cor:count-bars} and apply it repeatedly
  to the exact sequences from (\ref{eq:iterated-cones}). To simplify
  notation we will \zjnote{put
    $b(M,N)=\# \mathcal{B}_{\hom_{\C}(M,T^{\bullet}N)}$} where
  $M,N\in \mathrm{Obj}(\C)$.  For each sequence from
  (\ref{eq:iterated-cones}), and additional module $M$, we deduce from
  Corollary \ref{cor:count-bars}:

  $$b(M,Y_{i})\leq b(M,Y_{i-1}) +b(M,F_{i})
  \ , \ b(Y_{i},M)\leq b(Y_{i-1},M) +b(F_{i},M)~.~$$ \zjnote{where the
    second inequality can be deduced from the first one after applying
    the duality as in (\ref{equ-hom-duality}) in Section
    \ref{subsec:pairing}.} It is clear that, by taking $M=F_{s}$, we
  obtain:
  $$b(F_{s}, Y_{i})\leq ib \ , \ b(Y_{i},F_{s})\leq ib~.~$$
  Now suppose that $b(Y_{q},Y_{q})\leq q^{2}b$. We want to show
  $$b(Y_{q+1},Y_{q+1})\leq (q+1)^{2}b$$ as this concludes the proof of the lemma, by \zjnote{induction}. 
  We have $$b(Y_{q},Y_{q+1})\leq b(Y_{q},Y_{q}) + b(Y_{q}, F_{q+1})\leq q^{2}b+ qb$$ and further
\zjnote{\begin{align*}
b(Y_{q+1},Y_{q+1}) & \leq b(Y_{q}, Y_{q+1}) + b(F_{q+1},Y_{q+1})\\
& \leq q(q+1)b+ (q+1)b=(q+1)^{2}b.
\end{align*}}
\pbred{This concludes the proof of the lemma.}
\end{proof}

We now return to the proof of the theorem for the final part of the
argument.  \pbred{Assume} that $$ 4 k^{2} b r' < g$$ and \pbred{we
  will show} that this leads to a contradiction.  Here $k$ is at least
$1$ and $b$ is at least $2$.  This is because the max bar count of a
family of exact Lagrangians is at least the number of bars in
$\mathcal{B}_{HF(L, L)}$ for some $L\in \mathcal{X}$ and this contains
at least $2$ bars.  Thus $$g\geq 8r'~.~$$ Recall that there \zjred{is a
total of} \zjnote{ $k^2b$-many bars in
  $\mathcal B_{2}=\mathcal{B}_{\hom_{\C}(Y_{k},T^{\bullet}Y_{k})}$},
and each finite one is of length at most $r'$.  We will call the
values of the ends of the intervals in \zjnote{$\mathcal B_{2}$} the
spectral values of $\mathcal B_{2}$.

We focus on the spectral values that belong to
$H_{\delta}=[\alpha_{i+1}-\delta, \alpha_{i+1}+\delta]$ with
$0\leq \delta\leq \frac{g}{2}$.  We have
$\alpha_{i+1}-\delta-r'\geq g/2- r'\geq 3r'$.  As a result, only the
ends of finite intervals can be in $H_{\delta +r'}$ (recall that the
infinite intervals $[x,\infty)$ have $x\in [-2r',2r']$). \zjnote{Also,
  if} a finite interval $J\in B_{2}$ has an end in $H_{\delta}$, then
$J\subset H_{\delta+r'}$.  For each $H_{\delta}$ we consider the
polynomial $P'_{\delta}$ defined as in formula (\ref{eq:lambda-M-blt})
but taking into account only the intervals that have at least one end
in $H_{\delta}$.  This means that if
$B_{2}=\bigoplus_{j} I_{j}\oplus \bigoplus_{h}J_{h}$ where $I_{j}$ are
finite intervals $I_{j}=[u_{j}, v_{j})$ and $J_{h}$ are the infinite
intervals, then:
\zjnote{$$P'_{\delta}=\sum_{j} (-1)^{q_{j}}\epsilon_{\delta}(
  I_{j})(t^{u_{j}}-t^{v_{j}})~.~$$} Here $q_{j}$ is the degree of
$I_{j}$ \zjnote{in the graded barcode $\mathcal B_2$} and
$\epsilon_{\delta}(I_{j})=1$ if at least one of $u,v$ is in
$H_{\delta}$ and $\epsilon_{\delta}(I_{j})=0$ otherwise.  Because only
finite intervals contribute to $P'_{\delta}$, we have
$P'_{\delta}(1)=0$. Notice that if we truncate $P'_{\delta}$ to a
polynomial $P_{\delta}$ consisting of only the terms of $P'_{\delta}$
of the form $a t^{\alpha}$, $a\in \Z$ and with $\alpha\in H_{\delta}$,
then $P_{\delta}=a_{i}t^{\alpha_{i+1}}$.  It is clear that the number
of terms in $P'_{\delta}$ is at least that of $P_{\delta}$. The fact
that $P_{\delta}(1)=a_{i+1}$ and $P'_{\delta}(1)=0$ implies that
$P'_{\delta}$ has strictly more terms than $P_{\delta}$\zjred{, where
  within these extra terms, there is at least one term in the form of
  $at^{\alpha}$ with
  $\alpha \in [\alpha_i+\delta, \alpha_i+\delta +r')$ or
  $\alpha \in (\alpha_i-\delta-r', \alpha_i-\delta]$}. \zjred{Since
  $P_{\delta+r'}(1)=a_{i+1}$, this extra term $at^{\alpha}$}
\pbnote{has to be cancelled in} \zjred{$P'_{\delta + r'}$, which}
implies that: \zjnote{\begin{equation}\label{eq:indu-count} \#
    \mathcal B_2(H_{\delta +r'}) > \# \mathcal B_2(H_{\delta})
  \end{equation}}
where \zjnote{$\# \mathcal B_2(H_{\delta})$ denotes the number
  of the intervals in $\mathcal B_{2}$} that have at least  
one end in $H_{\delta}$.
Given that  $k^{2}br'\leq \frac{g}{2}$ we apply  (\ref{eq:indu-count}) to 
deduce the sequence of strict inequalities:
\zjnote{$$\# \mathcal B_2(H_{k^{2}br'+r'})
  >\# \mathcal B_2(H_{k^{2}br'}) > \ldots
  >\#\mathcal B_2(H_{r'}) > \#\mathcal B_2(H_{0}) \geq 1~.~$$} 
The existence of this sequence contradicts the fact that
the total number of intervals in \zjnote{$\mathcal B_{2}$} is at most $k^{2}b$.
\end{proof}

\subsection{Immersed Lagrangians}\label{subsec:immersed}
We will explain here how the statement of Theorem
\ref{thm:dist-immersed} provides an estimate on how much energy is
needed to resolve the singularities of certain immersed Lagrangians.

The immersed Lagrangians in question are (unobstructed) marked
immersed Lagrangians that geometrically coincide with a union of
elements in $\bar{\mathcal{F}}$.  To be more precise, let
$\bar{F}_{1},\ldots, \bar{F}_{k}\subset \bar{\mathcal{F}}$.  Recall
that $\bar{\mathcal{X}}$ are the elements of $\mathcal{X}$ with the
primitives and gradings forgotten, and \zjnote{the} same for
$\bar{\mathcal{F}}$. For $L\in \mathcal{X}$ we denote by $\bar{L}$ the
corresponding element of $\bar{\mathcal{X}}$. By hypothesis the family
$\bar{\mathcal{X}}$ is in general position.  Thus, geometrically, the
union $\bar{\mathcal{L}}=\cup_{i}\bar{F}_{i}$ has only double points
(to simplify the discussion \zjnote{we assume} that there are no
repetitions among the $\bar{F}_{i}$'s). A marked immersed Lagrangian -
in the sense of \cite{Bi-Co:LagPict} - with underlying geometric
Lagrangian $\bar{\mathcal{L}}$ - is a quadruple
$\mathcal{L}=(\bar{\mathcal{L}}, h_{\mathcal{L}},
\theta_{\mathcal{L}}, I_{\mathcal{L}})$ where $h_{\mathcal{L}}$,
$\theta_{\mathcal{L}}$ are a primitive and a grading, as in the
embedded case, and $I_{\mathcal{L}}$ is a subset of the
self-intersection points of $\bar{\mathcal{L}}$. Under certain
conditions imposed on $I_{\mathcal{L}}$ - see again
\cite{Bi-Co:LagPict} where this class of immersed Lagrangians is
discussed in detail - such a Lagrangian represents a twisted complex
$T_{\mathcal{L}}\in\C'\fuk(\mathcal{X})$ of the
form
$$T_{\mathcal{L}}=(\
\bigoplus^{k}_{i=1}\Sigma^{\alpha_{i}}F_{i}[k_{i}]\ ,\ \delta \ )$$
for some differential $\delta$.  The $\alpha_{i}$ and $k_{i}$ here are
determined by the choice of $I_{\mathcal{L}}$ and the initial
primitives and gradings on the $F_{i}$'s.

A resolution of the singularities of $\bar{\mathcal{L}}$ is a choice
of a set $I_{\mathcal{L}}$ of points of self-intersection of
$\mathcal{L}$ (as well as of primitive and grading) and an embedded
Lagrangian $N\subset \mathcal{X}$ together with an $r$-isomorphism in
$\C\fuk(\mathcal{X})$:
$$\phi:T_{\mathcal{L}}\to N~.~$$
By unwrapping the various definitions, the statement of the corollary
implies that
\zjnote{$$r \geq \frac{1}{4 k^{2}}\cdot \frac{\mathrm{gap}(\kappa
    (\oplus_{i} [\Sigma^{\alpha_{i}}F_{i}] , \oplus_{i}
    [\Sigma^{\alpha_{i}}F_{i}]))}{b_{\max}^{\#}(\{F_{1},\ldots,
    F_{k}\})}~.~$$}

As mentioned before, the right side of this inequality is often
approachable.  \

Of course, neglecting any reference to immersed Lagrangians, we can
talk directly about resolving the singularities of a twisted complex
$T\in \C'\fuk(\mathcal{X})$ of the form
$T=(\oplus_{i=1}^{k}\Sigma^{\alpha_{i}}F_{i}[k_{i}], \delta)$ and the
estimate above applies.  This fits very well with the notion of
reduced decomposition that was used before in the definition of the
number $\mathscr{R}^{\mathcal{F}}$.

\begin{rem}\label{rem:imme-no} (i) One geometric source for a
  resolution of the singularities of $\bar{\mathcal{L}}$ in the sense
  above is a surgery of $\bar{\mathcal{L}}$ at the points in
  $I_{\mathcal{L}}$, followed by a Hamiltonian homotopy. This is not
  quite so immediate to show because once the surgery at only {\em
    some} of the self intersection points of $\bar{\mathcal{L}}$ is
  performed, one is left with an immersed Lagrangian that \zjnote{may
    be} unobstructed but no longer a union of embedded components and,
  thus, no longer approachable with the Fukaya type machinery
  discussed till now. It is expected that marked immersed Lagrangians
  that are unobstructed can be included in a family like $\mathcal{X}$
  above and used as objects of a filtered Fukaya category
  $\fuk^{imm}(\mathcal{X};\mathscr{P})$. The construction of this
  $\fuk^{imm}(\mathcal{X};\mathscr{P})$, such as to allow in
  $\mathcal{X}$ marked, immersed, unobstructed Lagrangians that are
  not unions of embedded components, is somewhat delicate. The
  non-filtered version appears in \cite{Bi-Co:LagPict} but the
  \zjnote{adjustment} of the filtered construction from \cite[Chapter
  3]{BCZ23} requires some additional effort that will not be pursued
  here.

  (ii) Assuming $\fuk^{imm}(\mathcal{X};\mathscr{P})$ defined, we
  briefly mention a couple of properties of the pairing $\kappa$ for
  the resulting TPC, $\C\fuk^{imm}(\mathcal{X})$.  If a class
  $A\in K(\C\fuk^{imm}(\mathcal{X}))$ is represented by an immersed
  Lagrangian $L$ (with only transverse double points), then (compare
  with Corollary \ref{thm-emb-char}):
  \begin{itemize}
  \item[(a)] $\kappa(A, A)(t) = (-1)^n \kappa(A,A)(t^{-1})$; 
  \item[(b)] the number of self intersection points of $L$ has as
    lower bound $N^{+}[\kappa(A,A)]$ where for a polynomial
    $P= \sum k_{i}t^{\alpha_{i}}\in \La_{P}$, $k_{i}\in \Z$,
    $\alpha_{i}\in\R$, we have $N^{+}[P]=\sum_{\alpha_{i}>0} |k_{i}|$.
  \end{itemize}
  Both properties follow from the cluster, or pearly, model of Floer
  homology for immersed, unobstructed Lagrangians $\mathcal{L}$.  In
  this model the generators of the Floer complex
  $CF(\mathcal{L},\mathcal{L})$ are of two types: critical points of a
  Morse function on the domain of the immersion and self-intersection
  points of the immersion. Each self-intersection point appears twice,
  with opposite action levels \zjnote{(see \cite{Alston-Bao:imm2},
    \cite{Bi-Co:LagPict})}.

\end{rem}

\subsection{An example} \label{sb:exp} The aim here is to illustrate
how Corollary \ref{thm-emb-char} and Theorem \ref{thm:dist-immersed}
can be applied.

\begin{figure} [h]
  \includegraphics[scale=0.85]{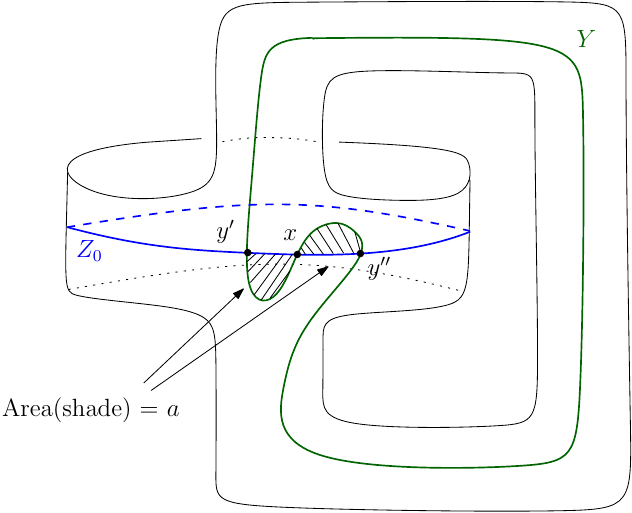}
  \centering
  \caption{A plumbing of two $T^*S^1$ where $Z_0$ and $Y$ are two
    embedded Lagrangians.} \label{f:plumbing}
\end{figure}

Consider the exact symplectic manifold $X$ obtained from the plumbing
of two cotangent bundles of $S^1$ as in Figure~\ref{f:plumbing}, where
$Z_0$ (in blue) is the zero-section of one (horizontal) $T^*S^1$ and
$Y$ (in green) is a Lagrangian in the other (vertical) $T^*S^1$,
perturbed from its zero-section. We orient $Z_0$ in the
counterclockwise orientation and $Y$ in the clockwise orientation
(with respect to how these two Lagrangian are depicted in
Figure~\ref{f:plumbing}), so that the intersection index of $Z_0, Y$
at $y'$ is $+1$. Next we fix a Liouville form $\lambda$ on $X$ for
which both $Z_0$ and $Y$ are exact.

In what follows we will view $Z_0$ and $Y$ as marked Lagrangians and
denote by $\bar{Z}_0$ and $\bar{Y}$ their underlying Lagrangians. We
will shortly be a bit more specific regarding the marking of $Z_0$ and
$Y$, namely the grading and the choice of primitives for the Liouville
form on these Lagrangians. Let $\bar{\mathcal{X}}$ be any collection
of exact Lagrangians in $X$ that contains $\bar{Z}_0$ and $\bar{Y}$
and such that $\bar{\mathcal{X}}$ satisfies the conditions from
\S\ref{subsubsec:filtr-fuk}. For the family $\mathcal{X}$ we take, as
before, the Lagrangians in $\bar{\mathcal{X}}$ together with all their
possible markings (i.e.~we endow each Lagrangian
$\bar{L} \in \bar{\mathcal{X}}$ with all possible primitives of
$\lambda|_{\bar{L}}$ and all possible grading).

Recall from \S\ref{subsec:Fuk-TPC} that in our realization of the
Fukaya category we take Floer data with $0$ Hamiltonians for pairs of
Lagrangians that intersect transversely.  This applies in particular
to $Z_0$ and $Y$, hence the Floer complex $CF_*(Z_0,Y)$ is generated
by the intersection points between $Z_0$ and $Y$, which are
$y', x, y''$.

We endow $Z_0$ and $Y$ with gradings such that:
$$|y'| = |y''| = 0, \text{ and } |x| = 1,$$
where $|\cdot|$ stands for homological grading. It follows that the
Floer complex (with $\mathbb{Z}_2$-coefficients) has the following
form:
\begin{equation} \label{eq:CF-Z0Y}
  \begin{aligned}
    & CF_0(Z_0, Y) = \Z_2 y' \oplus \Z_2 y'', \qquad & &dy'=dy''=0,
    \\
    & CF_1(Z_0, Y) = \Z_2 x, & &dx = y' + y''.
  \end{aligned}
\end{equation}

We now turn to the action filtration on the Floer complex. We choose
primitives for $\lambda|_{\bar{Z}_0}$ and $\lambda|_{\bar{Y}}$ such
that the actions of the intersection points of $Z_0 \cap Y$ are given
by:
\begin{equation} \label{eq:action-Z0Y} \mathcal A(y') = \mathcal
  A(y'') = 0, \quad \mathcal A(x) = a,
\end{equation}
where $a$ is the area depicted in Figure~\ref{f:plumbing}.

It follows that the persistence Floer homology $HF_*(Z_0, Y)$ is
concentrated in degree $0$ and moreover the persistence module
$HF_0(Z_0, Y)$ has the following barcode:
\begin{equation} \label{eq:barcode-HF-Z0Y} \mathcal{B}_{HF_0(Z_0, Y)) }=
  \bigl\{[0,a), [0, \infty)\bigr\}.
\end{equation}
We denote the resulting filtered Fukaya category by
$\fuk(\mathcal{X};\mathscr{P})$ - see \S\ref{subsubsec:filtr-fuk}.

\begin{figure} [h]
\includegraphics[scale=0.85]{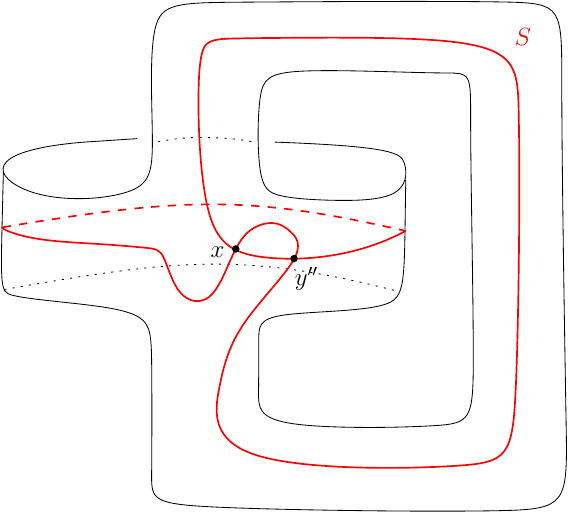}
\centering
\caption{An immersed Lagrangian from a surgery at
  $y'$.}\label{f:immersed}
\end{figure}

Consider a surgery of $Z_0 \cup Y$ at $y'$, so that the outcome
Lagrangian is represented by the immersed Lagrangian $S_{\epsilon}$
(depicted in red) in Figure~\ref{f:immersed}, with two
self-intersection points $x$, $y''$. Here $\epsilon$ is the area of
the surgery handle used in the surgery at $y'$ \zjnote{(see
  \cite{Pol91,LS91,Bi-Co:cob1} for a detailed description of
  surgery.)}  Note that $S_{\epsilon}$ is an exact immersed Lagrangian
(this is always the case when we perform surgery of two exact
Lagrangians at a single point). Since the choices of the primitives of
$\lambda|_{\bar{Z}_0}$ and $\lambda|_{Y}$ coincide at $y'$ (recall
that $\mathcal{A}(y')=0$) we can endow $S_{\epsilon}$ with a primitive
induced by those of $Z_0$ and $Y$ (more precisely, one needs to
parametrize $S_{\epsilon}$ by an immersion
$j: \hat{S} \longrightarrow X$ with $j(\hat{S})=S$ such that
$j^*\lambda$ is exact and define the primitive of the latter $1$-form
on $\hat{S}$\zjnote{)}. We also perturb slightly this $S_{\epsilon}$
by \zjred{a small} Hamiltonian perturbation of Hofer norm at most
$\epsilon'$ thus getting $S=S_{\epsilon,\epsilon'}$ that intersects
transversely all the elements in $\bar{\mathcal{X}}$.

It can easily be seen that $S$ is unobstructed (in the sense that it
carries no pseudo-holomorphic teardrops (see also
\cite{Bi-Co:LagPict})). Therefore, $S$ gives rise to a filtered
\zjnote{$A_{\infty}$-module} over $\fuk(\mathcal{X})$, which we denote
by $\mathscr{S}$. The existence of this module requires in fact
elements of the \zjnote{theory of} immersed Fukaya \zjnote{category}
mentioned in Remark \ref{rem:imme-no} (ii), however only the most
elementary part of this \zjnote{theory} is needed here. In essence,
because $S$ is unobstructed we can include $S$ in a larger family
$\mathcal{X}'$ that contains $\mathcal{X}$ and $S$. We may still
define a Fukaya filtered category associated to $\mathcal{X}'$,
$\fuk(\mathcal{X}';\mathscr{P}')$.  The immersed Lagrangian $S$ now
gives a (filtered) Yoneda module $\mathcal{Y}(S)$ over
$\fuk(\mathcal{X}';\mathscr{P}')$. The perturbation data
$\mathscr{P}'$ can be picked to be admissible also for the set
$\mathcal{X}$ and there is an embedding
$\fuk(\mathcal{X};\mathscr{P}')\to \fuk(\mathcal{X}';\mathscr{P}')$
\zjnote{such} that $\mathcal{Y}(S)$ pulls back to a module over
$\fuk(\mathcal{X};\mathscr{P}')$. Finally, using that
$\fuk(\mathcal{X})$ is \zjnote{independent of the perturbation data
  up} to filtered quasi-equivalence, we can pull-back this module to
the module $\mathscr{S}$ over $\fuk(\mathcal{X};\mathscr{P})$.  The
details of this process appear in \zjnote{\cite[Chapter 3, Theorem
  3.12]{BCZ23}}, but only in the case when $S$ is embedded. The
modifications needed for our $S$ here are easy to implement as in
\cite{Bi-Co:LagPict} (where is covered the construction of
non-filtered, Fukaya type categories with immersed objects in much
larger generality), and noting that the construction there of the
Floer complexes of the type $CF(S,S)$, $CF(S,X)$, $CF(X,S)$ (for
$X\in \mathcal{X}$) produces filtered structures, and similarly for
the higher $A_{\infty}$-operations.

Next, consider the twisted complex (over $\fuk(\mathcal{X})$) defined
as the mapping cone
\begin{equation} \label{eq:cone-Z0Y}
  \mathscr{C} := \text{Cone} (Z_0 \xrightarrow{\; y' \;} Y)
\end{equation}
of the morphism
$y' \in \hom_{\fuk(\mathcal{X})}(Z_0,Y) = CF_0(Z_0,Y)$. (Recall that
$y'$ is a cycle). Since $y'$ has filtration level $0$ and both $Z_0$,
$Y$ are filtered twisted complexes, $\mathscr{C}$ can be endowed with
the structure of a filtered twisted complex
$\in \C'\fuk(\mathcal{X})\subset \C\fuk(\mathcal{X})$ in the standard
way.

We \zjnote{claim} that $\mathscr{S}\in \C\fuk(\mathcal{X})$ and that
there exists an $r$-isomorphism $$\phi: \mathscr{C}\to \mathscr{S}$$
such that $r=r(\epsilon, \epsilon')$ tends to $0$ \zjnote{when} both
$\epsilon$ and $\epsilon'$ go to $0$.

To prove this we compare $\mathcal{Y}(S)$ and $\mathscr{C}$ over the
category $\fuk(\mathcal{X}';\mathcal{P}')$, and we notice the
existence of an $r$-isomorphism
$$\bar{\phi}: \mathscr{C}\to \mathcal{Y}(S)$$ 
of the form $\mu_{2}(\alpha, - )$ for an appropriate intersection
point $\alpha$ between $S$ and \zjnote{$Z_{0}\cup Y$}. To see this we
can interpret $\mathscr{C}=\mathcal{Y}(S')$ where $S'$ is a marked
immersed Lagrangian, \pbred{in the sense discussed
  in~\S\ref{subsec:immersed},} with obvious primitive and grading, and
with a marked point $I_{\mathscr{C}}= \{y'\}$. At this point it is
also useful to use the $\epsilon'$-small Hamiltonian perturbation to
be a bit special, in essence associated to a Morse function with two
critical points on $\hat{S}$ and in this case $\alpha$ is the
intersection point between $S_{\epsilon,\epsilon'}$ and $S'$ that
corresponds to the maximum of this function. By pulling back
$\bar{\phi}$ over $\fuk(\mathcal{X})$ we get $\phi$.

We now pass to calculations in \zjnote{the $K$-group}
$K(\C\fuk(\mathcal{X}))$. We have:
\begin{equation} \label{eq:K-Z0Y} [\mathscr{C}] = - [Z_0] + [Y] \in
  K(\C\fuk(\mathcal{X})).
\end{equation}

Next we compute the pairing $\kappa([\mathscr{C}], \mathscr{C}])$. By
linearity, and by Theorem~\ref{thm-emb-char} we obtain:
\begin{equation} \label{eq:kappa-CC}
  \begin{aligned}
    \kappa([\mathscr{C}], [\mathscr{C}]) & = \kappa([Y], [Y]) +
    \kappa([Z_0], [Z_0]) - \kappa( [Y], [Z_0])) -
    \kappa([Z_0], [Y]) \\
    & = -\chi(Y) -\chi(Z_0) -\lambda_{HF_*(Z_0, Y)}(t) +
    \lambda_{HF_*(Z_0, Y)}(t^{-1}) \\
    & = 0 + 0 - (1-t^a + 1) + (1 - t^{-a}+1) \\
    & = -t^{-a} + t^a,
  \end{aligned}
\end{equation}
where the third equality comes from the barcode
data~\eqref{eq:barcode-HF-Z0Y} and the identities in
Theorem~\ref{thm-emb-char}.

Since $\kappa([\mathscr{C}], [\mathscr{C}])$ is not a constant, by
Corollary~\ref{thm-emb-char}, the class
$[\mathscr{S}] = [\mathscr{C}] \in K(\C\fuk(\mathcal{X}))$ cannot be
represented by any Lagrangian in $\mathcal{X}$ .

Further, we now want apply Theorem \ref{thm:dist-immersed}.  First ,
the data before provide a reduced resolution $D$ of weight $r$ of
$\mathscr{S}$ with linearization $\ell(D)=(T^{-1}Y, Z_{0})$.

Clearly, $[D]=-[\mathscr{C}]$ so
$\mathrm{gap}(\kappa([D],[D]))=a$. Obviously $n_{D}=2$. It is also
immediate to calculate \zjnote{$b_{\max}^{\#}(\{T^{-1}Y,
  Z_{0}\})=2$}. We fix $\mathcal{F}=\{T^{-1}Y, Z_{0}\}$. Thus, we
obtain
$$\mathscr{R}^{\mathcal{F}}(\mathscr{S})\leq  \frac{8 r}{a}~.~$$
Theorem \ref{thm:dist-immersed} implies that if
\zjnote{$r= r(\epsilon,\epsilon')\leq \frac{a}{16}$}, then
$\mathscr{C}$ can not represent an embedded Lagrangian (in the sense
that it is not the Yoneda module of an element in
$\mathcal{X}$). Moreover, this result is independent of the set
$\mathcal{X}$ as long as $\mathcal{X}$ contains $Z_{0}, Y$.

Another way to apply Theorem \ref{thm:dist-immersed} is to deduce by
an argument very similar to the above that, assuming that both
$\epsilon$ and $\epsilon'$ are very small, then the Hofer norm of a
Hamiltonian homotopy needed to deform $S$ into an embedding is at
least \zjnote{$\frac{a}{16}$}.

\begin{rem}As noted in relation to the constant $4$ in Theorem
  \ref{thm:dist-immersed}, the constant $16$ here is \zjnote{not
    optimal}. In essence, the argument we discussed is very general
  but it produces constants that can be easily improved in our
  example. In this example it is possible to use a Hamiltonian
  homotopy of energy $a$ to deform $S$ into an embedding and it is
  likely that this is sharp.
\end{rem}



\bibliography{biblio_tpc}

%


\end{document}